\documentclass[final,1p,times,twocolumn]{elsarticle}

\usepackage{natbib}

\usepackage{graphicx}
\usepackage{epsfig}

\usepackage{subfigure}
\usepackage{color}

\usepackage{latexsym,bm}

\usepackage{amssymb}
\usepackage{amsmath}

\usepackage[linesnumbered,ruled]{algorithm2e}
\usepackage{algorithmic}

\usepackage{cases}

\usepackage{threeparttable}
\usepackage{multirow}

\usepackage{epstopdf}

\newtheorem{Def}{Definition}
\newtheorem{thm}{Theorem}
\newtheorem{lem}{Lemma}
\newtheorem{cor}{Corollary}
\newtheorem{rem}{Remark}

\newcommand{\norm}[1]{\left\Vert#1\right\Vert}
\newcommand{\abs}[1]{\left| #1 \right|}
\newcommand{\pref}[1]{(\ref{#1})}

\allowdisplaybreaks

\begin{document}

\begin{frontmatter}

\title{Isogeometric Least-squares Collocation Method with Consistency and Convergence Analysis}


\author[linadd]{Hongwei Lin\corref{cor1}}
    \cortext[cor1]{Corresponding author: phone number: 86-571-87951860-8304, fax number: 86-571-88206681, email:
    hwlin@zju.edu.cn}
\author[linadd]{Yunyang Xiong}
\author[linadd]{Xiao Wang}
\author[huadd]{Qianqian Hu}
\address[linadd]{Department of Mathematics, State Key Lab. of CAD\&CG, Zhejiang University, Hangzhou, 310027, China}
\address[huadd]{Department of Mathematics, Zhejiang Gongshang
                    University, Hangzhou, 310018, China}
\date{}



\begin{abstract}
In this paper, we present the isogeometric least-squares collocation (IGA-L) method, which determines the numerical solution by making the approximate differential operator fit the real differential operator in a least-squares sense. The number of collocation points employed in IGA-L can be larger than that of the unknowns. Theoretical analysis and numerical examples presented in this paper show the superiority of IGA-L over state-of-the-art collocation methods. First, a small increase in the number of collocation points in IGA-L leads to a large improvement in the accuracy of its numerical solution. Second, IGA-L method is more flexible and more stable, because the number of collocation points in IGA-L is variable. Third, IGA-L is convergent in some cases of singular parameterization. Moreover, the consistency and convergence analysis are also developed in this paper.
\end{abstract}


\end{frontmatter}

\section{Introduction}

 While classical Finite Element Analysis (FEA) methods,
    widely employed in physical simulations,
    are based on element-wise piecewise polynomials,
    computer-aided design (CAD) models
    are usually represented by non-uniform rational basis splines (NURBSs) with non-linear
    NURBS basis functions.
 Therefore, the first task in a CAD model simulation is to
    transform the non-linear NURBS-represented CAD model into a linear
    mesh representation.
 This mesh transformation is a very tedious
    operation, and it has become the most time-consuming task of the whole FEA
    procedure.
 To avoid the mesh transformation
    and advance the seamless integration of CAD and computer-aided engineering (CAE),
    isogeometric analysis (IGA) was invented by Hughes et al.~\cite{hughes2005isogeometric}.
 IGA is based on the NURBS basis functions of degree higher than $1$,
    hence it can deal with NURBS-represented CAD models directly.
 In this way, IGA not only saves a significant amount of computation,
    it also greatly improves the numerical accuracy of the solution.

 In IGA, the analytical solution $T$ to a boundary value problem
    is approximated by a NURBS function $T_r$ with unknown coefficients.
(For brevity, we only mention the boundary value problem in this
    paper.
    However, this method is also suitable for the initial value problem.)
 Solving the boundary value problem is then equivalent to
    determining the unknown coefficients of $T_r$.
 If the order of $T_r$ is higher
    than that of the differential operator $\mathcal{D}$ of the
    boundary value problem,
    $\mathcal{D}T_r$ can be represented explicitly by a NURBS derivative
    formula.
 Therefore, collocation methods can be applied to the strong form
    of the boundary value problem to determine the unknown coefficients of $T_r$.

 In Ref.~\cite{auricchio2010isogeometric},
    an isogeometric collocation (IGA-C) method was presented.
 Suppose $n$ is the number of unknown coefficients of $T_r$.
    IGA-C first samples $n$ values $\mathcal{D}T(\eta_j),\ j
    = 1,2,\cdots,n$,
    and then generates a linear system of equations
    by making $\mathcal{D}T_r$ interpolate the $n$
    values, i.e.,
    $\mathcal{D}T_r(\eta_j) = \mathcal{D}T(\eta_j),\
    j=1,2,\cdots,n$.
 In this way, the unknown coefficients of $T_r$ can be determined
    by solving the linear system.

 Essentially, IGA-C acquires the unknown coefficients by
    interpolation,
    where the number of the collocation points must
    be equal to that of the unknown coefficients.
 In this paper,
    we propose the isogeometric least-squares collocation (IGA-L)
    method,
    which allows the number of collocation points to be larger
    than that of the unknown coefficients.
 It yields some advantages over IGA-C.
 Instead of interpolation,
    IGA-L makes use of approximation to determine the
    unknown coefficients of $T_r$.
 Specifically,
    IGA-L first samples $m$ values $\mathcal{D}T(\eta_j), j =
    1,2,\cdots, m$,
    where $m$ is greater than the number of unknown coefficients, i.e., $m > n$.
 The coefficients of the unknown solution $T_r$ are then determined by
    solving the least-squares fitting problem
    \begin{equation*}
        \min{\sum_{j=1}^m \norm{\mathcal{D}T(\eta_j) - \mathcal{D}
            T_r(\eta_j)}^2}.
    \end{equation*}

 There are two advantages of IGA-L over IGA-C.
 First, numerical tests presented in this paper show that
    a small increase in the computation of IGA-L leads to
    a large improvement in the numerical accuracy of the solution,
    even though the computational cost of IGA-L is only slightly more
    than that of IGA-C.
 Second, IGA-L is more flexible and more stable than IGA-C,
    because IGA-L allows a variable number of collocation points,
    while the number of collocation points in IGA-C is fixed to be equal
    to the number of control points.

 The structure of this paper is
    as follows.
 In Section~\ref{sec:related_work},
    some related work is briefly reviewed.
 In Section~\ref{sec:generic_formulation},
    the generic formulation of IGA-L is presented.
 In Section~\ref{sec:analysis},
    we show the consistency and convergence properties of IGA-L.
 After thoroughly comparing IGA-L and
    IGA-C, both in theory and with numerical examples in Section~\ref{sec:comparison},
    this paper is concluded in Section~\ref{sec:conclusion}.

\subsection{Related Work}
\label{sec:related_work}

 \textbf{Least-squares collocation methods:}
 Although collocation-based meshless methods are efficient,
    equilibrium conditions are satisfied only at collocation points.
 Thus, collocation-based meshless methods can result in significant
    error.
 To improve computational accuracy,
    Zhang et al. developed a least-squares collocation
    method~\cite{zhang2001least},
    where equilibrium conditions hold at both collocation points and
    auxiliary points in a least-squares sense.
 In order to generate a better conditioned linear algebraic
    equations,
    a least-squares meshfree collocation method was proposed,
    based on first-order differential
    equations~\cite{jiang2012least}.
 In Ref.~\cite{kim2003point},
    a point collocation method was invented that
    employs the approximating derivatives based on the moving
    least-square reproducing kernel approximations.
 Moreover, several schemes using least-squares collocation methods were
    developed for two- and three-dimensional heat conduction
    problems~\cite{dai2011weighted},
    transient and steady-state hyperbolic
    problems~\cite{afshar2009collocated},
    and adaptive analysis problems in linear elasticity~\cite{kee2008least}.

 On the other hand, there are some meshfree methods which handle the
    improvement of stability.
 By eliminating the rank deficiency with stress points,
    a meshfree particle method was developed for large deformation,
    nonlinear problems that employs a Lagrangian kernel with correction of the derivatives~\cite{rabczuk2004stable}.
 In~\cite{rabczuk2004cracking}, a simplified meshfree method for arbitrary
    evolving cracks was proposed,
    where the crack is modelled by a discontinuous enrichment that can be arbitrarily aligned in the body at each particle.
 Moreover, an approach for modelling discrete cracks in meshfree particle
    methods in three dimension was devised,
    where the growth of a crack is represented by activation of crack surfaces with arbitrary orientation~\cite{rabczuk2007three}.

 However, the consistency and convergence properties of the
    aforementioned least-squares collocation methods were only validated
    using numerical examples,
    and theoretical numerical analyses were not reported.
 In this paper, we not only develop an isogeometric least-squares
    collocation method,
    but also show its consistency and convergence properties in
    theory.

 \textbf{Isogeometric analysis:}
 NURBS is a mathematical model
    for representing curves and surfaces
    by blending weighted control points and NURBS basis functions.
 Hence, the shape of the NURBS curves and surfaces can be easily modified
    by moving their control points.
 Because of the desirable traits of NURBS basis functions,
    NURBS curves and surfaces have many good properties,
    such as convex hull, affine invariance, and variation diminishing.
 Moreover, NURBSs can represent conic sections and quadric surfaces
    accurately.
 Therefore, NURBSs have been widely used in CAD,
    computer-aided manufacturing (CAM), and CAE,
    and have become a part of numerous industry standards,
    such as IGES, STEP, ACIS, and PHIGS.
 For more details on NURBSs,
    there are excellent books on the subject~\cite{piegl1997nurbs,
    de2001practical}.

 While a NURBS employs non-linear basis functions,
    classical FEA is based on element-wise piecewise polynomials.
 Hence, when analyzing NURBS-represented CAD models using classical
    FEA,
    the CAD model should be discretized into a mesh model.
 The mesh generation procedure not only yields an approximation,
    it is also very tedious, and hence has become a bottleneck in FEA.
 To overcome these difficulties,
    Hughes et al. invented the IGA
    technique~\cite{hughes2005isogeometric}.
 Because IGA is based on NURBS basis functions,
    it can handle NURBS-represented CAD models directly without
    generating a mesh.
 Moreover, because NURBSs can represent complex shapes (and physical
    fields) with significantly fewer control points than a mesh representation,
    the computational efficiency and numerical accuracy of IGA are higher than
    the classical FEA method.
 In addition, because of the knot insertion property of NURBSs,
    the original shape of the CAD model can be maintained exactly in the refinement
    procedure~\cite{hughes2005isogeometric}.

 In geometric design,
    the NURBS representation is usually employed to model curves and surfaces.
 There are a few effective methods in geometric design for modeling spline solids.
 To model a NURBS solid for IGA applications,
    Zhang et al. proposed a solid construction method from a boundary
    representation~\cite{zhang2012solid}.
 In the study of IGA,
    making the geometric representation
    more suitable for analysis purposes is a key research goal.
 Cohen et al. presented the analysis-aware modeling
    technique~\cite{cohen2010analysis}.
 Moreover, T-spline~\cite{bazilevs2010isogeometric,
    dorfel2010adaptive},
    trimmed surfaces~\cite{kim2009isogeometric},
    subdivision solids~\cite{burkhart2010iso},
    and splines on triangulations~\cite{speleers2012isogeometric,jaxon2014isogeometric}
    were also employed in the IGA method to model the physical domain.

 Currently, the IGA method has been successfully applied in various
    simulation
    problems, including elasticity~\cite{auricchio2007fully,
    elguedj2008view},
    structure~\cite{cottrell2006isogeometric, hughes2008duality,
    wall2008isogeometric},
    and fluids~\cite{bazilevs2008isogeometric, bazilevs2006fluid,
    bazilevs2009patient}.
 On the other hand,
    some work concerns the computational aspect of the IGA method,
    for instance,
    to improve the accuracy and efficiency by reparameterization and
    refinement~\cite{bazilevs2006isogeometric, cottrell2007studies,
    hughes2010efficient,aigner2009swept, xu2011optimal}.
 Moreover, fast solvers for both Galerkin and collocation approaches were
    developed in Refs.~\cite{donatelli2015robusta,donatelli2015robustb}.

 \textbf{Isogeometric collocation methods:}
 The collocation method is a simple numerical method for solving
    differential equations
    that can generate a solution that satisfies the differential
    equation at a set of discrete points,
    called collocation points~\cite{cottrell2009isogeometric}.
 If the order of the unknown NURBS function
    that is employed to approximate the solution of the differential
    equation
    is high enough,
    the collocation method can be applied to the strong form of the
    differential equation.
 In this way, the IGA-C method was
    presented in Ref.~\cite{auricchio2010isogeometric}.
 However, IGA-C is greatly influenced by the
    collocation points.
 Recently, a comprehensive study on IGA-C discovered its superior behavior
    over the Galerkin method in terms of
    its accuracy-to-computational-time ratio~\cite{schillinger2013isogeometric},
    and the consistency and convergence properties of the IGA-C
    method were developed in Ref.~\cite{lin2013consistency}.
 In Ref.~\cite{anitescu2015isogeometric},
    the isogeometric superconvergent collocation method (IGA-SC) was developed,
    where locations of the collocation points were derived from the superconvergence theory.
 Furthermore, an optimally convergent isogeometric collocation scheme
    for odd degrees was proposed in~\cite{montardini2016optimal},
    which are a subset of the Galerkin superconvergent points.
 It was shown that there exist the collocation points,
    called Cauchy-Galerkin points,
    that produce the Galerkin solution exactly~\cite{gomez2016variational}.
 In~\cite{casquero2016isogeometric},
    analysis-suitable T-splines were employed in combination with isogeometric collocation methods to solve second- and fourth-order boundary-value problems.

 Moreover, based on the local hierarchical refinement of NURBSs,
    adaptive IGA-Cs were developed and
    analyzed~\cite{schillinger2013isogeometric}.
 Meanwhile, IGA-C has also been extended to multi-patch NURBS configurations,
    various boundary and patch
    interface conditions, and explicit dynamic analysis~\cite{auricchio2012isogeometric}.
 Recently, IGA-C was successfully employed to solve the Timoshenko
    beam problem~\cite{beirao2012avoiding} and
    spatial Timoshenko rod problem~\cite{auricchio2013locking},
    showing that mixed collocation schemes are locking-free, independently
    of the choice of the polynomial degrees for the unknown fields.
 It was shown that IGA-C is particularly suitable for solving the system of
    ODEs governing the non-prismatic beam problem~\cite{Balduzzi2017}.
 Moreover, IGA-C was proposed for the linear static bending analysis of
    laminated composite plates governed by Reissner-Mindlin theory~\cite{Pavan2017}.


\section{Generic Formulation of IGA-L}
\label{sec:generic_formulation}

 Consider the following boundary value problem,
    \begin{equation} \label{eq:original_pde}
    \begin{cases}
    \mathcal{D}T = f,\ \text{in}\ \Omega \in \mathbb{R}^d,\\
    \mathcal{G}T = g,\ \text{on}\ \partial \Omega, \\
    \end{cases}
    \end{equation}
    where $\Omega$ is the physical domain in $\mathbb{R}^d$,
    $\mathcal{D}$ is a differential operator in the physical domain,
    $\mathcal{G}T = g$
    is the boundary condition,
    and $f : \Omega \rightarrow \mathbb{R}$ and
    $g : \partial \Omega \rightarrow \mathbb{R} $
    are given functions.
 Suppose $d_1$ is the maximum order of derivatives appearing in
    $\mathcal{D}: V \rightarrow W$,
    where $V$ and $W$ are two Hilbert spaces,
    and the analytical solution $T \in C^{d_2}(\Omega), d_2 > d_1 \geq 1$.

 In IGA,
    the physical domain $\Omega$ is represented by a NURBS mapping:
    \begin{equation}\label{eq:b_spline_mapping}
        \bm{F}: \Omega_0 \rightarrow \Omega,
    \end{equation}
    where $\Omega_0 \in \mathbb{R}^d $ is the parametric domain.
 By replacing the control points of $\bm{F}(\Omega_0)$ with
    unknown control coefficients,
    the representation of the unknown numerical solution $T_r$ is generated.

 Suppose there are $n$ unknown control coefficients in the unknown numerical solution $T_r$.
 We first sample $m_1$ points $\bm{\theta}_k$
    inside the parametric domain $\Omega_0$
    that correspond to $m_1$ values
    $\bm{\eta}_k = \bm{F}(\bm{\theta}_k),\ k = 1,2,\cdots,m_1,$
    inside the physical domain $\Omega$.
 Furthermore, we sample $m_2$ points $\bm{\theta}_l$
    on the parametric domain boundary $\partial \Omega_0$
    that correspond to $m_2$ values
    $\bm{\eta}_l = \bm{F}(\bm{\theta}_l),\ l = m_1+1,m_1+2,\cdots,m_1+m_2,$ on the physical domain boundary $\partial \Omega$.
 The total number of these points, i.e., $m = m_1 + m_2$,
    is greater than the number of unknown coefficients of the numerical solution $T_r$,
    namely, $m = m_1 + m_2 > n$.
 Just as in IGA-C,
    these sampling points are also called \textbf{collocation points}.

 Substituting these collocation points into the boundary value
    problem~\pref{eq:original_pde},
    a system of equations with $m$ equations and $n$
    unknowns is obtained (where $m = m_1 + m_2 > n$),
    \begin{equation} \label{eq:system_of_equations}
    \begin{cases}
    \mathcal{D}T_r(\bm{\eta}_k) = f(\bm{\eta}_k),\ k = 1,2,\cdots m_1,\\
    \mathcal{G}T_r(\bm{\eta}_l) = g(\bm{\eta}_l),\ l = m_1+1, m_1+2, \cdots, m_1+m_2. \\
    \end{cases}
    \end{equation}

 Arranging the unknowns of the numerical solution $T_r$ into an
    $n \times 1$ matrix, i.e., $X = [x_1\ x_2\ \cdots\ x_n]^T$,
    the system of equations~\pref{eq:system_of_equations}
    can be represented in matrix form by
    \begin{equation*}
        A X = b.
    \end{equation*}
 Because the number of equations is greater than the number of
    unknowns,
    the solution is sought in the least-squares sense, i.e.,
    \begin{equation} \label{eq:least_squares}
        \min_X \norm{AX - b}^2.
    \end{equation}

 \textbf{Computation of the least-squares problem~\pref{eq:least_squares}:}
 The least-squares problem~\pref{eq:least_squares} is very important in practice,
    and there have been developed lots of robust and efficient methods for solving it~\cite{golub1996matrix}.
 One frequently employed method is to solve the normal
    equation~\pref{eq:normal_equ},
    \begin{equation}\label{eq:normal_equ}
        A^T A X = A^T b.
    \end{equation}
 Although the condition number of the matrix $A^TA$ is the square of that of the matrix $A$,
    $A^TA$ is a symmetric positive definite matrix,
    and the normal equation~\pref{eq:normal_equ} can be solved
    efficiently by Cholesky decomposition~\cite{golub1996matrix}.
 Moreover, Householder orthogonalization and Given
    orthogonalization are also two efficient methods~\cite{golub1996matrix} which often employed
    in solving the least-squares problem~\pref{eq:least_squares}.
 For more methods on solving~\pref{eq:least_squares},
    please refer to~\cite{golub1996matrix}.

 \begin{rem} \label{rem:assumption}
    We assume that the matrix $A$ is of full rank,
        and then $A^T A$ is non-singular,
        and the linear system~\pref{eq:normal_equ} has unique solution.
 \end{rem}

 In the following, we will show the consistency and convergence
    properties of IGA-L, and compare it with IGA-C and IGA-SC~\cite{anitescu2015isogeometric},
    both in theory and with numerical examples.

\section{Numerical Analysis}
\label{sec:analysis}

 In the IGA-L method developed above,
    a NURBS function $T_r$ is employed to approximate the analytical solution $T$ of the
    boundary value problem~\pref{eq:original_pde},
    and hence the real differential operator $\mathcal{D}T$ is approximated by
    $\mathcal{D}T_r$.
 In Ref.~\cite{lin2013consistency},
    it was shown that both $\mathcal{D}T_r$ and $T_r$ are defined
    on the same knot intervals, i.e.,

 \begin{lem}\label{lem:same_knot_intervals}
 If $T_r$ is a NURBS function
    and $\mathcal{D}$ is a differential operator,
    then both $\mathcal{D}T_r$ and $T_r$ are defined
    on the same knot intervals~\cite{lin2013consistency}.
 \end{lem}

 Moreover, given a set $\Phi$, the \textbf{diameter} of $\Phi$ is
    defined as
    \begin{equation*}
        diam(\Phi) = \max\{d(x,y), x, y \in \Phi\},
    \end{equation*}
    where $d(x,y)$ is the Euclidean distance between $x$ and $y$.
 In the following, we suppose the NURBS function $T_r$ to be defined on a knot grid
    $\mathcal{T}^h$,
    where $h$ is the \textbf{knot grid size}:
 In the one-dimensional case, $\mathcal{T}^h$ is a knot sequence,
    and $h = \max_i\{diam([u_i,u_{i+1}])\}$;
 in the two-dimensional case, $\mathcal{T}^h$ is a rectangular grid,
    and $h = \max_{ij}\{diam([u_i,u_{i+1}] \times [v_j,v_{j+1}])\}$;
 in the three-dimensional case, $\mathcal{T}^h$ is a hexahedral grid,
    and $h = \max_{ijk}\{diam([u_i,u_{i+1}] \times [v_j,v_{j+1}] \times
    [w_k,w_{k+1}])\}$.

 In this section, we study the consistency and convergence
    properties of IGA-L,
    i.e., when $h \rightarrow 0$,
    not only will approximate differential operator $\mathcal{D}T_r$ tend to
    real operator $\mathcal{D}T$,
    but numerical solution $T_r$ will also tend to analytical
    solution $T$.

 \subsection{Consistency}
 \label{sec:consistency}

 In this section, we explore the consistency property of the
    IGA-L method.
 In Section~\ref{ssub:consistency_special},
    the consistency of the IGA-L method with special operators which have polynomial coefficients is developed.
 In Section~\ref{ssub:consistency_generic},
    the consistency of the IGA-L method in the generic case is studied.

 \subsubsection{Consistency with special operators}
 \label{ssub:consistency_special}

 In this section, we deal with special operators $\mathcal{D}$ and
    $\mathcal{G}$,
    which have polynomial coefficients.

 Suppose the NURBS function
    $T_r(\bm{\eta}),\ \bm{\eta} \in \Omega_0 \subset \mathbb{R}^d$
    defined on the knot grid $\mathcal{T}^h$,
    has $n$ unknown control coefficients $p_{\bm{i}}$, i.e.,
    \begin{equation}\label{eq:rewrite-nurbs}
    T_r(\bm{\eta}) = \sum_{\bm{i}} p_{\bm{i}} \frac{w_{\bm{i}} B_{\bm{i}}(\bm{\eta})}{W(\bm{\eta})}
                   = \frac{P(\bm{\eta})}{W(\bm{\eta})},
                    \quad \bm{\eta}=(\eta_1,\eta_2,\cdots,\eta_d) \in \Omega_p \subset \mathbb{R}^d,
 \end{equation}
    where the subscript $\bm{i}$ is an index vector,
    $\bm{i} = (i_1, i_2, \cdots, i_d)$,
    $w_{\bm{i}} > 0$ are known weights,
 \begin{equation} \label{eq:tensor_basis}
    B_{\bm{i}}(\bm{\eta})=B_{i_1}(\eta_1) B_{i_2}(\eta_2) \cdots B_{i_d}(\eta_d),
 \end{equation}
    are the B-spline basis functions.
 Moreover, the weight function $W(\bm{\eta})$ is a known polynomial
    spline function,
    and $P(\bm{\eta})$ is a polynomial spline function with
    $n$ unknown control coefficients $p_{\bm{i}}$.

 According to the result developed in Ref.~\cite{lin2013consistency},
    $\mathcal{D}T_r$ can be represented by
    \begin{equation} \label{eq:dtr_representation}
        \mathcal{D}T_r(\bm{\eta})
            = \sum_{\bm{i}} p_{\bm{i}} \mathcal{D} \frac{w_{\bm{i}}B_{\bm{i}}(\bm{\eta})}{W(\bm{\eta})}
            = \sum_{\bm{i}} p_{\bm{i}} \frac{\bar{B}_{\bm{i}}(\bm{\eta})}{\bar{W}(\bm{\eta})}
            = \frac{\bar{P}(\bm{\eta})}{\bar{W}(\bm{\eta})},
    \end{equation}
 where $\bar{B}_i(\bm{\eta})$, the result by applying the
    differential operator $\mathcal{D}$ to
    $\frac{w_{\bm{i}}B_{\bm{i}}(\bm{\eta})}{W(\bm{\eta})}$,
    is a polynomial spline function,
    $\bar{W}(\bm{\eta}) \neq 0$ is the power of $W(\bm{\eta})$,
    and
    \begin{equation} \label{eq:poly_DT}
        \bar{P}(\bm{\eta}) = \sum_{\bm{i}} p_{\bm{i}} \bar{B}_i(\bm{\eta}),
    \end{equation}
    is a polynomial B-spline function with $n$ unknowns $p_{\bm{i}}$.

 Similarly, $\mathcal{G}T_r (\bm{\eta})$ in
    Eq.~\pref{eq:original_pde} can be written as
    \begin{equation} \label{eq:gtr_representation}
        \mathcal{G}T_r (\bm{\eta})
        = \sum_{\bm{i}} p_{\bm{i}} \mathcal{G}
            \frac{w_{\bm{i}}B_{\bm{i}}(\bm{\eta})}{W(\bm{\eta})}
        = \sum_{\bm{i}} p_{\bm{i}} \frac{\tilde{B}_{\bm{i}}(\bm{\eta})}{\tilde{W}(\bm{\eta})}
        = \frac{\tilde{P}(\bm{\eta})}{\tilde{W}(\bm{\eta})},
    \end{equation}
    where $\tilde{B}_{\bm{i}}(\bm{\eta})$,
    the result generated by applying the operator $\mathcal{G}$ to
    $\frac{w_{\bm{i}}B_{\bm{i}}(\bm{\eta})}{W(\bm{\eta})}$,
    is a polynomial spline function,
    $\tilde{W}(\bm{\eta}) \neq 0$ is a known B-spline function,
    and
    \begin{equation} \label{eq:poly_GT}
    \tilde{P}(\bm{\eta}) = \sum_{\bm{i}} p_{\bm{i}} \tilde{B}_{\bm{i}}(\bm{\eta}),
    \end{equation}
    is an unknown B-spline function with $n$ unknowns $p_{\bm{i}}$.

 By the result developed in Ref.~\cite{lin2013consistency}, $\bar{P}(\bm{\eta})$~\pref{eq:poly_DT} and
    $\tilde{P}(\bm{\eta})$~\pref{eq:poly_GT} both have the same
    break point sequence and the same knot intervals as $T_r(\bm{\eta})$.
 So they can be made to be defined on the same knot sequence by
    knot insertion and degree elevation.

 Therefore, based on Eqs.~\pref{eq:poly_DT} and~\pref{eq:poly_GT},
    the linear system~\pref{eq:system_of_equations} becomes
  \begin{equation}\label{eq:equ_linear_system}
        \begin{cases}
        & \bar{P}(\bm{\eta}_k)
            = \sum_{\bm{i}} p_{\bm{i}}\bar{B}_{\bm{i}}(\bm{\eta}_k)
            = \bar{W}(\bm{\eta}_k)f(\bm{\eta}_k),\ k = 0,1,\cdots,m_1, \\
        & \tilde{P}(\bm{\eta}_l)
            = \sum_{\bm{i}} p_{\bm{i}}\tilde{B}_{\bm{i}}(\bm{\eta}_l)
            = \tilde{W}(\bm{\eta}_l)
            g(\bm{\eta}_l),\ l = m_1+1, m_1+2, \cdots, m_1+m_2=m.
        \end{cases}
    \end{equation}
  By Remark~\ref{rem:assumption}, the coefficient matrix of the linear
    system~\pref{eq:equ_linear_system} is of full rank.
  As a consequence, the polynomial spline functions in
    $\Phi = \{\bar{B}_{\bm{i}}(\bm{\eta}); \tilde{B}_{\bm{j}}(\bm{\eta})\}$
    are linear independent.
  Otherwise, the coefficient matrix of system ~\pref{eq:equ_linear_system}
    is not of full rank.
  Because each of polynomial spline functions in $\Phi$ is a linear
    combination of B-spline basis functions,
    the spline space generated by the combination of functions in $\Phi$ is a B-spline sub-space $\mathbb{S}$,
    defined on the knot grid $\mathcal{T}^h$.
  Therefore, the least-squares solution to the linear
    system~\pref{eq:equ_linear_system} is actually the least-squares projection to the B-spline sub-space $\mathbb{S}$.
  We denote the least-squares projector as $\mathcal{P}$.
  Thus, the following lemma is reached.
  \begin{lem} \label{lem:ls_projection}
    The IGA-L solution to the boundary problem~\pref{eq:original_pde} is the least-squares projection to a B-spline sub-space $\mathbb{S}$.
  \end{lem}

  In Ref.~\cite{shadrin2001thel}, Shadrin proved the
    ``de Boor's conjecture'', i.e.,
    the $L_\infty$-norm of the $L_2$-spline projector is bounded independently of the knot sequence in univariate case.
  Moreover, Passenbrunner et. al. extended this result to tensor product
    spline projections~\cite{passenbrunner2013almost}.
  That is,
  \begin{lem} \label{lem:ls_projector_bound}
    Let $\mathcal{\bar{P}}: C(\Omega) \rightarrow \mathbb{\bar{S}}$ is
        a $L_2$ projector from the continuous function space $C(\Omega)$ to a B-spline space $\mathbb{\bar{S}}$.
    There exists a constant $c_{d,\bm{k}}$, such that,
    $$ \norm{\mathcal{\bar{P}}}_{\infty} \leq c_{d,\bm{k}},$$
    where, $c_{d,\bm{k}}$ is related to $d$, the dimension of the parametric domain $\Omega$, and $\bm{k} = (k_1, k_2, \cdots, k_d)$.
    Here, $k_j$ is the degree of the B-spline basis function $B_{i_j}(\eta_j), j=1,2,\cdots,d$~\pref{eq:tensor_basis}.
  \end{lem}

  Owing to Lemma~\ref{lem:ls_projector_bound}, the $L_2$ projector
    $\mathcal{P}: C(\Omega) \rightarrow \mathbb{S}$,
    from the continuous function space $C(\Omega)$ to the B-spline sub-space $\mathbb{S}$,
    is also bounded in $L_\infty$ norm.

  Thus, we have,
  \begin{lem}
    Denoting $dist(f, \mathbb{S})$ as the distance from $f$ to the B-spline
        sub-space $S$, we have,
     \begin{equation} \label{eq:distance}
        \norm{\mathcal{D}T - \mathcal{D}T_r}_\infty =
        \norm{f - \mathcal{P}f}_\infty \leq
        (1 + \norm{\mathcal{P}}_\infty) dist(f,\mathbb{S}).
     \end{equation}
  \end{lem}

  \textbf{Proof:} For an arbitrary $s \in \mathbb{S}$,
    it holds that $\mathcal{P}s = s$.
    Letting $\mathcal{I}$ is the identity operator, we have,
    \begin{equation*}
      \norm{f - \mathcal{P}f}_\infty =
      \norm{f - s + \mathcal{P} s - \mathcal{P}f}_\infty
      = \norm{(\mathcal{I} - \mathcal{P})(f-s)}_\infty
      \leq (1 + \norm{\mathcal{P}}_\infty) \norm{f-s}_\infty.
    \end{equation*}
  Because $s$ is an arbitrary function in $\mathbb{S}$,
    it can be so chosen that $\norm{f-s}_\infty = dist(f,\mathbb{S})$.
  Thus, Eq.~\pref{eq:distance} is proved. $\Box$

  \vspace{0.3cm}

  Therefore, if $dist(f, \mathbb{S}) \rightarrow 0$, when $h \rightarrow 0$,
    we get $\norm{\mathcal{D}T - \mathcal{D}T_r}_\infty \rightarrow 0$,
    when $h \rightarrow 0$.
  In other words, the IGA-L method is consistency.
  Here, $h$ is the knot grid size of $\mathcal{T}^h$,
    where the splines in $\mathbb{S}$ are defined on.
  In conclusion, the theorem follows.

  \begin{thm} \label{thm:consistency_special}
    Suppose the operators $\mathcal{D}$ and $\mathcal{G}$ have polynomial
        coefficients.
    If $dist(f, \mathbb{S}) \rightarrow 0$, when $h \rightarrow 0$,
        the IGA-L method is consistency.
  \end{thm}

 \subsubsection{Consistency in the generic case}
 \label{ssub:consistency_generic}

 In this section, we consider the case that the operators $\mathcal{D}$ and
    $\mathcal{G}$ in Eq.~\pref{eq:original_pde} are generic operators.

 Give three knot sequences,
 \begin{align}
        &\{\underbrace{u_0,u_0,\cdots,u_0}_{l_u+1},u_1,\cdots,u_{n_u-1},\underbrace{u_{n_u},u_{n_u},\cdots,u_{n_u}}_{l_u+1}\},
                \label{eq:knot_u} \\
        &\{\underbrace{v_0,v_0,\cdots,v_0}_{l_v+1},v_1,\cdots,v_{n_v-1},\underbrace{v_{n_v},v_{n_v},\cdots,v_{n_v}}_{l_v+1}\},
                \label{eq:knot_v} \\
        &\{\underbrace{w_0,w_0,\cdots,w_0}_{l_w+1},w_1,\cdots,w_{n_w-1},\underbrace{w_{n_w},w_{n_w},\cdots,w_{n_w}}_{l_w+1}\}.
                \label{eq:knot_w}
    \end{align}
 In one dimensional case, the numerical solution $T_r$ is defined on the knot
    sequence~\pref{eq:knot_u};
    in two dimensional case, $T_r$ is defined on the knot sequences~\pref{eq:knot_u} and~\pref{eq:knot_v};
    in three dimensional case, $T_r$ is defined on the knot sequences~\pref{eq:knot_u}, \pref{eq:knot_v}, and~\pref{eq:knot_w}.
 Denote $h$ as the knot grid size of the knot grid where $T_r$ is defined on.

 Let $R(\bm{\eta}) = (\mathcal{D}T(\bm{\eta}) - \mathcal{D}T_r(\bm{\eta}))^2,
    \bm{\eta} \in \Omega_0$.
 Denote $H_u = u_{n_u}-u_0, H_v = v_{n_v}-v_0, H_w = w_{n_w}-w_0$,
    and $e_h = \sum_k R(\bm{\vartheta}_k)$ as the least-squares fitting error,
    where $\bm{\vartheta}_k \in \Omega_0$ are collocation points.
 The following theorem holds.

 \begin{thm} \label{thm:error_formula}
  In the IGA-L method, if
  \begin{enumerate}
    \item [(1)] each knot interval of the NURBS function $T_r$
            defined on knot grid $\mathcal{T}^h$
            contains at least one collocation point, and,
    \item [(2)] the degree of each variable in $T_r$ is
            larger than the maximum order of the partial derivatives to the variables appearing
            in $\mathcal{D}$ (see~Eq.~\pref{eq:original_pde}),
  \end{enumerate}
    the fitting error of $\mathcal{D}T_r$ to $\mathcal{D}T$ in $L_2$ norm
        can be deduced as,
  \begin{enumerate}
    \item [(1)] in one dimensional case,
        \begin{equation} \label{eq:one_dim_error}
          \norm{\mathcal{D}T - \mathcal{D}T_r}_2^2 \leq
          h e_h + h H_u \abs{R'(\eta^*)},\ \text{where}\
          \eta^* \in (u_0, u_{n_u});
        \end{equation}
    \item [(2)] in two dimensional case,
        \begin{equation} \label{eq:two_dim_error}
        \begin{split}
          \norm{\mathcal{D}T - \mathcal{D}T_r}_2^2 \leq
          h^2 e_h + & h H_u H_v
          (\abs{R'_u(\bar{\eta}_1, \bar{\xi}_1)} +
           \abs{R'_v(\bar{\eta}_2, \bar{\xi}_2)}),\\
           & \qquad \text{where}\
          (\bar{\eta}_1, \bar{\xi}_1), (\bar{\eta}_2, \bar{\xi}_2)\ \text{are points in}\ (u_0, u_{n_u}) \times (v_0, v_{n_v});
        \end{split}
        \end{equation}
    \item [(3)] in three dimensional case,
        \begin{equation} \label{eq:three_dim_error}
        \begin{split}
          &\norm{\mathcal{D}T - \mathcal{D}T_r}_2^2  \leq
          h^3 e_h +  h H_u H_v H_w
          (\abs{R'_u(\tilde{\eta}_1, \tilde{\xi}_1)} +
           \abs{R'_v(\tilde{\eta}_2, \tilde{\xi}_2)} +
           \abs{R'_w(\tilde{\eta}_3, \tilde{\eta}_3}),\\
           & \text{where}\
          (\tilde{\eta}_1, \tilde{\xi}_1), (\tilde{\eta}_2, \tilde{\xi}_2), \text{and}\ (\tilde{\eta}_3, \tilde{\xi}_3)\ \text{are points in}\ (u_0, u_{n_u}) \times (v_0, v_{n_v}) \times (w_0, w_{n_w}).
        \end{split}
        \end{equation}
  \end{enumerate}
  \end{thm}

  The proof to the three
    formulae~\pref{eq:one_dim_error}-~\pref{eq:three_dim_error} in one-, two-, and three-dimensional cases are similar.
  We present the proof to the two-dimensional case
    (Eq.~\pref{eq:two_dim_error}) in Appendix A6.

  \vspace{0.3cm}

  Moreover, if the derivative or partial derivative of $R(\bm{\eta})$ is
    continuous, (then it is bounded in its domain),
    and the least-squares fitting error $e_h$ is also bounded,
    we have $\norm{\mathcal{D}T - \mathcal{D}T_r}_2^2 \rightarrow 0, (h \rightarrow 0)$,
    based on Theorem~\ref{thm:error_formula}.
  That is, the IGA-L method is consistency.
  This leads to the following theorem.
 \begin{thm} \label{thm:consistency}
    If the least-square fitting error $e_h$ is bounded,
        $R(\bm{\eta}) \in C^1(\Omega_0)$,
        and the conditions in Theorem~\ref{thm:error_formula} are satisfied,
        then the IGA-L method is consistency.
 \end{thm}

 \subsection{Convergence}
 \label{subsec:convergence}

 Based on the consistency property of IGA-L,
    we can show that IGA-L is convergent if the
    differential operator is stable or strongly monotonic,
    similarly as in Ref.~\cite{lin2013consistency}.

 Let $V$ and $W$ be two Hilbert spaces
    and $\norm{\cdot}_V$ and $\norm{\cdot}_W$ be two norms defined
    on $V$ and $W$, respectively.
 Suppose $\norm{\cdot}_V$ and $\norm{\cdot}_W$ are equivalent to the
    $L^2$ norm, i.e.,
    there exist positive constants $\alpha_V,\ \beta_V,\ \alpha_W,$
    and $\beta_W$ such that,
    \begin{align} \label{eq:equivalence}
        & \alpha_V \norm{\cdot}_{L^2} \leq \norm{\cdot}_V \leq
            \beta_V \norm{\cdot}_{L^2}, \\
        & \alpha_W \norm{\cdot}_{L^2} \leq \norm{\cdot}_W \leq
            \beta_W \norm{\cdot}_{L^2}.
    \end{align}

 We first give the definitions for stable and strongly monotonic operators.

 \begin{Def}[Stability estimate and stable operator~\cite{solin2006partial}]
    \label{def:stable_ope}
 Let $V,W$ be Hilbert spaces and $\mathcal{D}: V \rightarrow W$
    be a differential operator.
 If there exists a constant $C_{\mathcal{S}}>0$ such that
 \begin{equation} \label{eq:stability_estimate}
     \norm{\mathcal{D}v}_W \geq C_{\mathcal{S}} \norm{v}_V,\ \quad \text{for all}\ v \in D(\mathcal{D}),
\end{equation}
    where $D(\mathcal{D})$ represents the domain of $\mathcal{D}$,
    differential operator $\mathcal{D}$ is called the stable operator,
    and the inequality~\pref{eq:stability_estimate} is called the
    stability estimate.
 \end{Def}

 \begin{Def}[Strongly monotonic operator~\cite{solin2006partial}]
 Let $V$ be a Hilbert space and $\mathcal{D} \in \mathcal{L}(V,V')$.
 Operator $\mathcal{D}$ is said to be a strongly monotonic operator,
    if there exists a constant $C_{\mathcal{D}} > 0$, such that
    \begin{equation} \label{eq:strong_monotone}
        \langle \mathcal{D}v,v \rangle \geq C_{\mathcal{D}} \norm{v}_V^2,\ \text{for
        all}\ v \in V.
    \end{equation}
 For every $v \in V$, element $\mathcal{D}v \in V'$ is of a linear
    form.
 The symbol $\langle \mathcal{D}v,v \rangle$,
    which denotes the application of $\mathcal{D}v$ to $v \in V$,
    is called a duality pairing.
 \end{Def}

 Clearly, if a differential operator $\mathcal{D}$ is strongly
    monotonic, it is stable.

 \begin{lem} \label{lem:strong_monotone}
 Let $V$ be a Hilbert space and $\mathcal{D} \in \mathcal{L}(V,V')$ be a
    continuous strongly monotonic linear operator.
 Then, there exists a constant $C_{\mathcal{D}} > 0$ such that $\mathcal{D}$
 satisfies the stability estimate~\pref{eq:stability_estimate}~\cite{solin2006partial}.
\end{lem}

 The proof can be found in Ref.~\cite{lin2013consistency}.

 \vspace{0.3cm}

  Therefore, we have the convergence property of IGA-L as
    follows.

 \begin{thm} \label{thm:stable_operator_convergence}
 Suppose NURBS function $T_r$, defined on knot grid $\mathcal{T}^h$, is
    the numerical solution to the boundary value
    problem~\pref{eq:original_pde},
    generated by IGA-L,
    and the conditions presented in
    Theorem~\ref{thm:consistency}
    are satisfied in one, two, and three dimensions, respectively.
 If differential operator $\mathcal{D}: V \rightarrow W$ in~\pref{eq:original_pde}
    is a stable operator,
    $T_r$ will converge to analytic solution $T$
    when the knot grid size $h \rightarrow 0$.
\end{thm}

 \textbf{Proof:} Differential operator
    $\mathcal{D}$ in~\pref{eq:original_pde} is a stable operator,
    so there exists a constant $C_{\mathcal{S}}>0$, such that
    \begin{equation*}
        \norm{\mathcal{D}(T-T_r)}_W \geq C_{\mathcal{S}} \norm{T-T_r}_V.
    \end{equation*}
 And it is equivalent to
 \begin{equation*}
    \norm{T-T_r}_V \leq \frac{1}{C_{\mathcal{S}}} \norm{\mathcal{D}T-\mathcal{D}T_r}_W.
 \end{equation*}
 Due to the equivalence of $\norm{\cdot}_{L^2}$ and
    $\norm{\cdot}_W$~\pref{eq:equivalence},
    we have,
    \begin{equation*}
        \norm{T-T_r}_V \leq \frac{1}{C_{\mathcal{S}}} \norm{\mathcal{D}T-\mathcal{D}T_r}_W
        \leq
        \frac{\beta_W}{C_{\mathcal{S}}}\norm{\mathcal{D}T-\mathcal{D}T_r}_{L^2}.
    \end{equation*}
 Because of the consistency of IGA-L
    (Theorem~\ref{thm:consistency}),
    $\norm{T-T_r}_V$ will converge to $0$ when $h \rightarrow 0$.
 And this theorem is proved.

    \rightline{$\Box$}

 Moreover, Theorem~\ref{thm:stable_operator_convergence}
    and Lemma~\ref{lem:strong_monotone} lead to the direct
    corollary.

 \begin{cor}
 Suppose NURBS function $T_r$ defined on knot grid $\mathcal{T}^h$ is
    the numerical solution to the boundary value problem~\pref{eq:original_pde},
    generated by IGA-L,
    and the conditions presented in
    Theorem~\ref{thm:consistency}
    are satisfied in one, two, and three dimensions, respectively.
 Additionally, suppose norm $\norm{\cdot}_{L^2}$ bounds norm $\norm{\cdot}_{V'}$.
 If differential operator $\mathcal{D}$ in~\pref{eq:original_pde}
    is a strongly monotonic operator,
    then $T_r$ will converge to analytic solution $T$
    when $h \rightarrow 0$.
 \end{cor}

 \vspace{0.3cm}

 It is well known that a wide class of elliptic
    differential operators are stable or strongly monotonic.
    Hence, IGA-L is convergent for equations that have these elliptic
    differential operators.
 The examples of a PDE whose differential
    operators are strongly monotonic can be found in
    Refs.~\cite{solin2006partial,lin2013consistency}.

 \section{Comparisons and discussions}
\label{sec:comparison}

 \subsection{Theoretical comparisons}
 In the following, we compare IGA-L with IGA-C and IGA-SC
    in terms of their computational efficiency at solving a scalar problem
    (Laplace equation~\cite{schillinger2013isogeometric})
    and vector problem (elasticity equation~\cite{schillinger2013isogeometric}).
 We consider model discretizations in one, two, and three dimensions that
    are characterized by the degree of the basis functions and
    the numbers of control and collocation points in each parametric direction.
 For the sake of simplicity, we assume that the model discretizations in
    two and three dimensions have the same number of collocation points
    (and control points) in each parametric direction,
    and the degrees of basis functions in one, two, and three dimensions are
    $p,\ p \times p$, and $p \times p \times p$, respectively.

\begin{table}[!htb]
 \caption{Cost in flops for formation at one collocation point in
 IGA-C, IGA-SC, and IGA-L.} \label{tbl:cost_of_col_pt}
 \begin{center}
 \begin{tabular}{c | c | c }
 \hline
 Dimension & A scalar problem (Laplace) & A vector problem
 (elasticity)\\
 \hline
  & \multicolumn{2}{c}{Solve for 1st derivatives}\\
 \hline
 $1$ & \multicolumn{2}{c}{$(p+1)$} \\
 $2$ & \multicolumn{2}{c}{$5(p+1)^2+4$}\\
 $3$ & \multicolumn{2}{c}{$12(p+1)^3+16$}\\
 \hline
  & \multicolumn{2}{c}{Compute right hand side vectors and solve for 2nd derivatives}\\
 \hline
 $1$ & \multicolumn{2}{c}{$3(p+1)$} \\
 $2$ & \multicolumn{2}{c}{$24(p+1)^2+16$}\\
 $3$ & \multicolumn{2}{c}{$87(p+1)^3+140$}\\
 \hline
 & \multicolumn{2}{c}{Total number of flops for basis function}\\
 \hline
 $1$ & \multicolumn{2}{c}{$35(p+1)+1$} \\
 $2$ & \multicolumn{2}{c}{$124(p+1)^2+33$}\\
 $3$ & \multicolumn{2}{c}{$302(p+1)^3+219$}\\
 \hline
 & & Evaluate Navier's eqs. on global level\\
 \hline
 $1$ & & $(p+1)$ \\
 $2$ & & $12(p+1)^2$\\
 $3$ & & $21(p+1)^3$\\
 \hline
 & \multicolumn{2}{c}{Total number of flops to evaluate the local stiffness matrix}\\
 \hline
 $1$ & $35(p+1)+1$ & $36(p+1)+1$ \\
 $2$ & $125(p+1)^2+33$ & $134(p+1)^2+33$\\
 $3$ & $304(p+1)^3+219$ & $323(p+1)^3+219$\\
 \hline
 \end{tabular}
 \end{center}
 \end{table}

 First, the costs in floating point operations (flops) for the
    formation at one collocation point are the same for
    IGA-L, IGA-C, and IGA-SC.
 These costs are listed in Table~\ref{tbl:cost_of_col_pt}.

 Second, we compared the costs (in flops) of solving the linear
    system of equations in IGA-C IGA-SC, and IGA-L,
    and present the comparison in Table~\ref{tbl:cost_comp}.
 In this table, the first column is the dimension of the problem
    solved,
    and the second column is the number of control
    points (equal to the number of collocation points in IGA-C).
 The third column is the cost in flops to solve the linear system of
    equations using Gaussian elimination in IGA-C~\cite{golub1996matrix}.
 Moreover, the fourth column is the number of collocation points in IGA-L,
    and the fifth column is the cost in flops to solve the normal
    equation~\pref{eq:normal_equ} using Cholesky decomposition in
    IGA-L~\cite{golub1996matrix}.
 Finally, the last two columns are the number of collocation points
    in IGA-SC,
    and the cost in flops to solve the normal equation using Cholesky decomposition in IGA-SC~\cite{golub1996matrix}, respectively.
 We can see that the cost to solve Eq.~\pref{eq:normal_equ} in IGA-L linearly
    increases with the number of
    collocation points ($m,\ m \times m,\ m \times m \times m$).

\begin{table}[!htb]
\begin{center}
\small
\begin{threeparttable}
 \caption{Cost comparison of IGA-C, IGA-SC, and IGA-L.}
 \label{tbl:cost_comp}
\begin{tabular}{c | c | c | c | c | c | c }
 \hline
 \multirow{2}{*}{Dim.\tnote{1}} & \multirow{2}{*}{\#Cont.}\tnote{2} & IGA-C & \multicolumn{2}{|c}{IGA-L} & \multicolumn{2}{|c}{IGA-SC}\\
 \cline{3-7}
  &    & Cost\tnote{3}& \#Col.\tnote{4} & Cost\tnote{3} & \#Col.\tnote{4} & Cost\tnote{3} \\
 \hline
 \multirow{2}{*}{$d=1$} & \multirow{2}{*}{$n$} & \multirow{2}{*}{$2n^3/3$} & \multirow{2}{*}{$m$} & \multirow{2}{*}{$ n^3/3 + mn^2 $}
 & $2(n-p)$ (odd $p$) & $n^3/3 + (2n-p)n^2$ \\
 \cline{6-7}
     &    &    &     &    & $2(n-p)-1$ (even $p$) & $n^3/3 + (2(n-p)-1)n^2$  \\
 \hline
 \multirow{2}{*}{$d=2$} & \multirow{2}{*}{$n^2$} &\multirow{2}{*}{$2n^6/3$}   & \multirow{2}{*}{$m^2$} & \multirow{2}{*}{$n^6/3 + m^2n^4$}
 & $(2(n-p))^2$ (odd $p$) & $n^6/3+(2(n-p))^2 n^4$\\
 \cline{6-7}
     &   &    &     &     & $(2(n-p)-1)^2$ (even $p$)
     & $n^6/3 + (2(n-p)-1)^2 n^4$ \\
 \hline
 \multirow{2}{*}{$d=3$} & \multirow{2}{*}{$n^3$} &\multirow{2}{*}{$2n^9/3$}  & \multirow{2}{*}{$m^3$} & \multirow{2}{*}{$n^9/3 + m^3 n^6$}
 & $(2(n-p))^3$ (odd $p$) & $n^9/3 + (2(n-p))^3 n^6$ \\
 \cline{6-7}
    &    &    &     &     & $(2(n-p)-1)^3$ (even $p$)
    & $n^9/3 + (2(n-p)-1)^3 n^6$ \\
 \hline
 \end{tabular}
\begin{tablenotes}
 \item [1] Dimension.
 \item [2] Number of control points.
 \item [3] Cost in flops.
 \item [4] Number of collocation points.
 \end{tablenotes}
 \end{threeparttable}
 \end{center}
 \end{table}

\subsection{Numerical comparisons}

 In this section, we compare IGA-L with IGA-C and IGA-SC using
    some numerical examples.
 To measure the approximation accuracy,
    we define the error formulae,
    i.e., the \emph{relative error for the solution $T_r$},
    \begin{equation} \label{eq:rel_err_func}
        e_T = \sqrt{\frac{\int_{\Omega} (T-T_r)^t (T-T_r) \rm d \Omega}
                    {\int_{\Omega} T^t T \rm d \Omega}}.
    \end{equation}
 Additionally, to illustrate the error distribution of the
    numerical solution,
    the following \emph{absolute errors} $e_a$ are employed, i.e.,
    \begin{equation} \label{eq:abs_err}
    \begin{split}
        & e_a(u)  = \abs{T(u) - T_r(u)},\ \text{for one-dimensional
        case},\\
        & e_a(u,v)  = \abs{T(u,v) - T_r(u,v)},\ \text{for two-dimensional
        case},\\
        & e_a(u,v,w)  = \abs{T(u,v,w) - T_r(u,v,w)},\ \text{for three-dimensional
        case}.
    \end{split}
    \end{equation}

 In the following, six numerical examples are presented.
 These examples are implemented with MATLAB and run on a PC
    with a 2.66-GHz Intel Core2 Quad CPU
    Q9400 and 3 GB memory.
 Examples I--III are three source problems in one, two, and
    three dimensions, respectively,
    Example IV is a linear elasticity problem,
    and Example V demonstrates the stability of the IGA-L method
    with respect to that of the IGA-C method.
 Moreover, Example VI illustrates the capability of IGA-L method for solving
    a 2D source problem on a frame-corner-like domain which contains a $C^0$ line.
 The problems in Examples I--V are solved by IGA-L, IGA-C, and IGA-SC
    methods.
 With the IGA-C method, the control points are increased gradually,
    and the collocation points are the Greville
    abscissae~\cite{auricchio2010isogeometric} of the knot vectors,
    also called the \emph{Greville collocation points}.
 The collocation manner for IGA-SC follows the method developed
    in~\cite{anitescu2015isogeometric}.
 With the IGA-L method,
    the control points are variable and
    are increased gradually at the same rate as those of IGA-C.
 In each computation round of the IGA-L variable strategy,
    supposing the numbers of the control points are $n$, $n \times n$,
    and $n \times n \times n$ in the one-, two-, and three-dimensional
    cases, respectively,
    the numbers of the corresponding collocation points are taken as $n+2$, $(n+2) \times (n+2)$,
    and $(n+2) \times (n+2) \times (n+2)$, respectively.
 In this paper, we take the following collocation manner for IGA-L method.

 \textbf{Collocation manner for IGA-L:}
 The collocation points for IGA-L are taken as the
    Greville abscissae of a knot sequence.
 To produce $n$ Greville collocation points for a NURBS curve of degree
    $k,\ (n \geq k)$,
    we first uniformly insert $n-k-1$ numbers into the interval $[0,1]$, resulting in the knot sequence,
    \begin{equation*}
      \underbrace{0,0,\cdots,0}_k, \tfrac{1}{n-k}, \tfrac{2}{n-k}, \cdots,
      \tfrac{n-k-1}{n-k}, \underbrace{1,1,\cdots,1}_k,
    \end{equation*}
 and then, $n$ Greville collocation points for a NURBS curve of degree
    $k,\ (n \geq k)$ can be generated from the above knot sequence.
 The collocation points for NURBS surfaces and solids can be produced by
    the aforementioned manner for each parameter.

 \begin{figure}[!htb]
    \centering
  \subfigure[]{
    \label{subfig:one_dim_solution}
    \includegraphics[width = 0.41\textwidth]
        {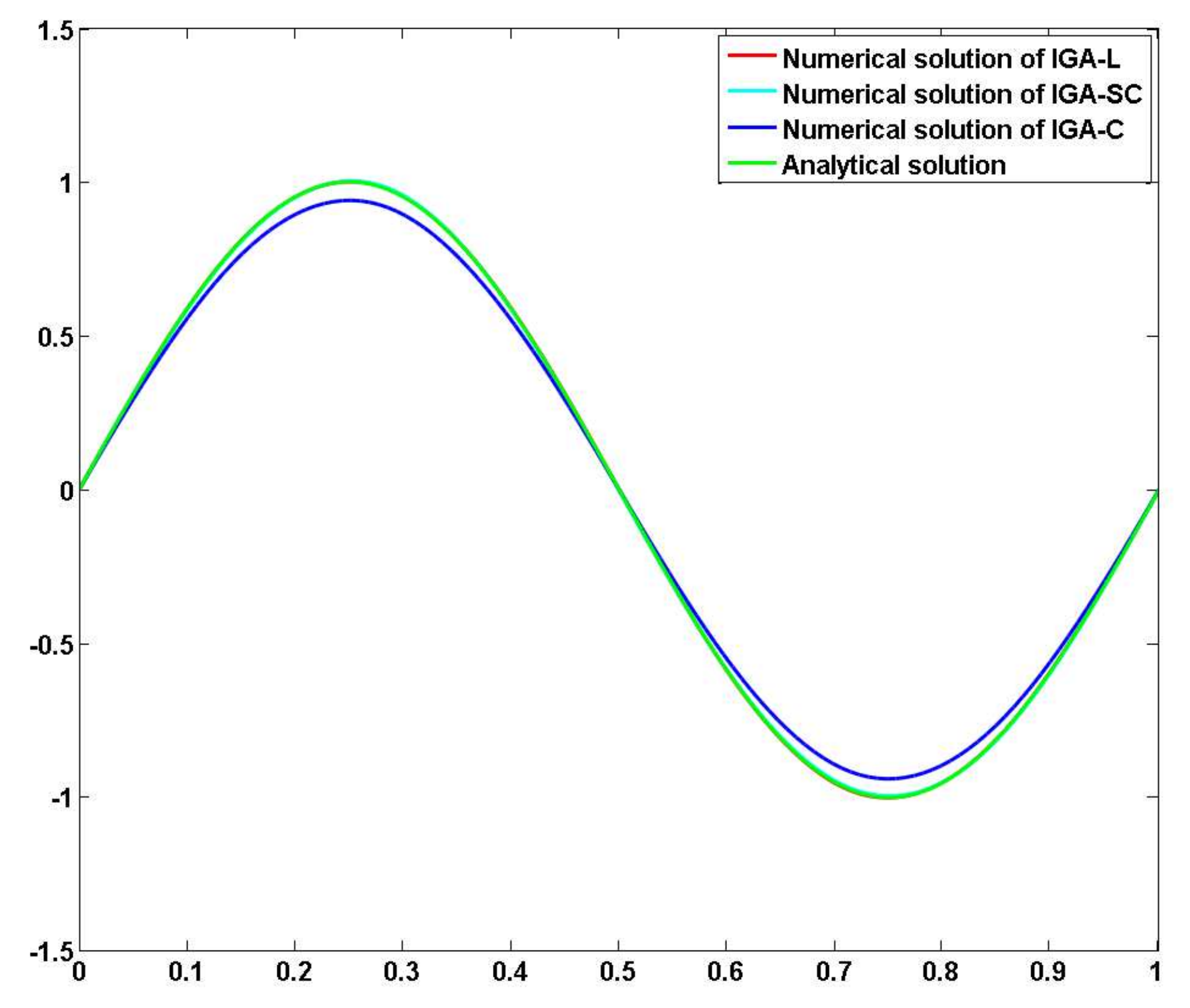}}
  \subfigure[]{\label{subfig:one_dim_absolute_error}
    \includegraphics[width = 0.41\textwidth]
        {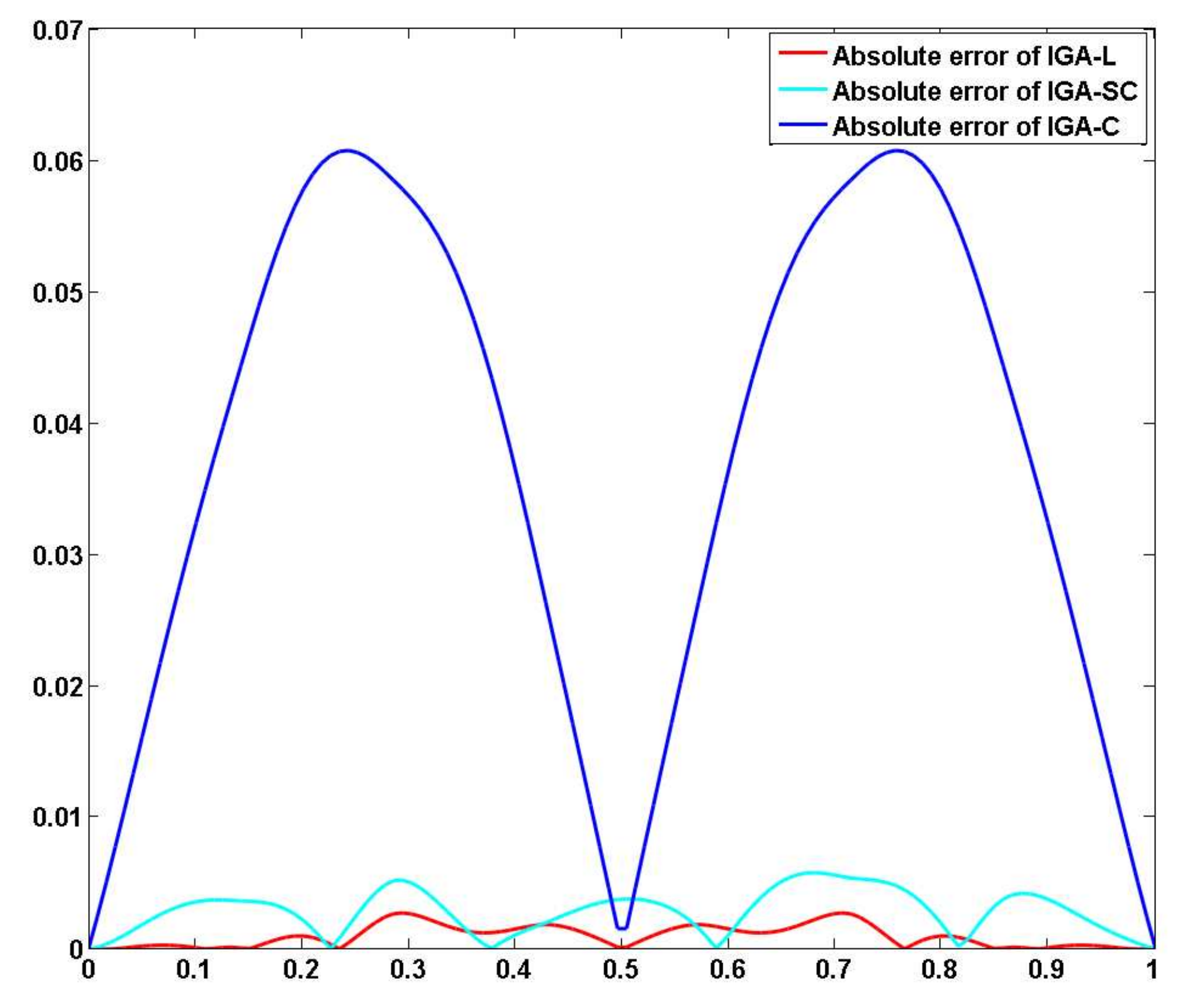}}
  \caption{Comparison of the analytical, IGA-L, IGA-SC and IGA-C
    solutions of Eq.~\pref{eq:exmp_one_dim}.
  (a) Analytical solution, IGA-L, IGA-SC and IGA-C solutions.
    Note that the analytical solution almost overlaps the IGA-L and IGA-SC solutions.
  (b) Absolute error distribution curves of the IGA-L, IGA-SC and IGA-C
    solutions.}
  \label{fig:one_dim_solutions}
\end{figure}

 \begin{figure}[!htb]
    \centering
  \subfigure[]{
    \label{subfig:1d_iga_l}
    \includegraphics[width = 0.46\textwidth, height = 0.31\textwidth]
        {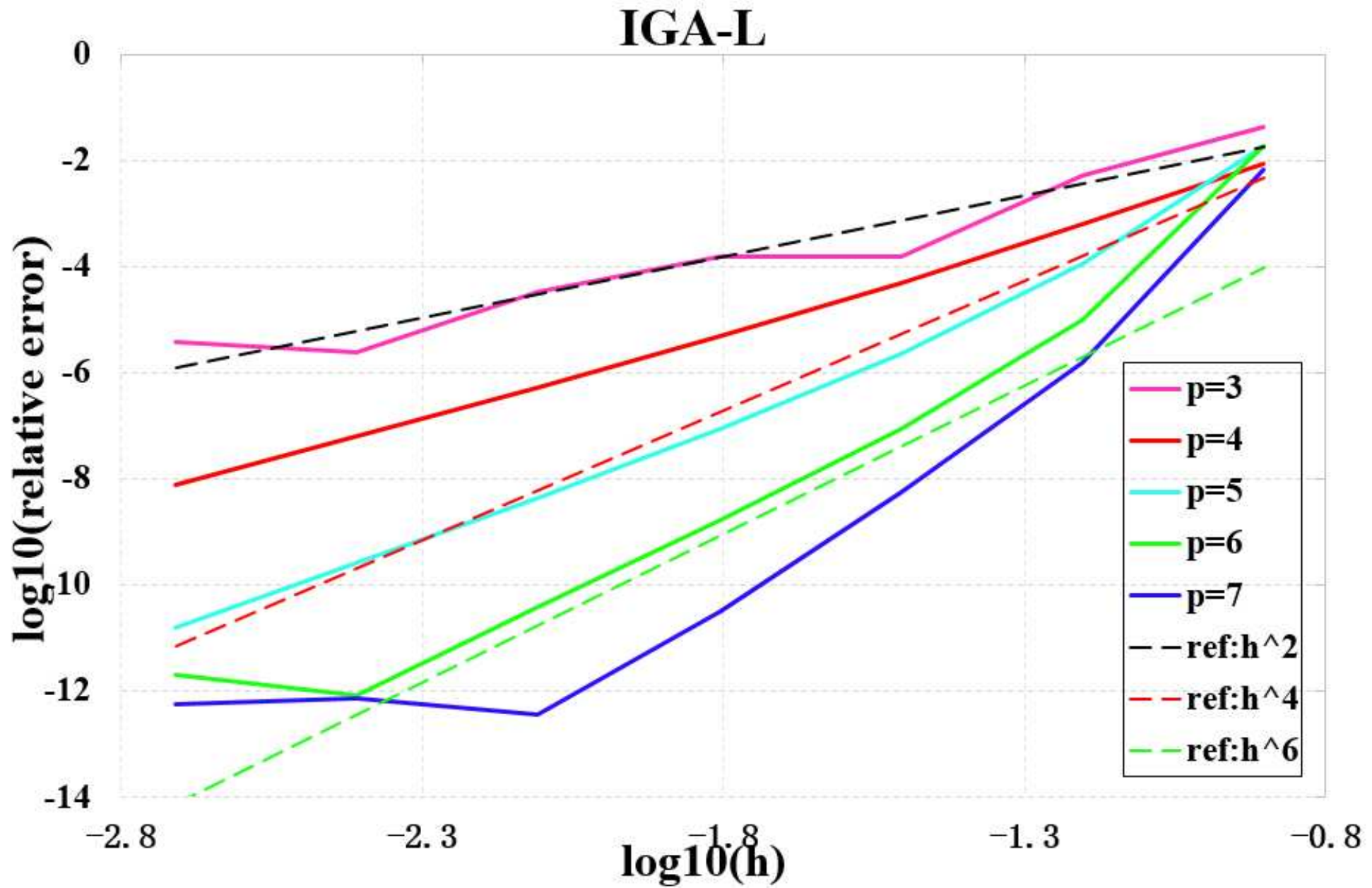}}
  \subfigure[]{\label{subfig:1d_iga_c}
    \includegraphics[width = 0.46\textwidth, height = 0.31\textwidth]
        {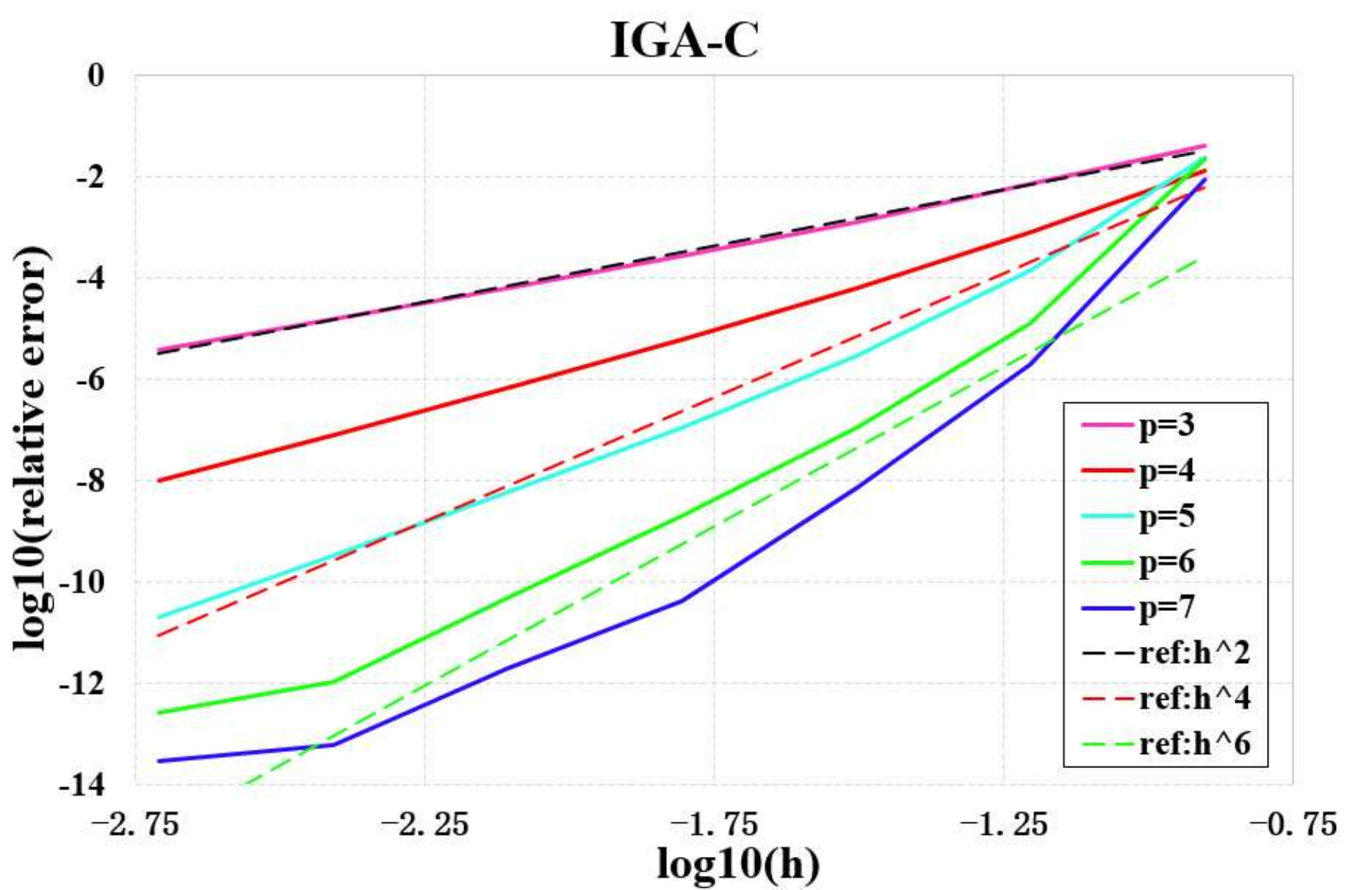}}
  \subfigure[]{\label{subfig:1d_iga_sc}
    \includegraphics[width = 0.46\textwidth, height = 0.31\textwidth]
        {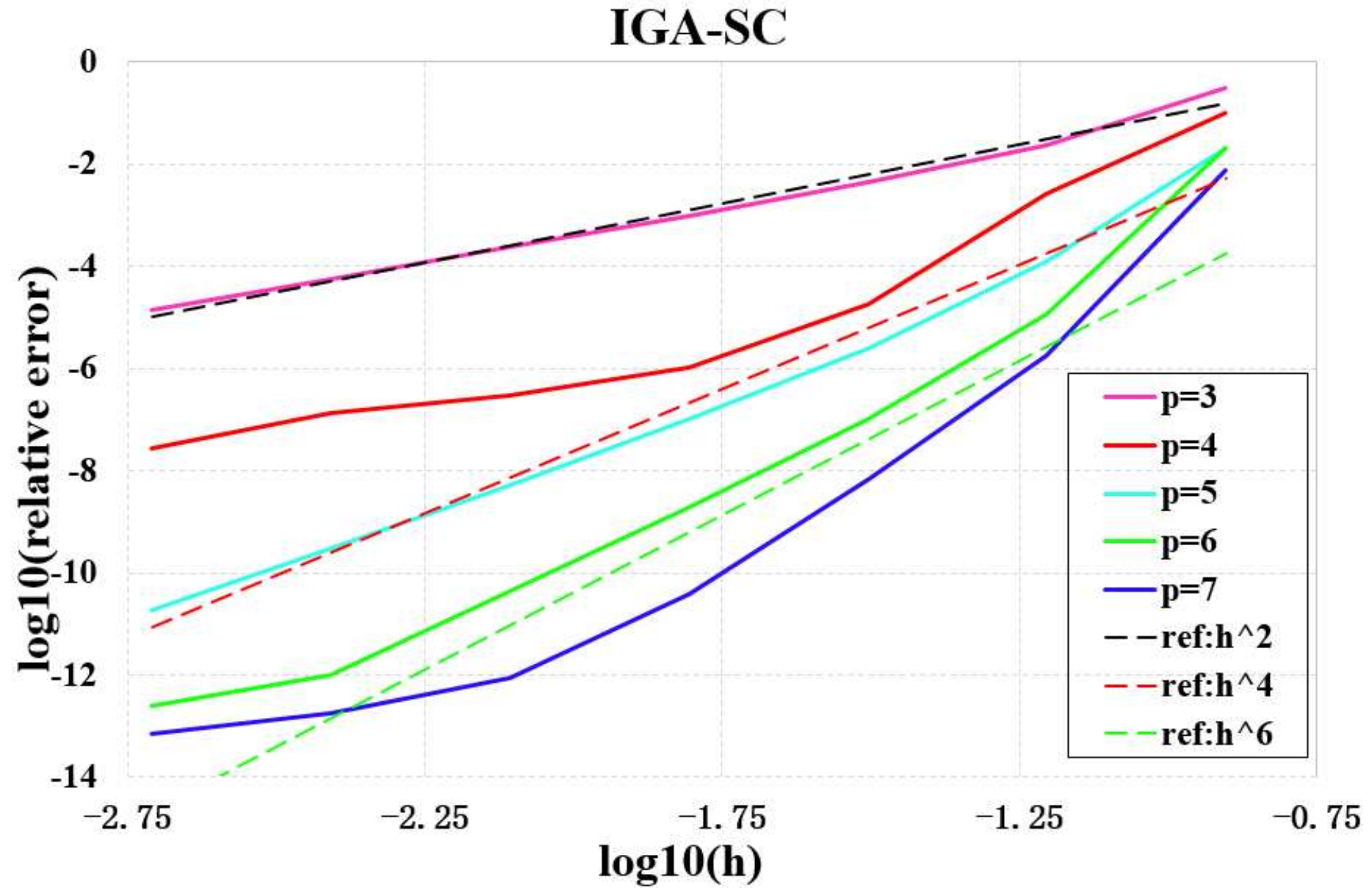}}
  \subfigure[]{\label{subfig:1d_time_error}
    \includegraphics[width = 0.46\textwidth, height = 0.30\textwidth]
        {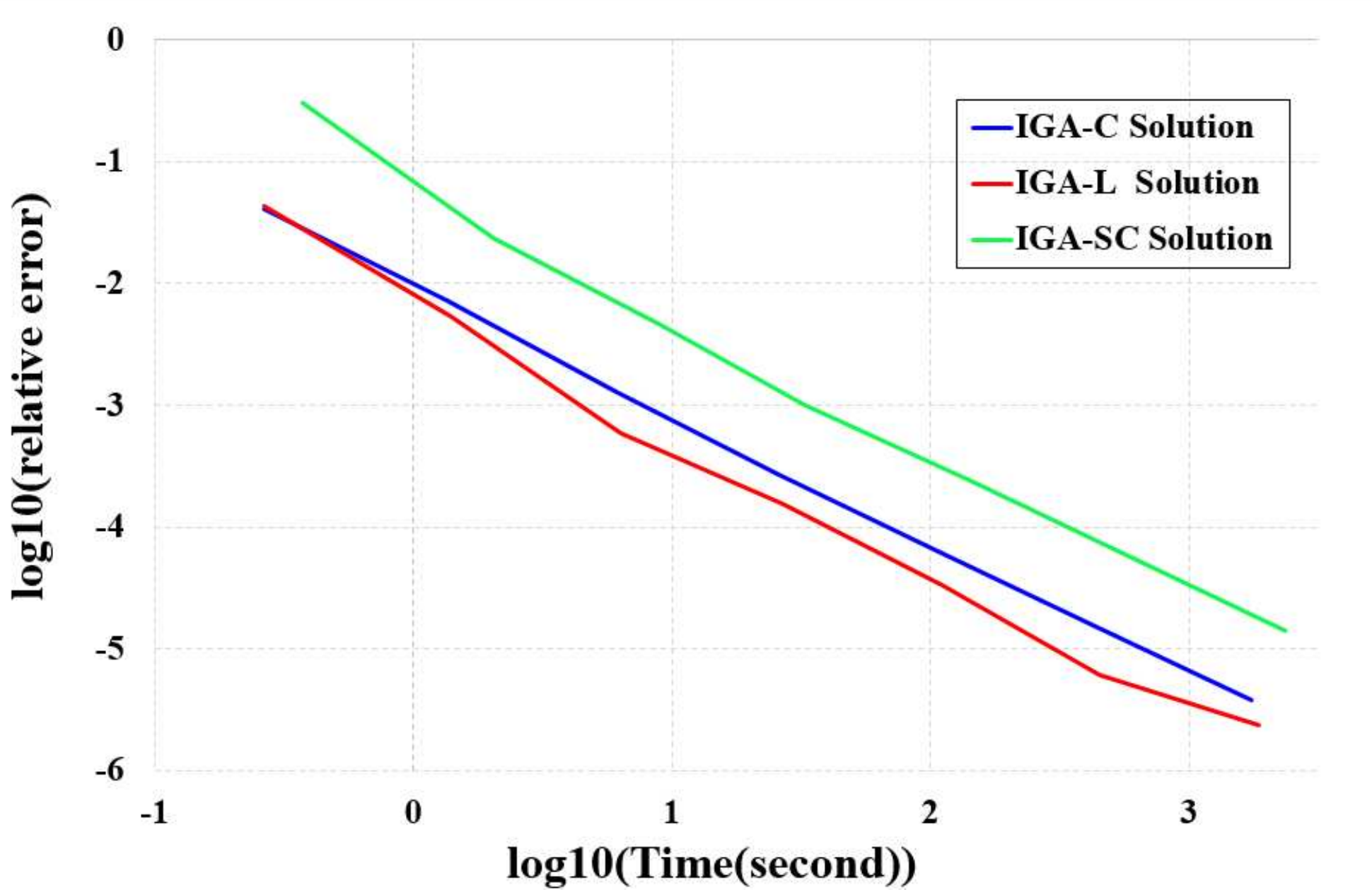}}
  \caption{Numerical results for the 1D source
        problem~\pref{eq:exmp_one_dim}.
        Diagrams of $\log_{10}(h)$ v.s. $\log_{10}$(relative error) for IGA-L (a), IGA-C (b), and IGA-SC (c), respectively. And, diagram of $\log_{10}$(Time) v.s. $\log_{10}$(relative error) (d).}
  \label{fig:one_dim_example}
\end{figure}

 \vspace{0.3cm}

 \textbf{Example I:}
 one-dimensional source problem
    with Dirichlet boundary condition:
 \begin{equation}
    \label{eq:exmp_one_dim}
    \begin{cases}
        &-T'' + T = (1 + 4 \pi^2) \sin(2 \pi x),\ x \in \Omega = [0,1],\\
        &T(0) = 0,\ T(1) = 0.
    \end{cases}
 \end{equation}
 This problem's analytical solution is $T = \sin(2 \pi x)$.
 The physical domain $\Omega = [0,1]$ of the boundary problem~\pref{eq:exmp_one_dim}
    is represented by a cubic B-spline curve.
 For the mathematical representation of the cubic B-spline curve,
    please refer to Appendix A1.

 The analytical, IGA-L, IGA-SC and IGA-C solutions for the one-dimensional
    source problem~\pref{eq:exmp_one_dim} with cubic B-spline,
    are illustrated in Fig.~\ref{subfig:one_dim_solution},
    where the IGA-L solution was generated with $10$ control points and $14$
    Greville collocation points,
    IGA-SC solution was generated with $10$ control point and $14$ collocation points,
    and IGA-C solution was generated with $10$ control points.
 The relative errors for the IGA-L, IGA-SC, and IGA-C solutions are $0.0018$,
    $0.0023$, and $0.0598$.
 The relative error for the IGA-L solution is one order of
    magnitude less than that of the IGA-C solution.
 In addition, Fig.~\ref{subfig:one_dim_absolute_error} demonstrates the
    absolute error distribution curves of the IGA-L, IGA-SC, and IGA-C solutions, respectively.
 The maximum absolute errors of IGA-L, IGA-SC, and IGA-C solutions are
    $0.0028$, $0.0034$, and $0.0607$, respectively.

 Moreover, diagrams of $\log_{10}(h)$ v.s. $\log_{10}$(relative error)
    for IGA-L, IGA-C, and IGA-SC, are illustrated in Figs.~\ref{subfig:1d_iga_l}-~\ref{subfig:1d_iga_sc}, respectively.
 From these diagrams, it can be seen that, in solving the one-dimensional
    source problem~\pref{eq:exmp_one_dim},
    the convergence rates of IGA-L, IGA-C, and IGA-SC are all $O(h^p)$ for even $p$, and $O(h^{p-1})$ for odd $p$.
 Additionally, the diagrams of
    $\log_{10}$(time) v.s. $\log_{10}$(relative error) for IGA-L, IGA-SC, and IGA-C methods are demonstrated in Fig.~\ref{subfig:1d_time_error},
    and the performance of the IGA-L method is the best.


\begin{figure}[!h]
\centering
 \subfigure[IGA-L solution of Eq.~\pref{eq:exmp_two_dim} with relative error $3.84 \times 10^{-4}$.]{
    \label{subfig:2d-igl-solution}
  \includegraphics[width = 0.30\textwidth]
        {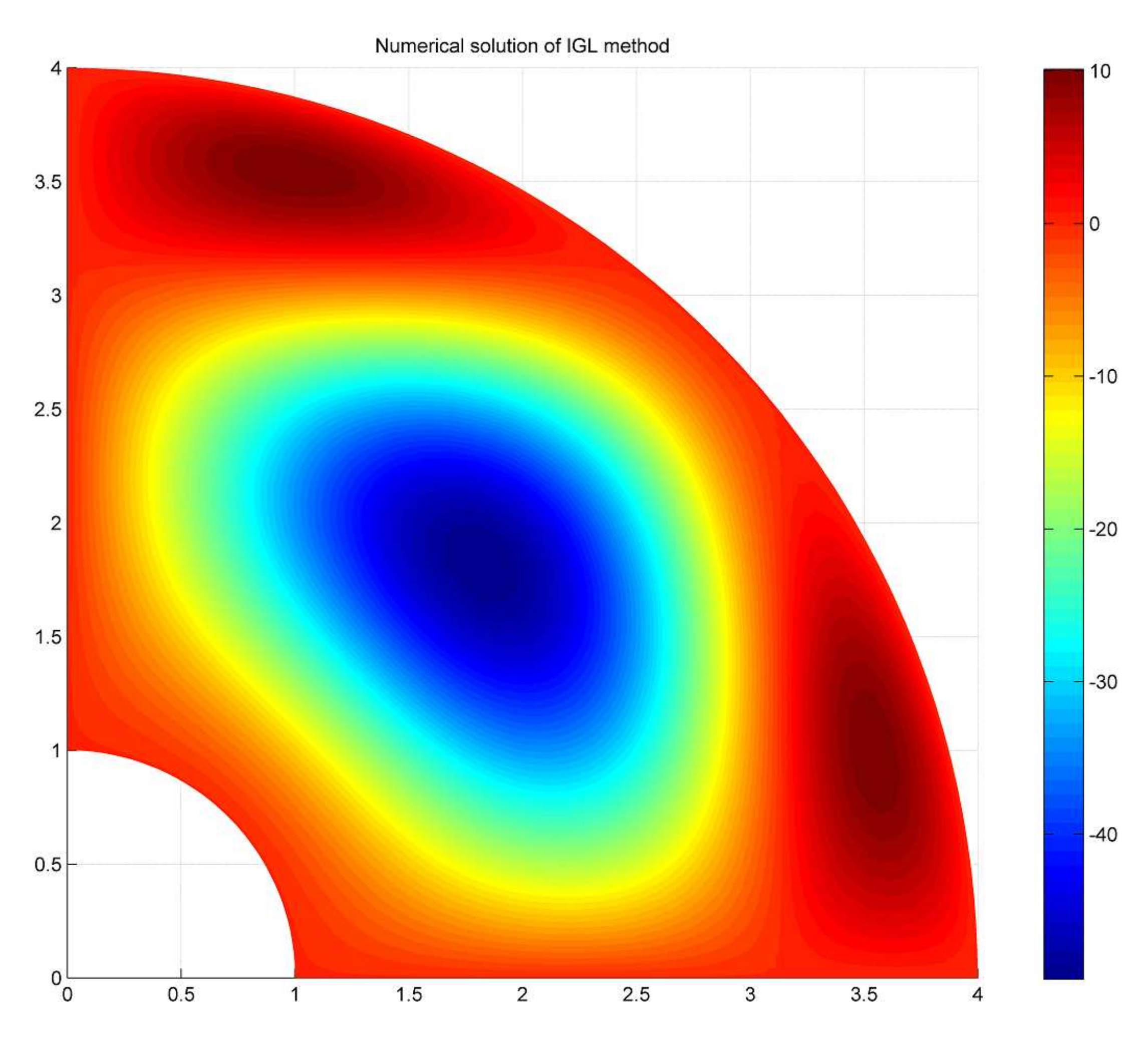}}
  \subfigure[IGA-C solution of Eq.~\pref{eq:exmp_two_dim} with relative error $1.59 \times 10^{-2}$.]{
    \label{subfig:2d-igc-solution}
  \includegraphics[width = 0.31\textwidth]
        {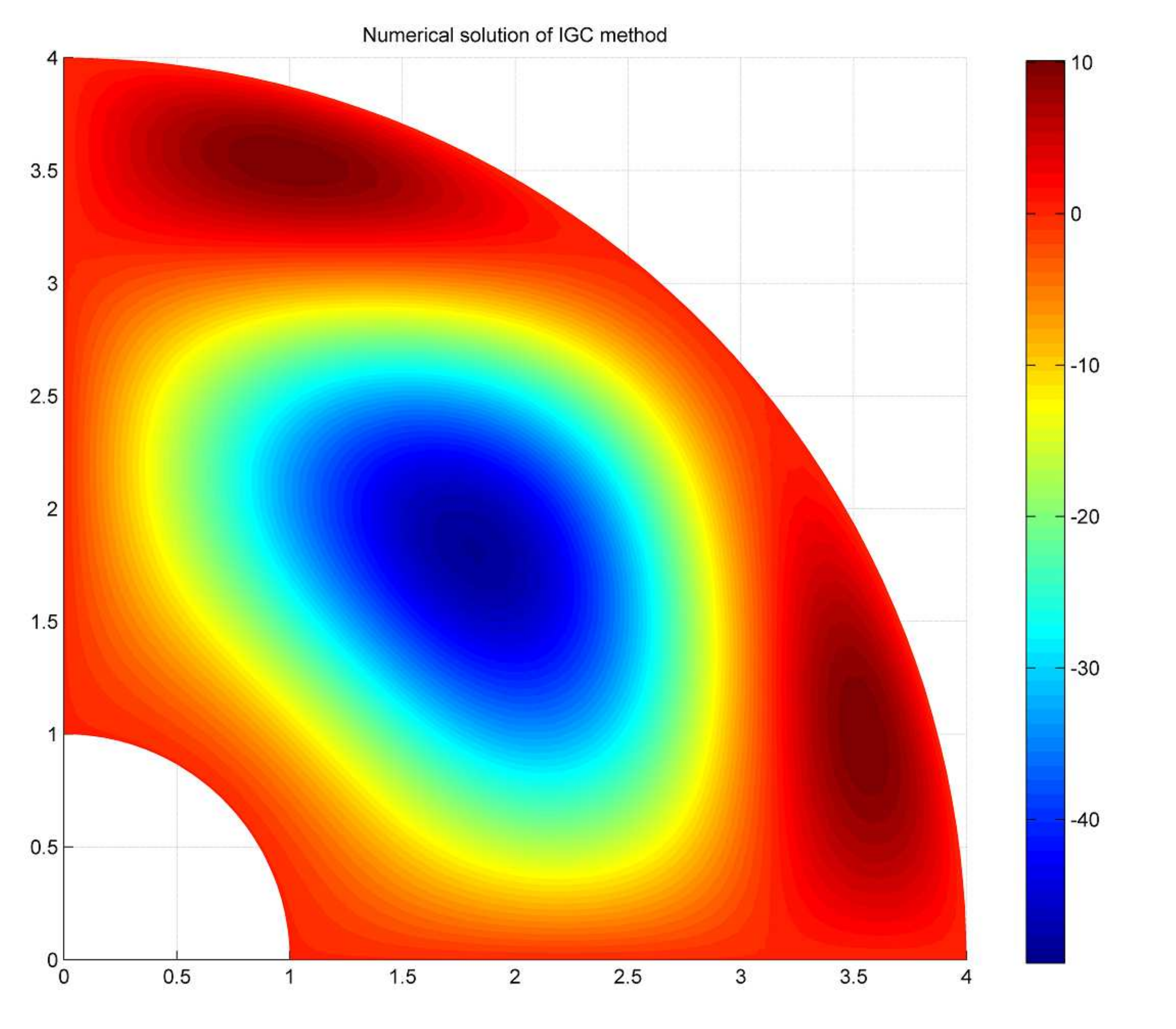}}
  \subfigure[IGA-SC solution of Eq.~\pref{eq:exmp_two_dim} with relative error $4.62 \times 10^{-4}$.]{
    \label{subfig:2d-igsc-solution}
  \includegraphics[width = 0.33\textwidth]
        {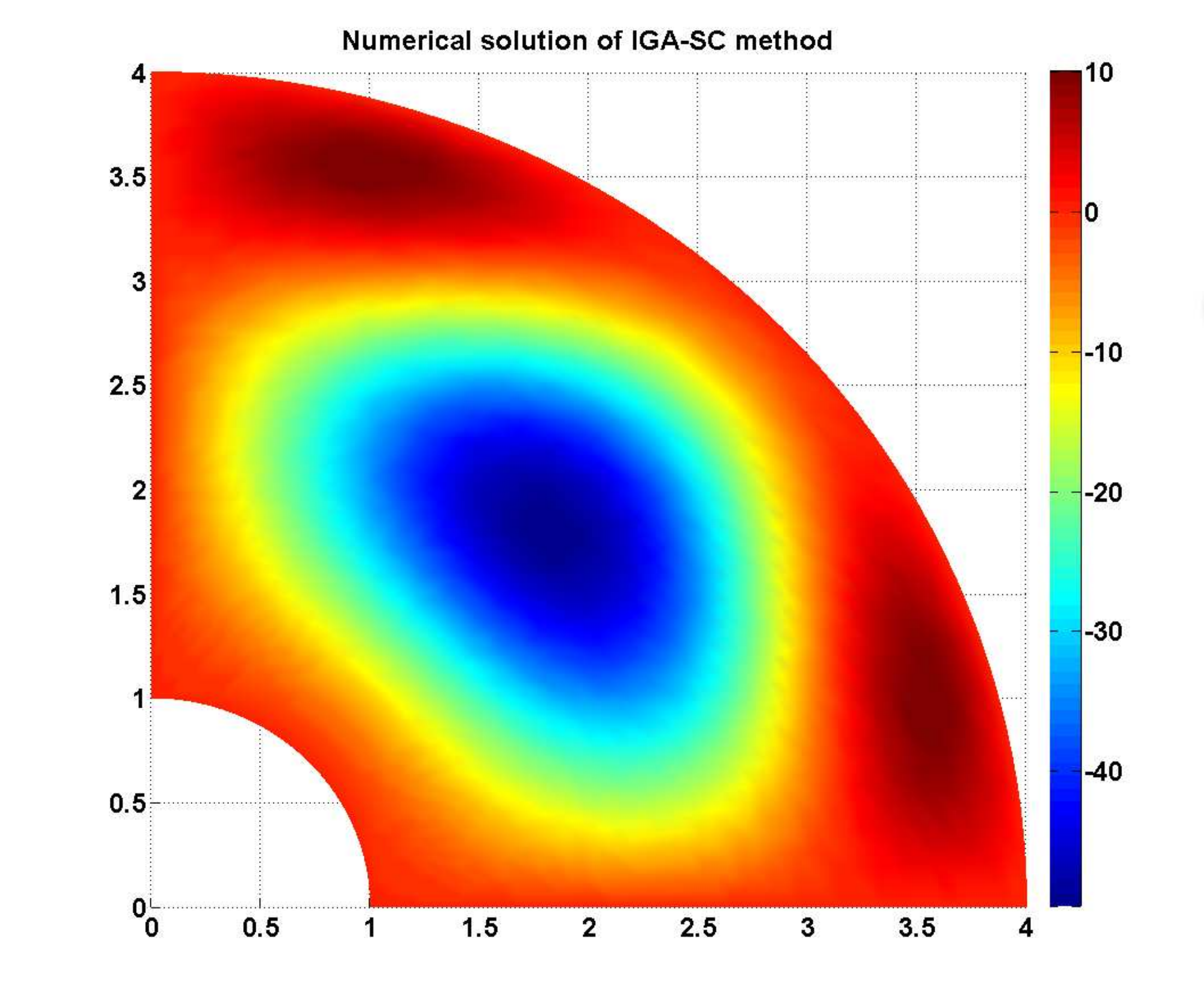}}
  \subfigure[Absolute error distribution of the IGA-L solution.]{
    \label{subfig:2d-igl-err}
  \includegraphics[width = 0.30\textwidth]
        {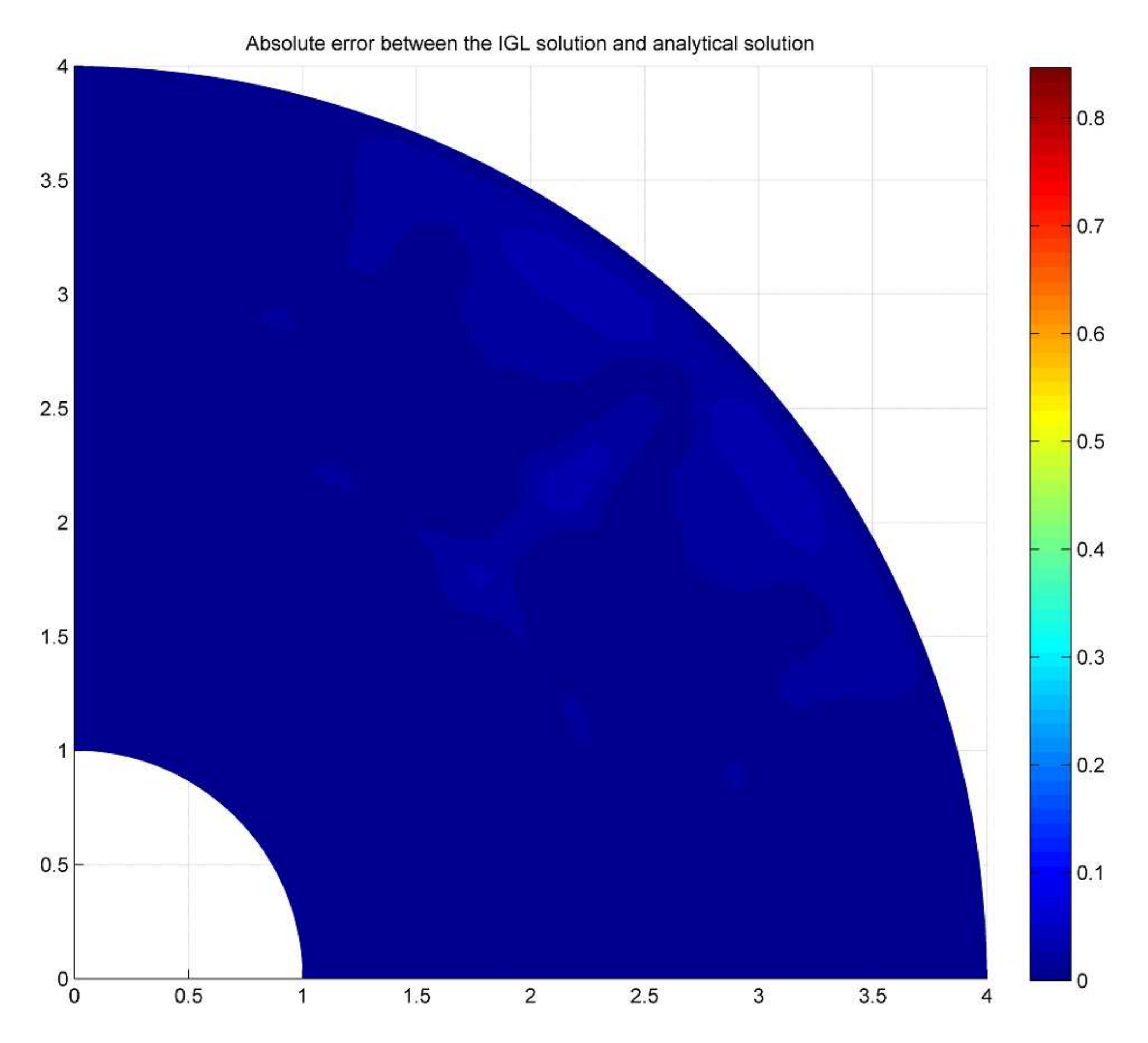}}
  \subfigure[Absolute error distribution of the IGA-C solution.]{
    \label{subfig:2d-igc-err}
  \includegraphics[width = 0.30\textwidth]
        {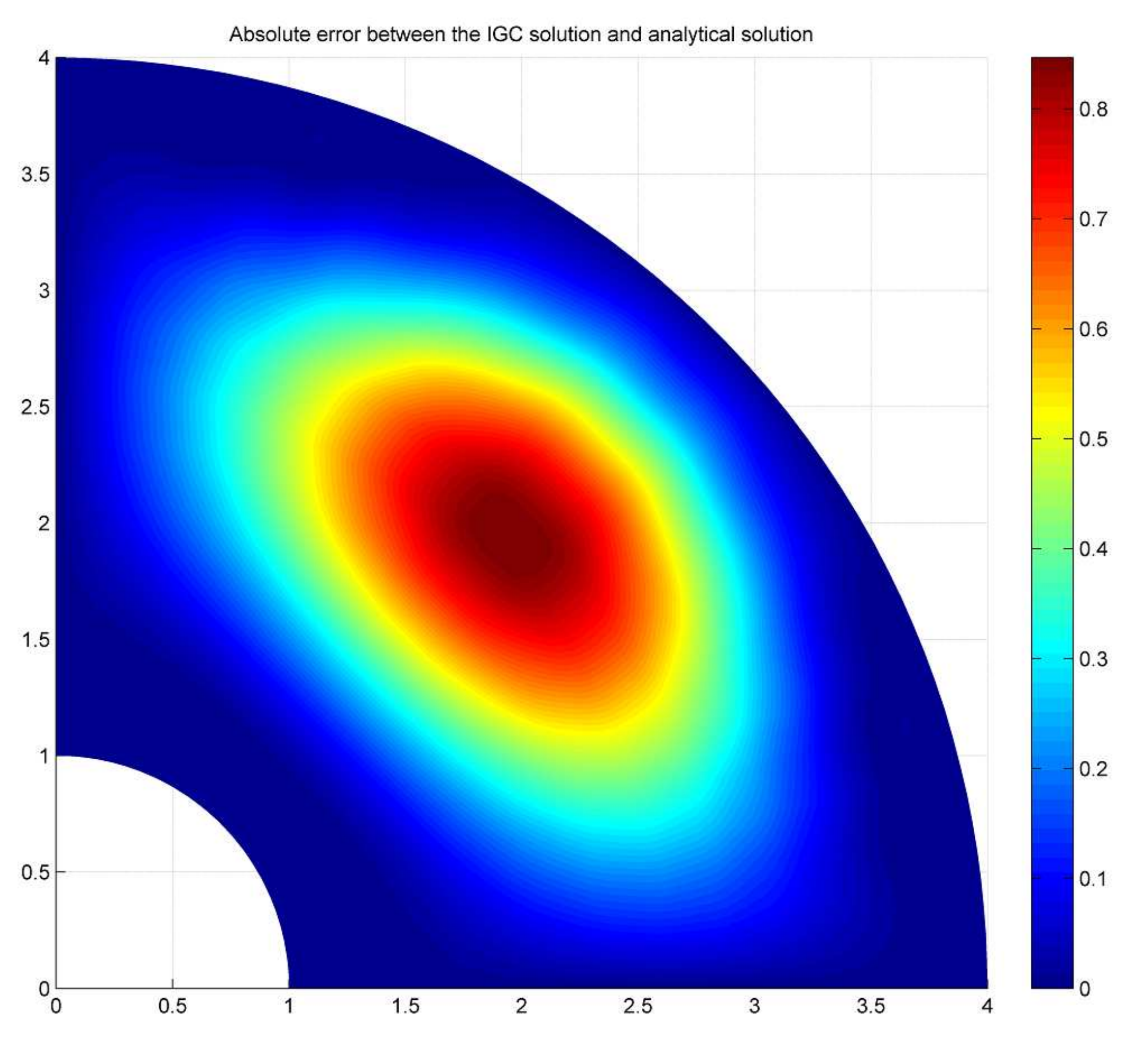}}
  \subfigure[Absolute error distribution of the IGA-SC solution.]{
    \label{subfig:2d-igsc-err}
  \includegraphics[width = 0.32\textwidth]
        {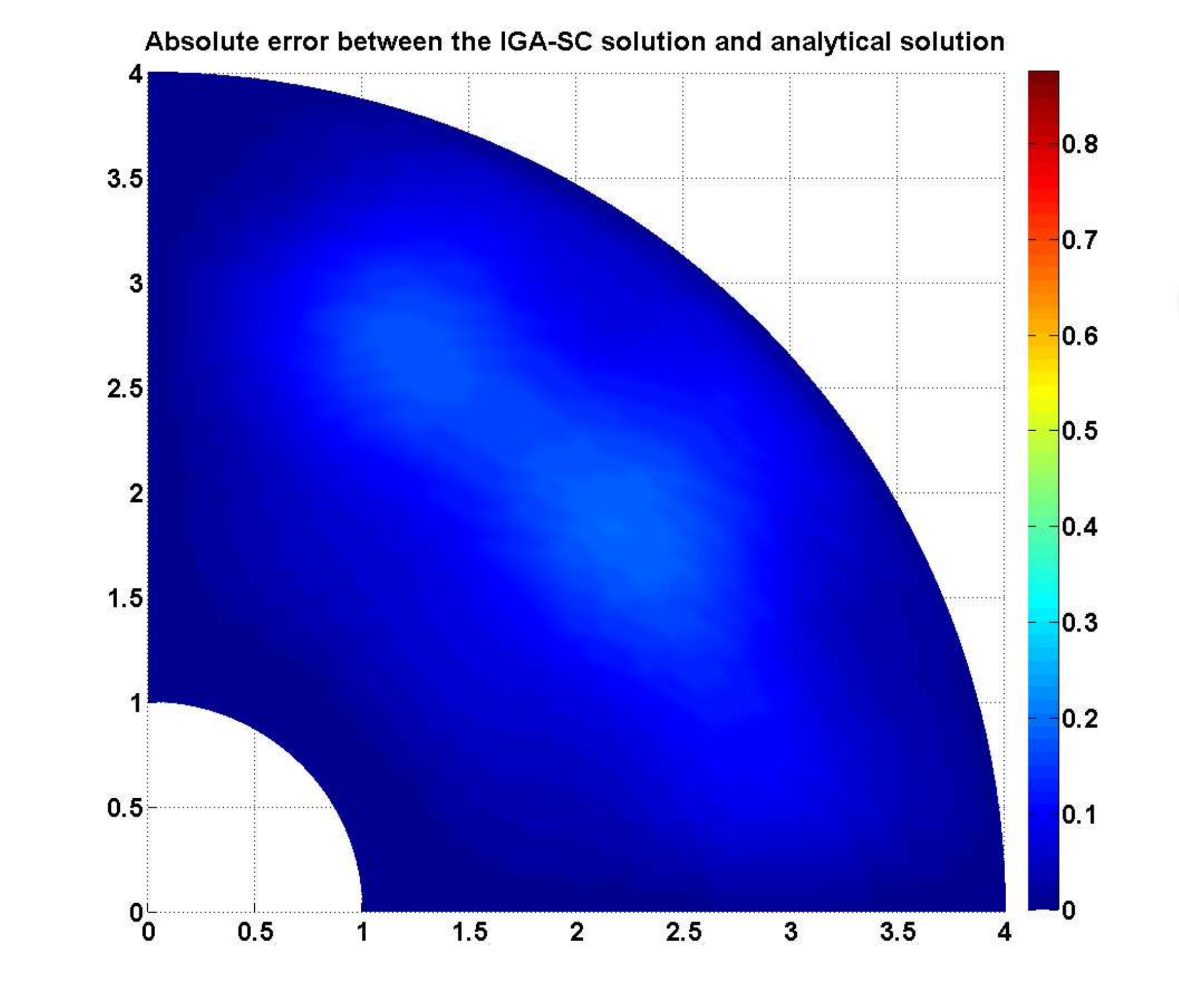}}
  \caption{Comparison of the analytical, IGA-L, IGA-C, and IGA-SC
    solutions of Eq.~\pref{eq:exmp_two_dim}.
    The relative error of the IGA-L solution is nearly two orders of magnitude less than the IGA-C solution.}
  \label{fig:two-dim-solutions}
\end{figure}

 \begin{figure}[!h]
    \centering
  \subfigure[]{
    \label{subfig:2d_iga_l}
    \includegraphics[width = 0.46\textwidth]{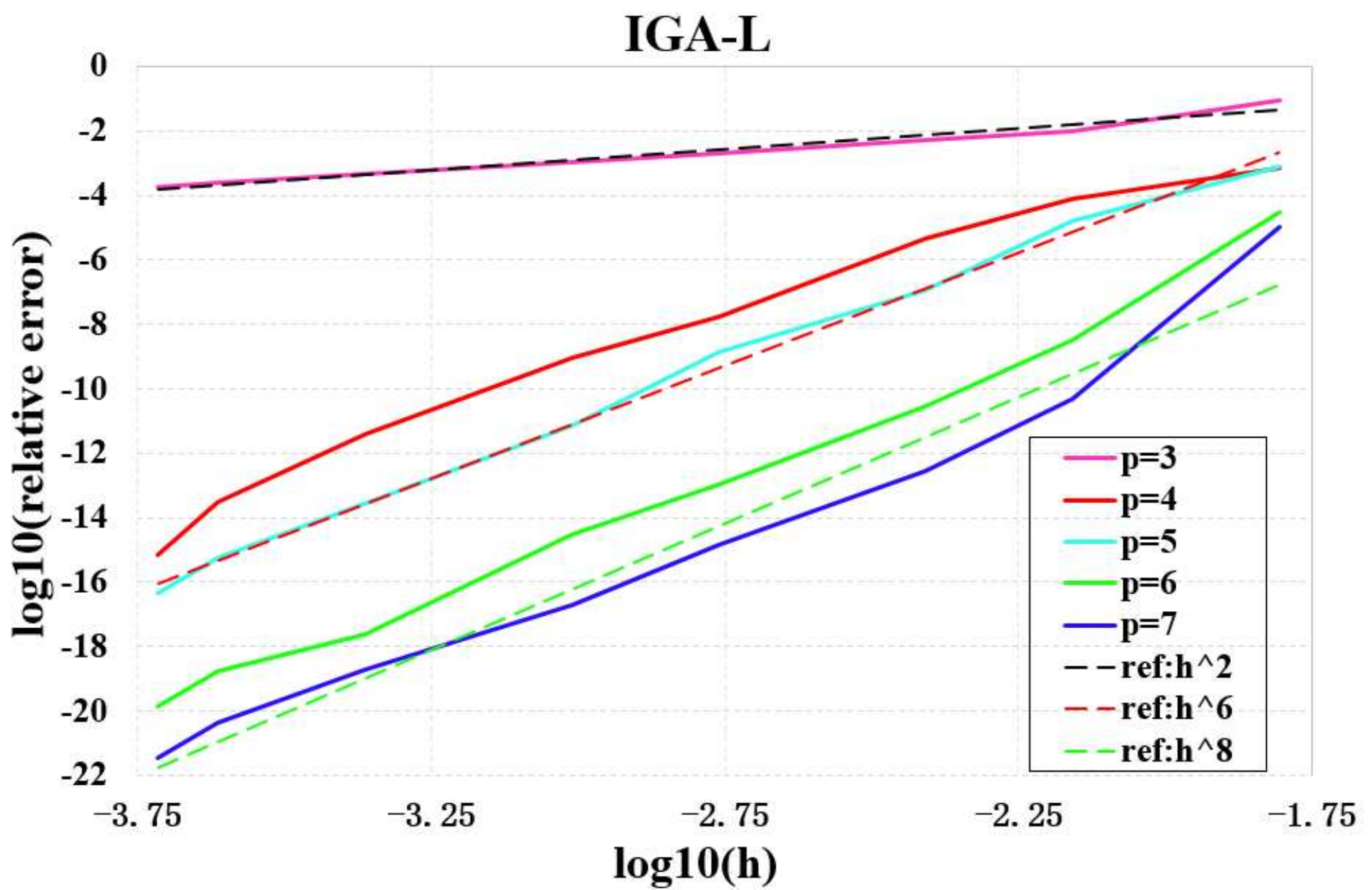}}
  \subfigure[]{\label{subfig:2d_iga_c}
    \includegraphics[width = 0.46\textwidth]{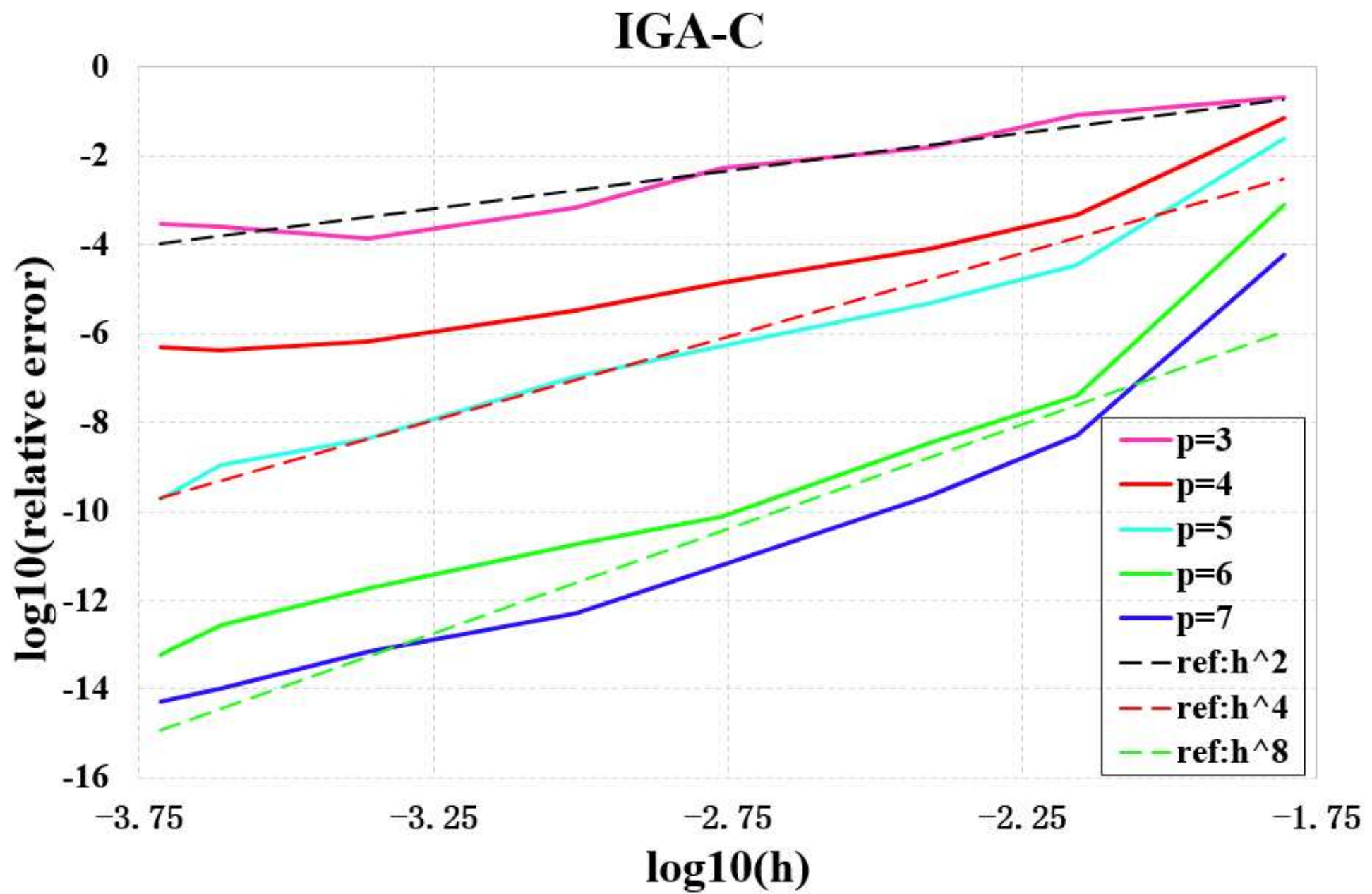}}
  \subfigure[]{\label{subfig:2d_iga_sc}
    \includegraphics[width = 0.46\textwidth]{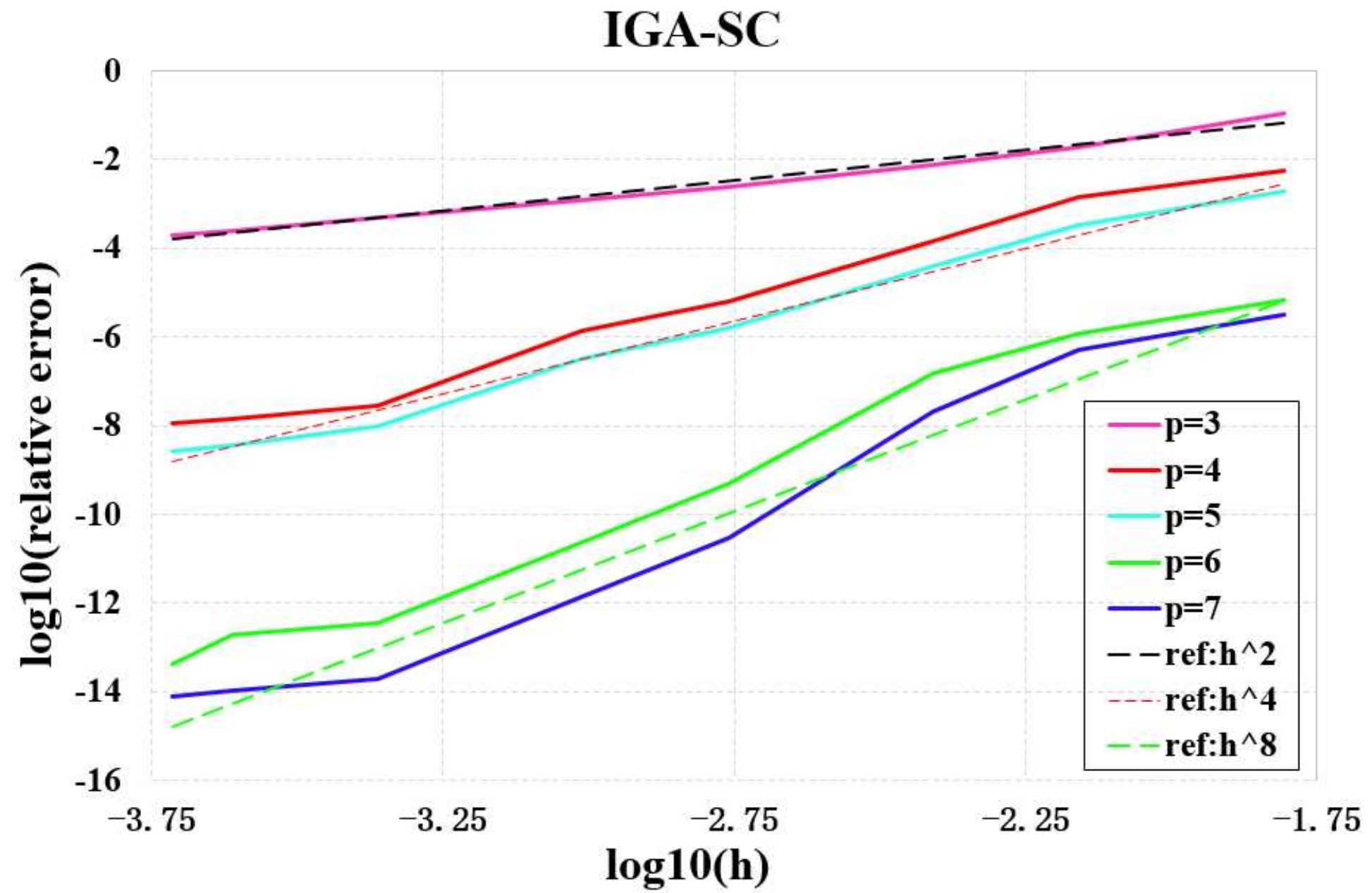}}
  \subfigure[]{\label{subfig:2d_time_err}
    \includegraphics[width = 0.46\textwidth]{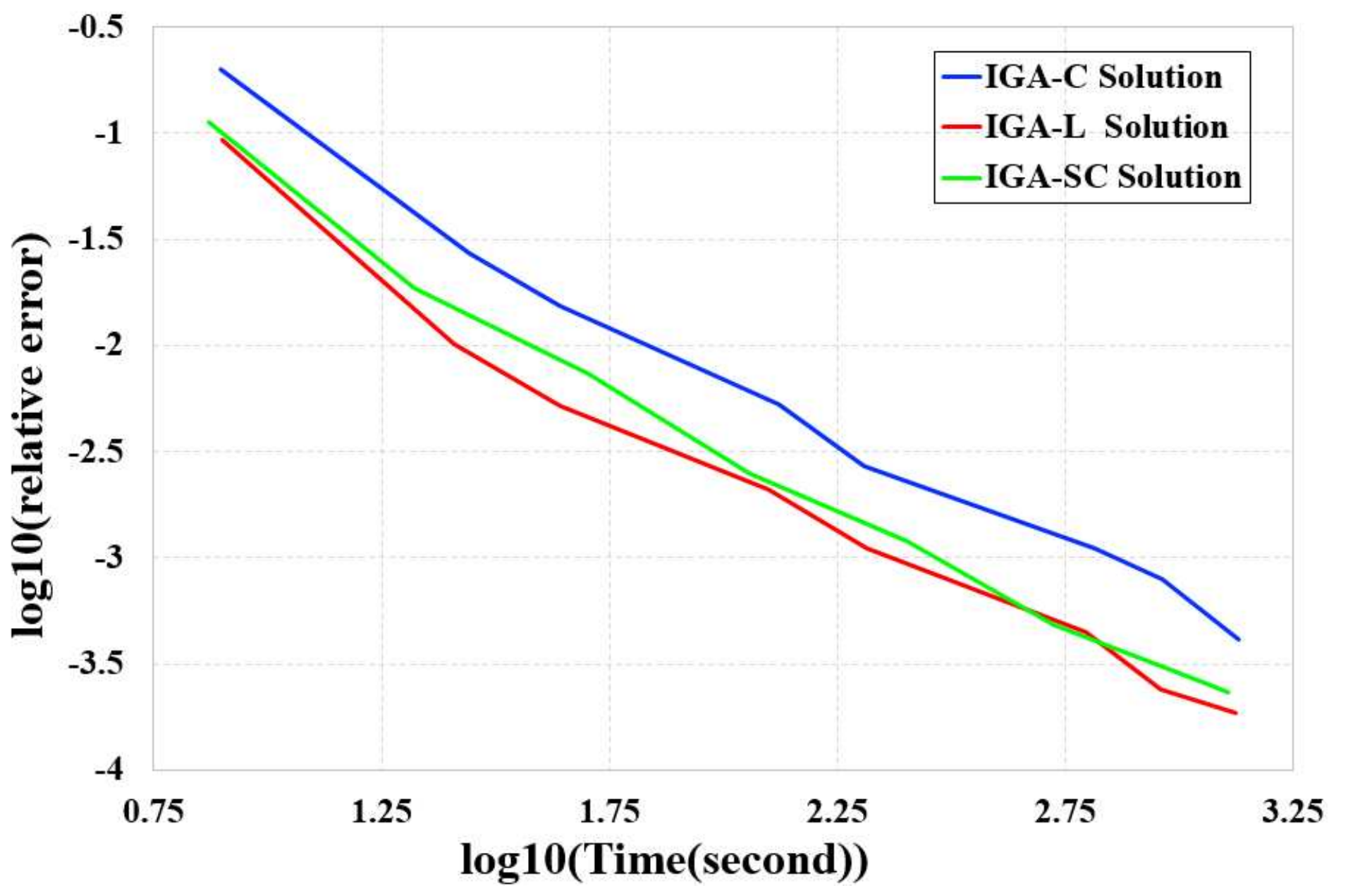}}
  \caption{Numerical results for the 2D source
        problem~\pref{eq:exmp_two_dim}.
        Diagrams of $\log_{10}(h)$ v.s. $\log_{10}$(relative error) for IGA-L (a), IGA-C (b), and IGA-SC (c), respectively. And, diagram of $\log_{10}$(Time) v.s. $\log_{10}$(relative error) (d).}
  \label{fig:two_dim_example}
\end{figure}

\textbf{Example II:}
 source problem in the
    two-dimensional domain $\Omega$,
    \begin{equation}
        \label{eq:exmp_two_dim}
        \begin{cases}
        &-\Delta T + T = f,\  (x,y) \in \Omega \\
        &T|_{\partial {\Omega}} = 0,
        \end{cases}
    \end{equation}
    where $\Omega$ is a quarter of an annulus,
    which can be exactly represented by a cubic NURBS
    patch with $4 \times 4$ control points,
    as presented in Appendix A2,
    and where
    \begin{equation*}
        \begin{aligned}
        f = (3x^{4} - 67x^{2} - 67y^{2} + 3y^{4} + 6x^{2}y^{2} + 116)\sin(x)\sin(y) \\
        + (68x - 8x^{3} - 8xy^{2})\cos(x)\sin(y) \\
        + (68y - 8y^{3}-8yx^{2})\cos(y)\sin(x).
        \end{aligned}
    \end{equation*}
 The analytical solution of the source problem~\pref{eq:exmp_two_dim}
    is
    \begin{displaymath}
        T = (x^2 + y^2 - 1)(x^2 + y^2 - 16)\sin(x)\sin(y).
    \end{displaymath}

 Fig.~\ref{fig:two-dim-solutions} presents numerical solutions of the
    two-dimensional source problem~\pref{eq:exmp_two_dim}, as
    generated by the IGA-L, IGA-C and IGA-SC methods.
 To produce the numerical solutions,
    we uniformly inserted $11$ knots along the $u-$ and $v-$ directions,
    respectively, to the cubic NURBS patch presented in Appendix A2,
    resulting in a cubic NURBS patch with $15 \times 15$ control
    points.
 With $20 \times 20$ Greville collocation points,
    the IGA-L method was employed to solve Eq.~\pref{eq:exmp_two_dim}.
 The relative error of the IGA-L solution (see
    Fig.~\ref{subfig:2d-igl-solution}) is $3.84 \times 10^{-4}$,
    and Fig.~\ref{subfig:2d-igl-err} illustrates the absolute error
    distribution of the IGA-L solution.
 Moreover, the source problem~\pref{eq:exmp_two_dim} was also solved
    by the IGA-C method using the same NURBS patch of $15 \times 15$ control
    points (see Fig.~\ref{subfig:2d-igc-solution})
    with Greville collocation points.
 The relative error of the IGA-C solution is $1.59 \times 10^{-2}$,
    and its absolute error distribution is illustrated in
    Fig.~\ref{subfig:2d-igc-err}.
 In this example,
    the relative error of the IGA-L solution is two orders of
    magnitude less than that of the IGA-C solution.
 In addition, we employed the IGA-SC method to solve the two-dimensional
    source problem~\pref{eq:exmp_two_dim},
    with the same bi-cubic NURBS patch of $15 \times 15$ control
    points (see Fig.~\ref{subfig:2d-igc-solution}),
    and $24 \times 24$ collocation points.
 Note that, even the number of collocation points ($24 \times 24$) for the
    IGA-SC method is larger than that for the IGA-L method ($20 \times 20$),
    the relative error of the IGA-L solution ($3.84 \times 10^{-4}$) is still less than that of the IGA-SC solution ($4.62 \times 10^{-4}$).
 Fig.~\ref{subfig:2d-igsc-err} demonstrates the absolute error distribution
    of the IGA-SC solution.

 Diagrams of $\log_{10}(h)$ v.s. $\log_{10}$(relative error) for
    IGA-L, IGA-C, and IGA-SC methods are illustrated in Figs.~\ref{subfig:2d_iga_l}-~\ref{subfig:2d_iga_sc}, respectively.
 It should be pointed out that, while the convergence rate of the IGA-SC and
    IGA-C methods with degree $p = 4,5$ is $O(h^4)$,
    that of the IGA-L method with degree $p = 4,5$ reaches $O(h^6)$.
 Moreover, Fig.~\ref{subfig:2d_time_err} presents the diagrams of
    $\log_{10}$(time) v.s. $\log_{10}$(relative error) for the three methods.
 In these diagrams, the performance of the IGA-L and IGA-SC methods are
    comparable,
    both better than that of the IGA-C method.

\begin{figure}[!h]
\centering
    \subfigure[IGA-L solution of Eq.~\pref{eq:exmp_three_dim} with relative error $0.0232$.]{
    \label{subfig:3d-igl-solution}
    \includegraphics[width = 0.27\textwidth]
    {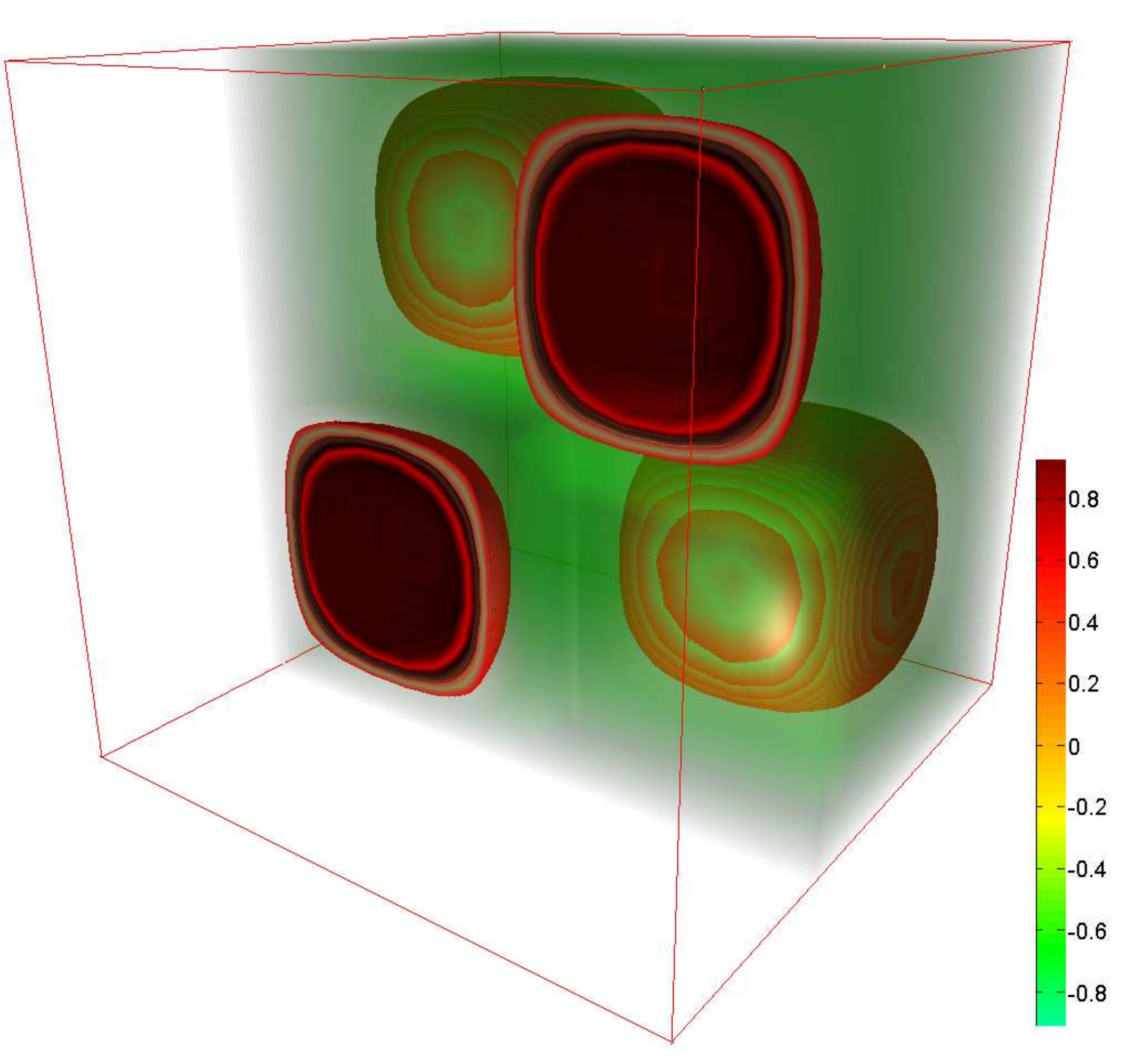}}
  \subfigure[IGA-C solution of Eq.~\pref{eq:exmp_three_dim} with relative error $0.1546$.]{
    \label{subfig:3d-igc-solution}
    \includegraphics[width = 0.27\textwidth]
    {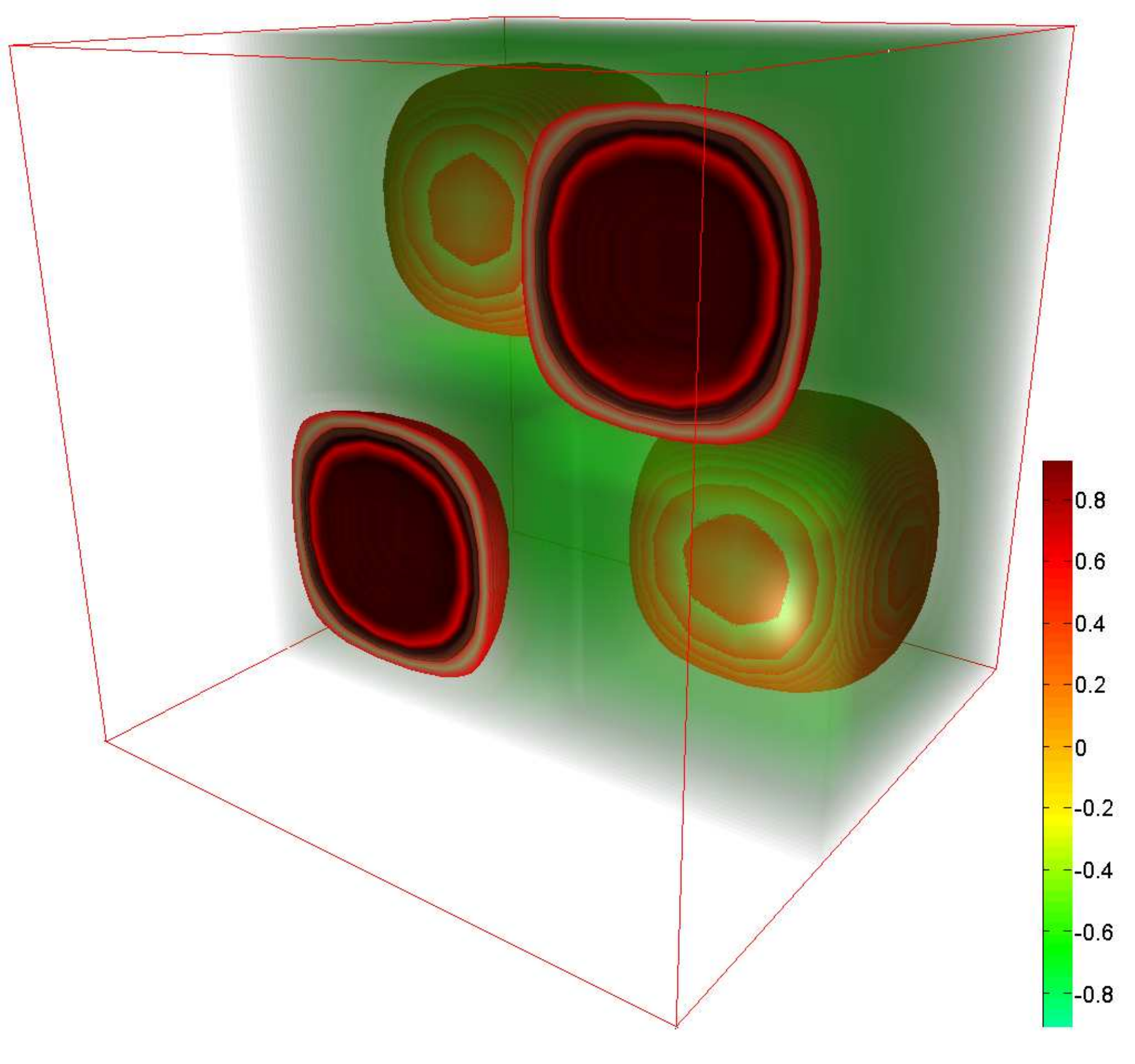}}
  \subfigure[IGA-SC solution of Eq.~\pref{eq:exmp_three_dim} with relative error $0.0347$.]{
    \label{subfig:3d-igsc-solution}
    \includegraphics[width = 0.27\textwidth]
    {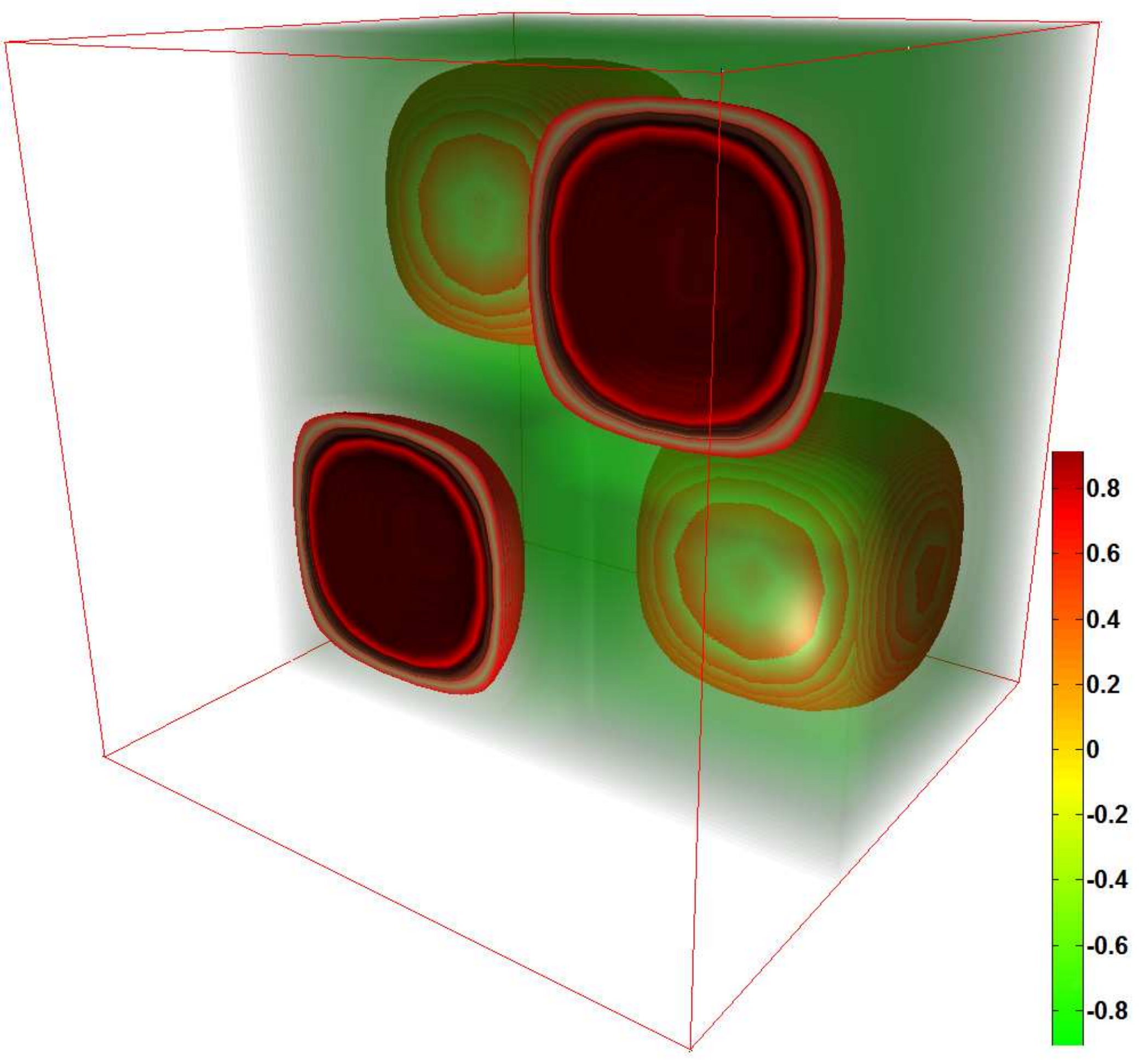}}
  \subfigure[Absolute error distribution of the IGA-L solution.]{
    \label{subfig:3d-igl-err}
  \includegraphics[width = 0.27\textwidth]
    {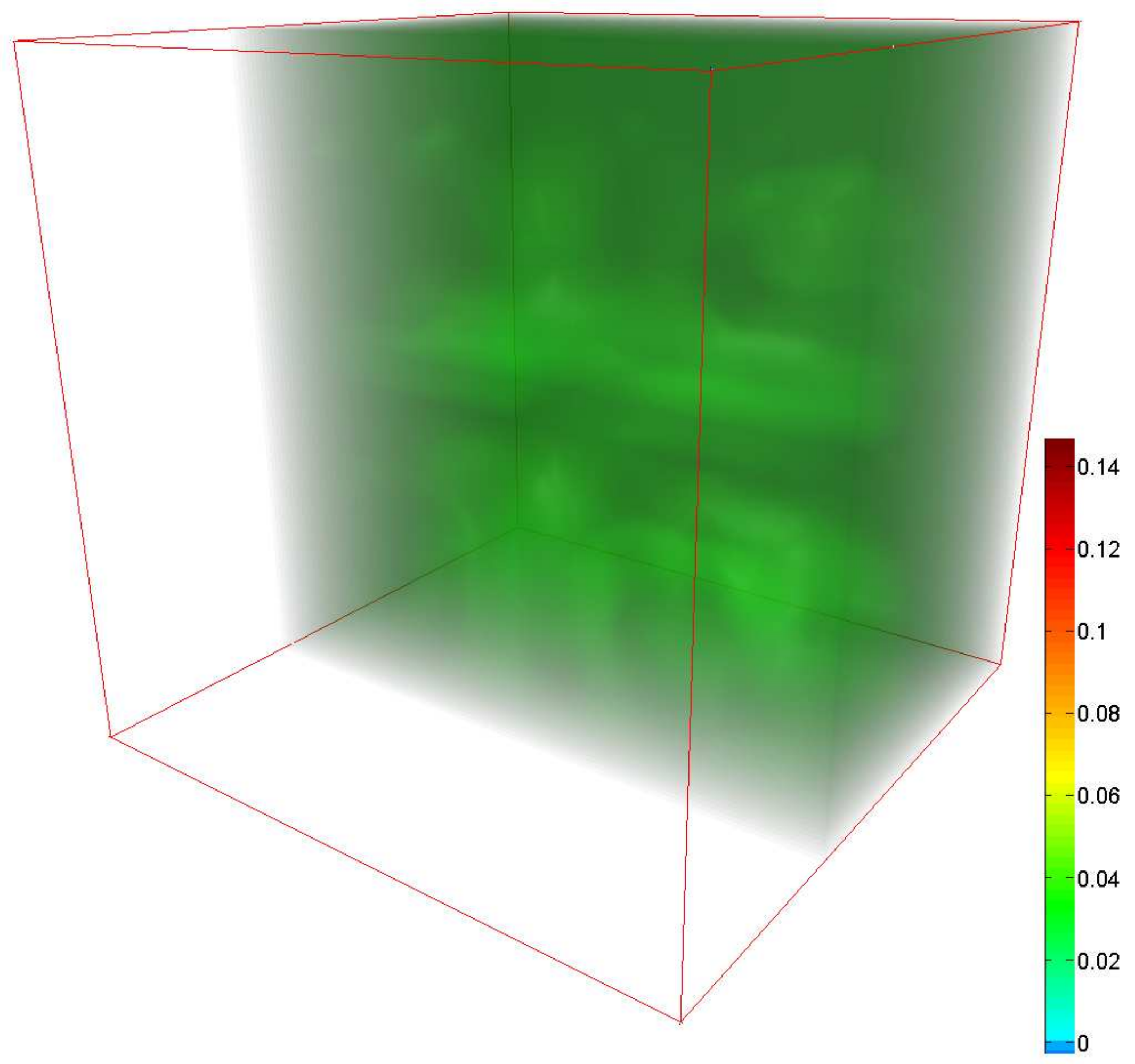}}
  \subfigure[Absolute error distribution of the IGA-C solution.]{
    \label{subfig:3d-igc-err}
    \includegraphics[width = 0.27\textwidth]
    {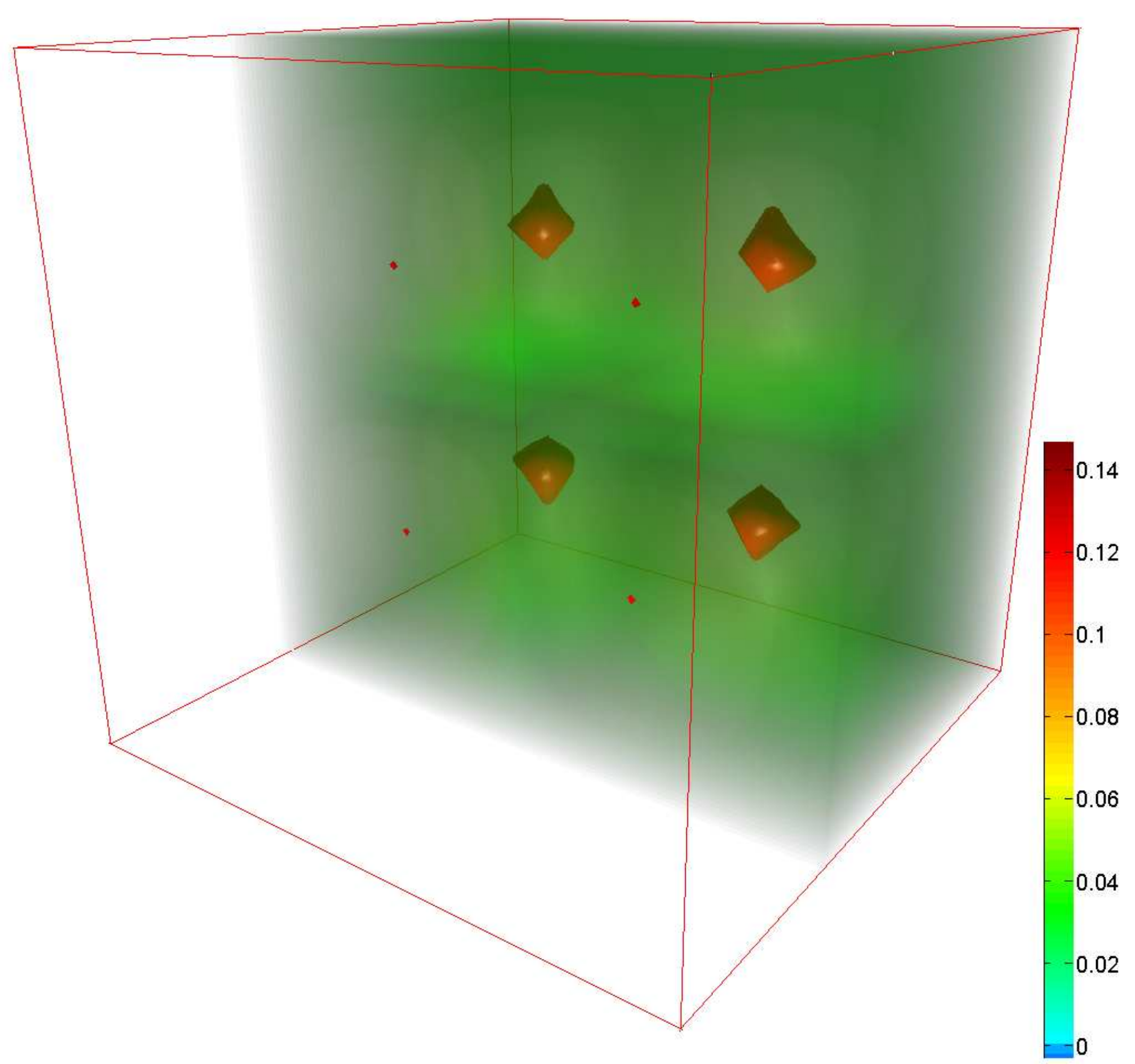}}
  \subfigure[Absolute error distribution of the IGA-SC solution.]{
    \label{subfig:3d-igsc-err}
    \includegraphics[width = 0.27\textwidth]
    {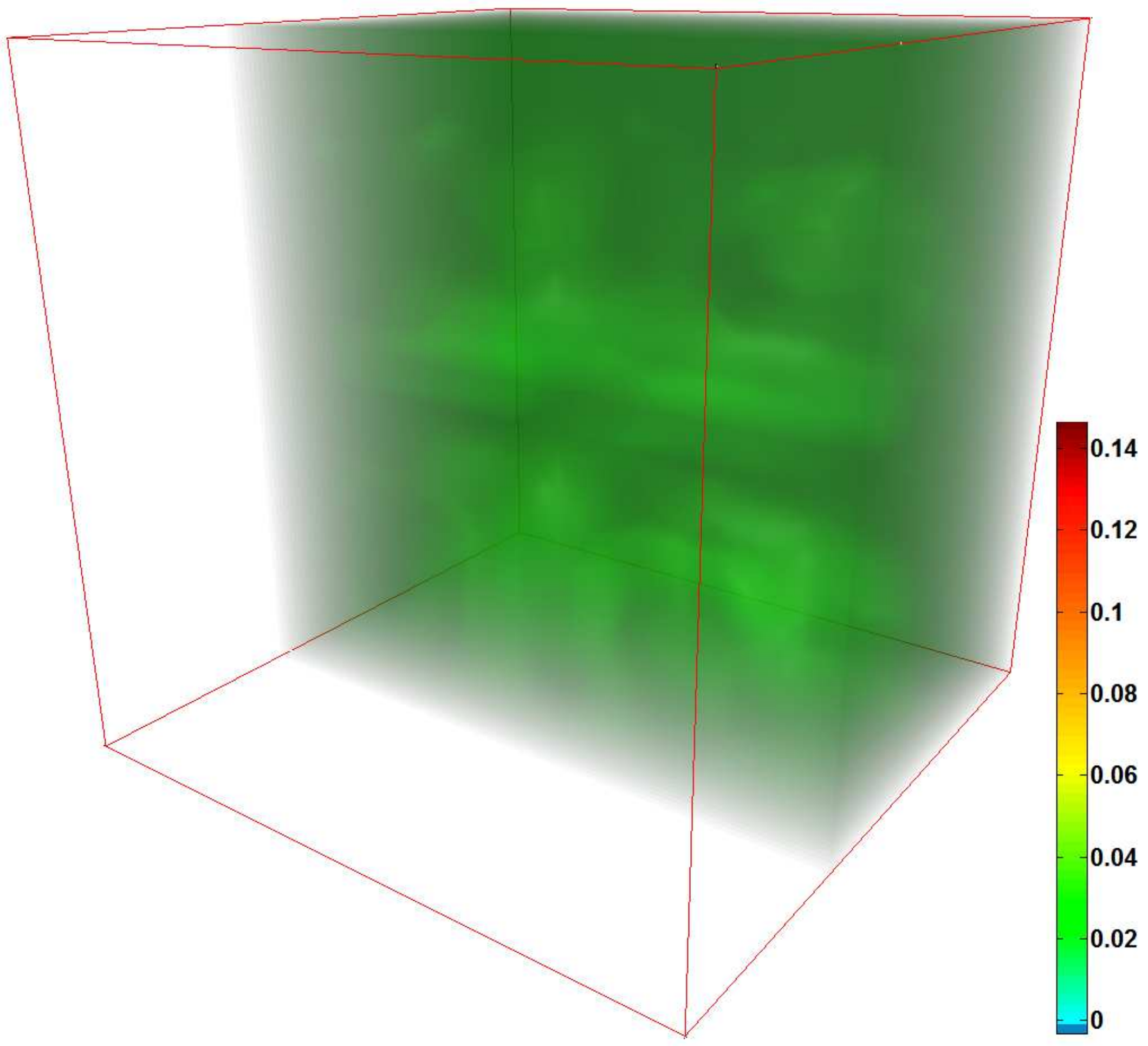}}
  \caption{Comparison of the IGA-L, IGA-C and IGA-SC solutions of Eq.~\pref{eq:exmp_three_dim}.
  The relative error of the IGA-L solution
    is nearly one order of magnitude less than that of the IGA-C solution.}
  \label{fig:three-dim-solution}
\end{figure}

 \begin{figure}[!h]
    \centering
  \subfigure[]{
    \label{subfig:3d-iga-l}
    \includegraphics[width = 0.46\textwidth]
        {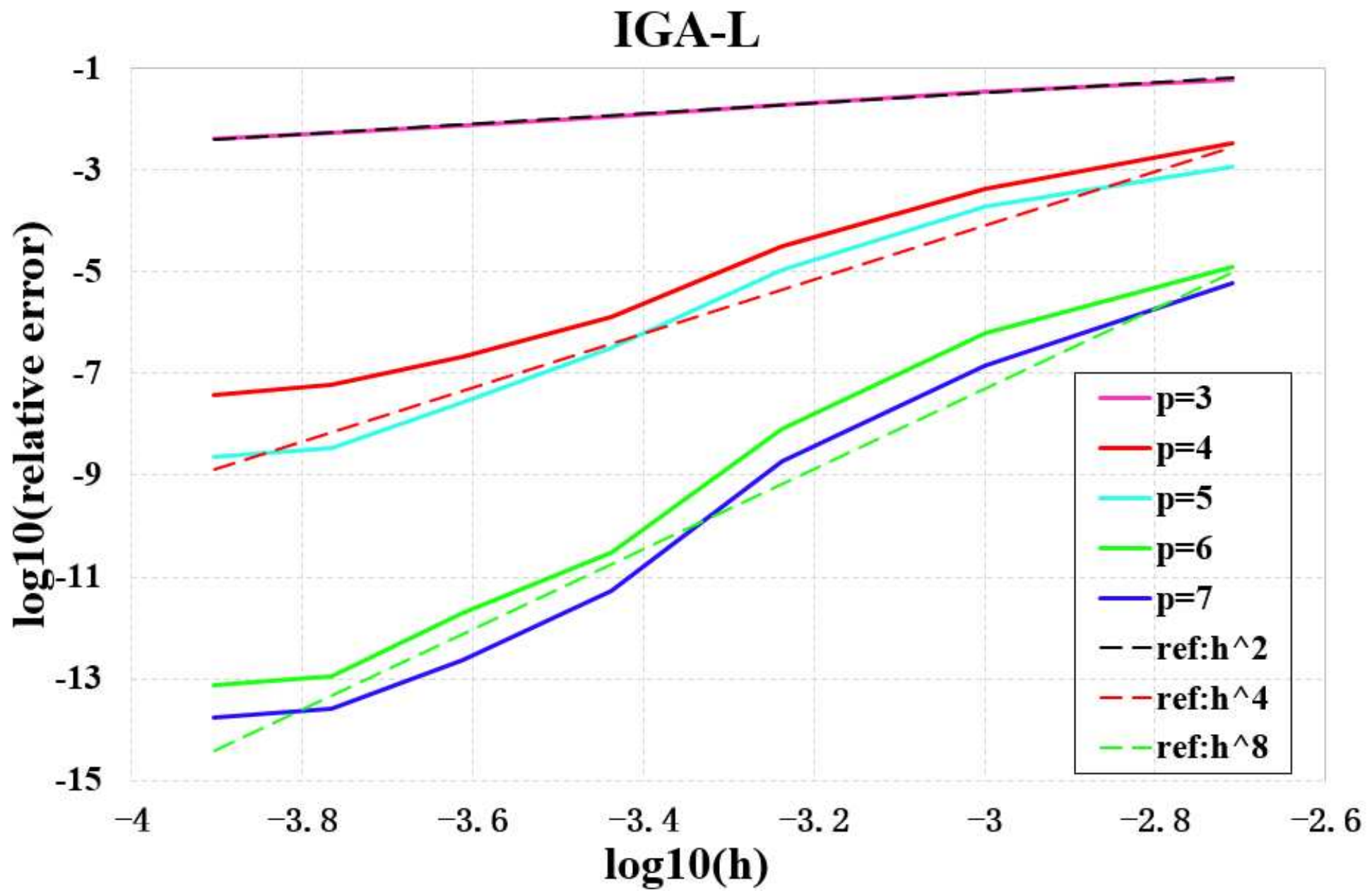}}
  \subfigure[]{\label{subfig:3d-iga-c}
    \includegraphics[width = 0.46\textwidth]
        {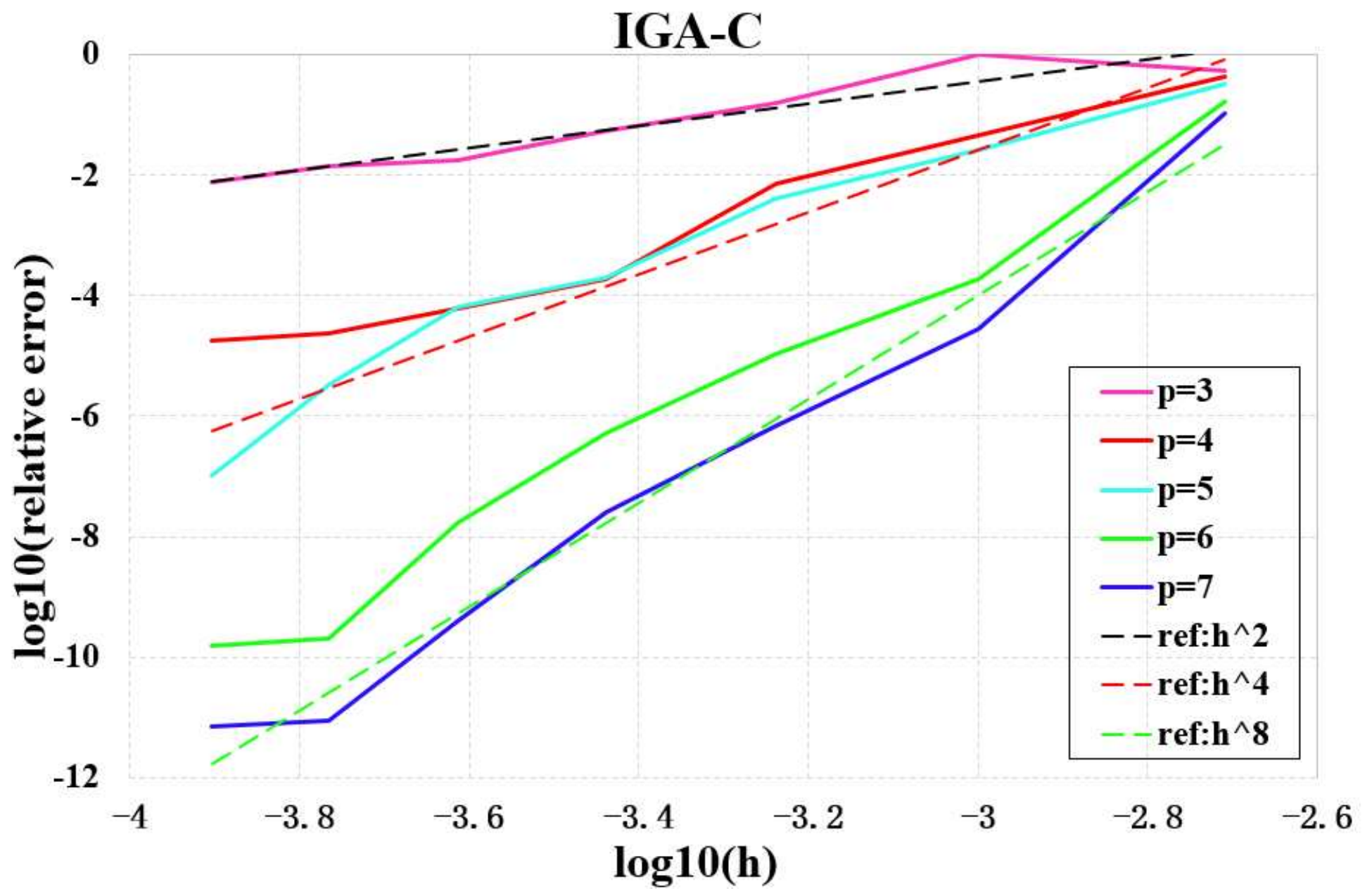}}
  \subfigure[]{\label{subfig:3d-iga-sc}
    \includegraphics[width = 0.46\textwidth]
        {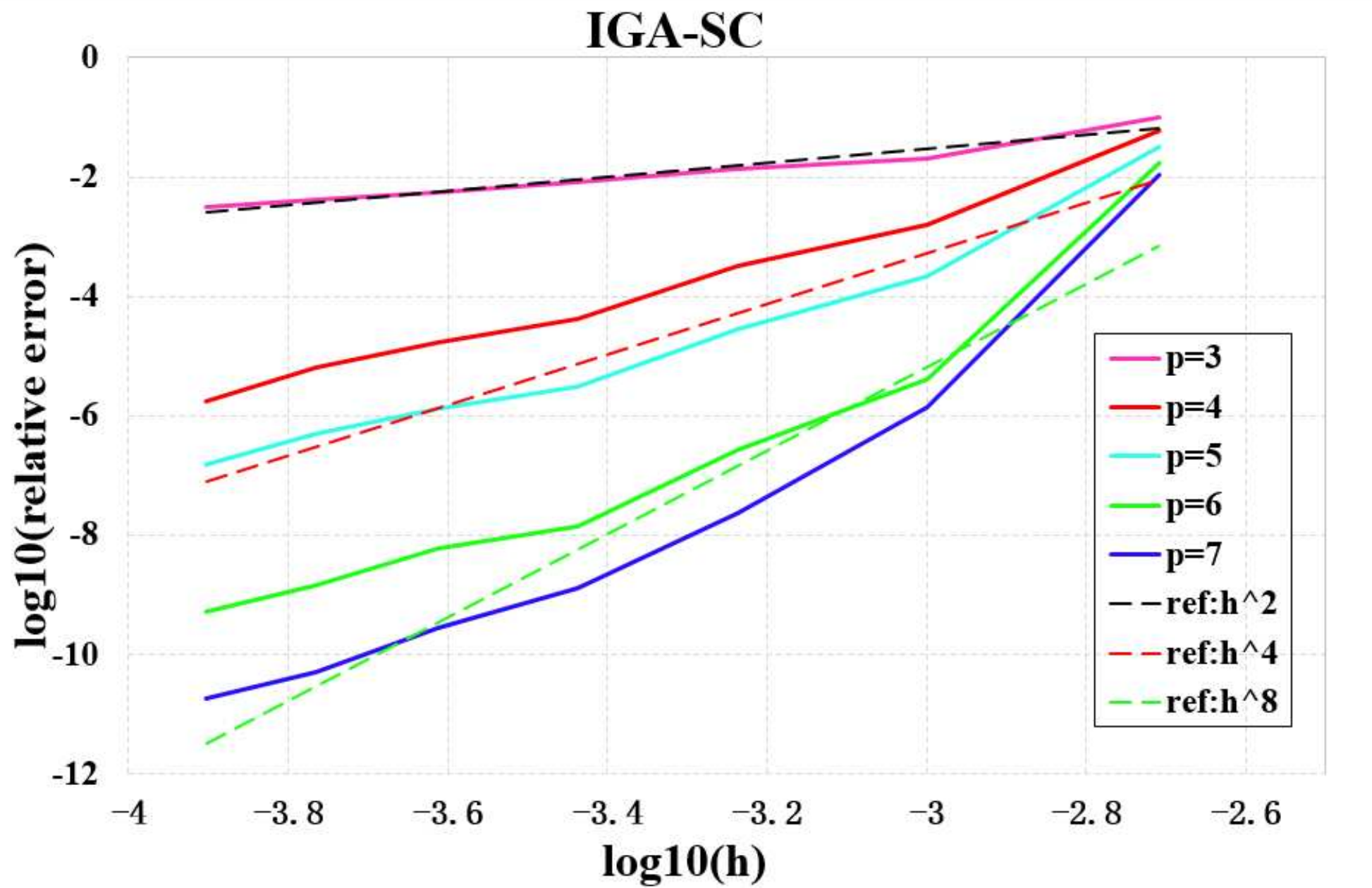}}
  \subfigure[]{\label{subfig:3d-time-err}
    \includegraphics[width = 0.46\textwidth, height = 0.30\textwidth]
        {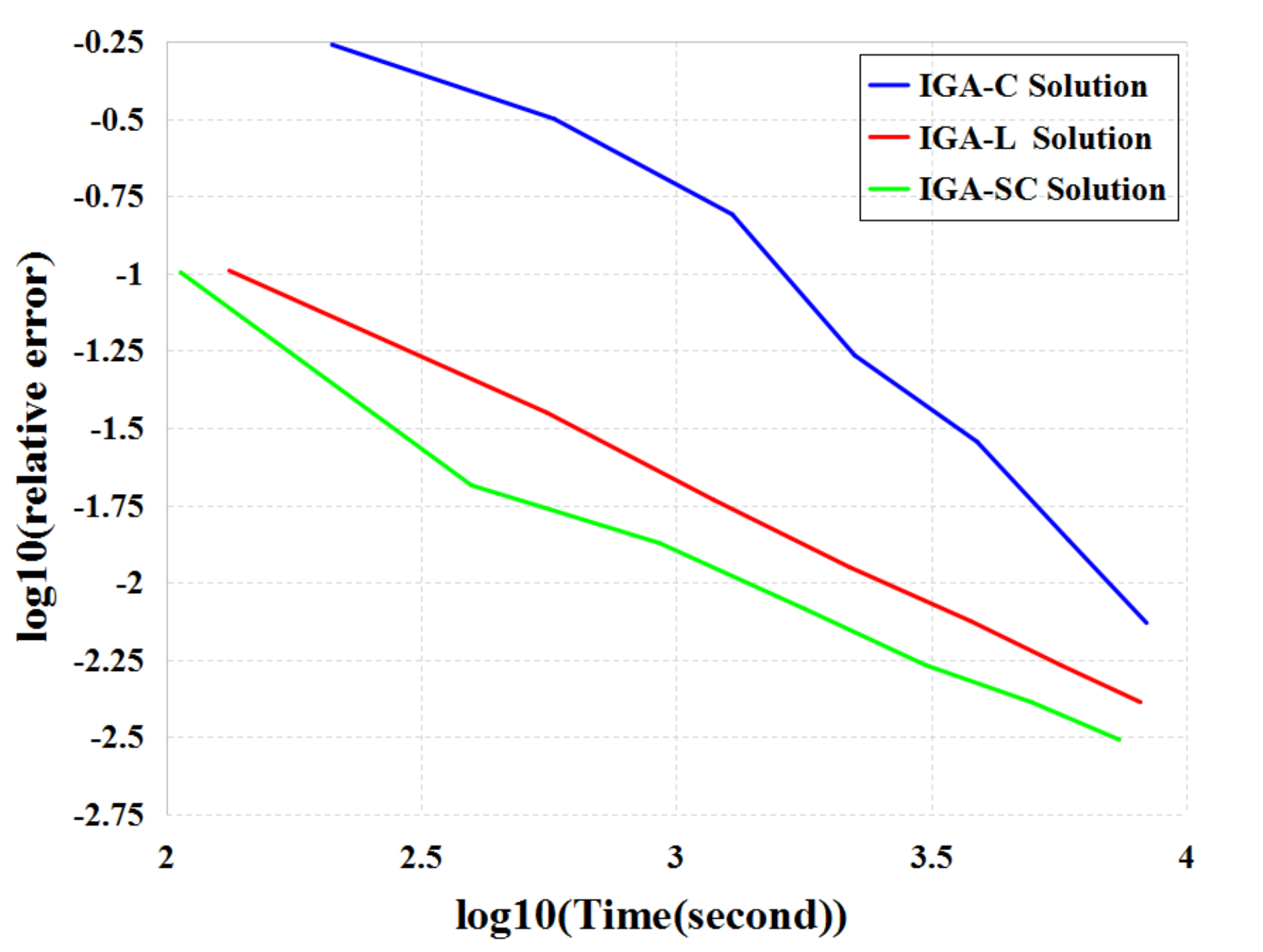}}
  \caption{Numerical results for the 3D source
            problem~\pref{eq:exmp_three_dim}.
  Diagrams of $\log_{10}(h)$ v.s. $\log_{10}$(relative error) for IGA-L (a), IGA-C (b), and IGA-SC (c), respectively. And, diagram of $\log_{10}$(Time) v.s. $\log_{10}$(relative error) (d).
}
  \label{fig:three_dim_example}
\end{figure}

 \textbf{Example III:}
 source problem defined on the three-dimensional cubic domain
    $\Omega = [0,1] \times [0,1] \times [0,1]$, i.e.,
    \begin{equation}
    \label{eq:exmp_three_dim}
    \begin{cases}
    & -\Delta T + T = f,\  (x,y,z) \in \Omega, \\
    & T|_{\partial {\Omega}} = 0,
    \end{cases}
    \end{equation}
    where
    \begin{equation*}
    f = (1 + 12\pi^{2})\sin(2 \pi x)\sin(2 \pi y)\sin(2 \pi z).
    \end{equation*}
 Its analytical solution is,
    \begin{equation*}
    T = \sin(2 \pi x)\sin(2 \pi y)\sin(2 \pi z).
    \end{equation*}
 The three-dimensional physical domain $\Omega$ is modeled as a cubic
    B-spline solid with $4 \times 4 \times 4$ control points,
    as listed in Appendix A3.

 In Fig.~\ref{fig:three-dim-solution}, the IGA-L, IGA-C, and IGA-SC solutions
    for Eq.~\pref{eq:exmp_three_dim} are illustrated,
    where the solutions were generated with tri-cubic B-spline solid.
 Specifically, to get the IGA-L solution,
    $7 \times 7 \times 7$ control points and $10 \times 10 \times 10$ Greville collocation points were employed,
    and the relative error is $0.0232$ (see Fig.~\ref{subfig:3d-igl-solution}).
 On the other hand, the relative error for the IGA-C solution with
    $7 \times 7 \times 7$ control points and Greville collocation points is $0.1456$ (Fig.~\ref{subfig:3d-igc-solution}).
 In this example,
    the relative error of the IGA-L solution is one order of magnitude smaller than that of the IGA-C solution.
 Moreover, in Fig.~\ref{subfig:3d-igsc-solution},
    $10 \times 10 \times 10$ control points and $14 \times 14 \times 14$ collocation points were used to generate the IGA-SC solution with relative error $0.0347$.
 Similar as the two-dimensional case, the relative error of the IGA-SC
    solution is larger than that of the IGA-L solution,
    though the numbers of control points and collocation points of IGA-L method are both less than those of the IGA-SC method.
 Additionally, Figs.~\ref{subfig:3d-igl-err}-~\ref{subfig:3d-igsc-err}
    present the absolute error distribution for the IGA-L, IGA-C,
    and IGA-SC solutions, respectively.

 Furthermore, diagrams of the numerical results for
    Eq.~\pref{eq:exmp_three_dim} are
    illustrated in Fig.~\ref{fig:three_dim_example}.
 Specifically, diagrams of $\log_{10}(h)$ v.s. $\log_{10}$(relative error)
    for IGA-L, IGA-C, and IGA-SC are illustrated in Figs~\ref{subfig:3d-iga-l}-~\ref{subfig:3d-iga-sc}, respectively.
 From these diagrams, it can be seen that,
    the convergence rates of the three methods are all $O(h^2)$ for $p=3$,
    $O(h^4)$ for $p=4,5$, and $O(h^8)$ for $p=6,7$.
 The diagrams of $\log_{10}$(Time) v.s. $\log_{10}$(relative error) for the
    the three methods are presented in Fig.~\ref{subfig:3d-time-err},
    where the performance of IGA-SC method is the best.

\begin{figure}[!htb]
  \centering
  \includegraphics[width=0.4\textwidth, bb = 0 0 1082 603]
    {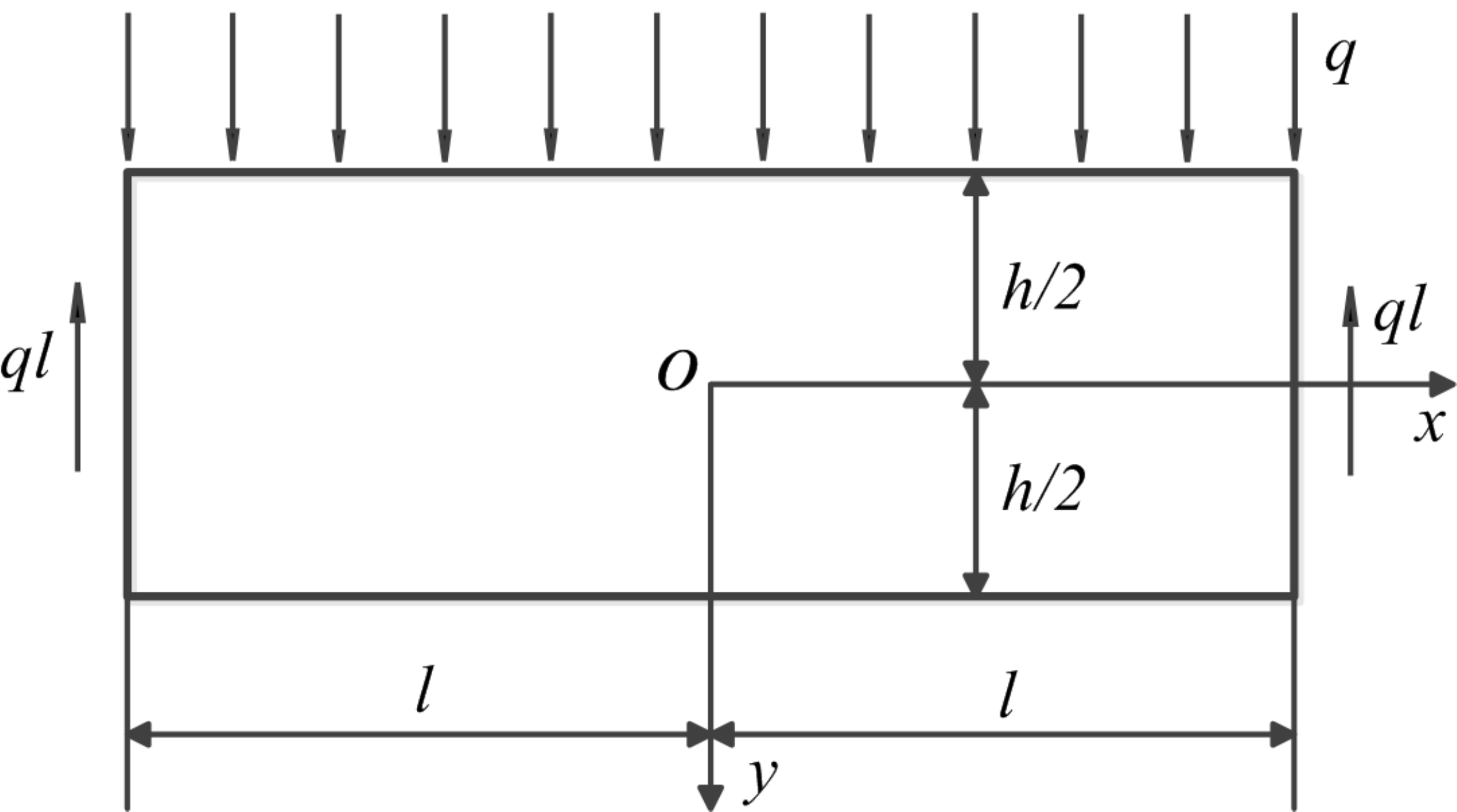}
  \caption{The simply supported beam.}
  \label{fig:beam}
\end{figure}

\begin{figure}[!htb]
\centering
  \subfigure[]{
    \label{fig:ana_delta_x}
    \includegraphics[width = 0.30\textwidth, bb = 0 0 1289 681]
    {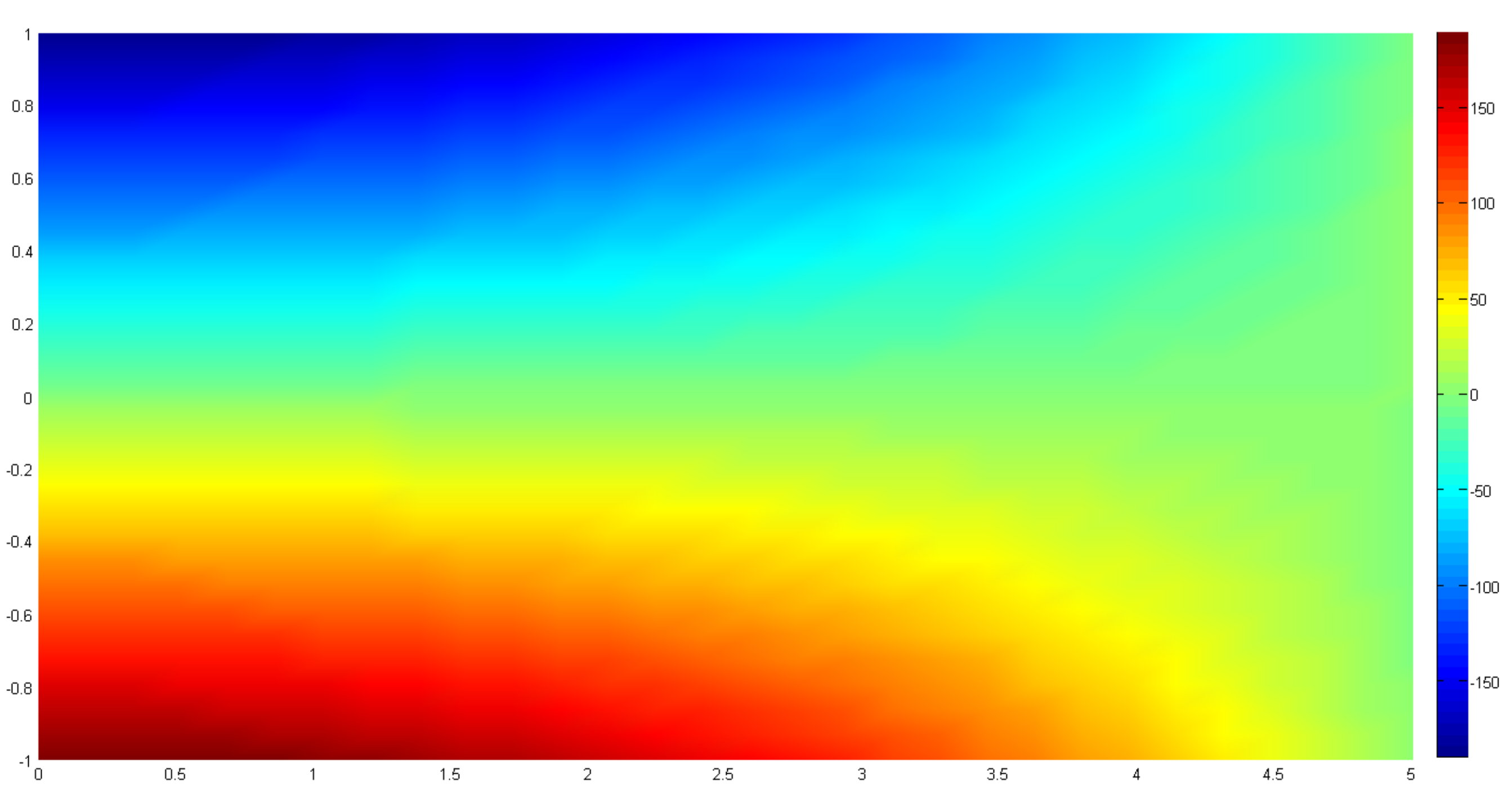}}
  \subfigure[]{
    \label{fig:ana_delta_y}
    \includegraphics[width = 0.30\textwidth, bb = 0 0 1293 669]
    {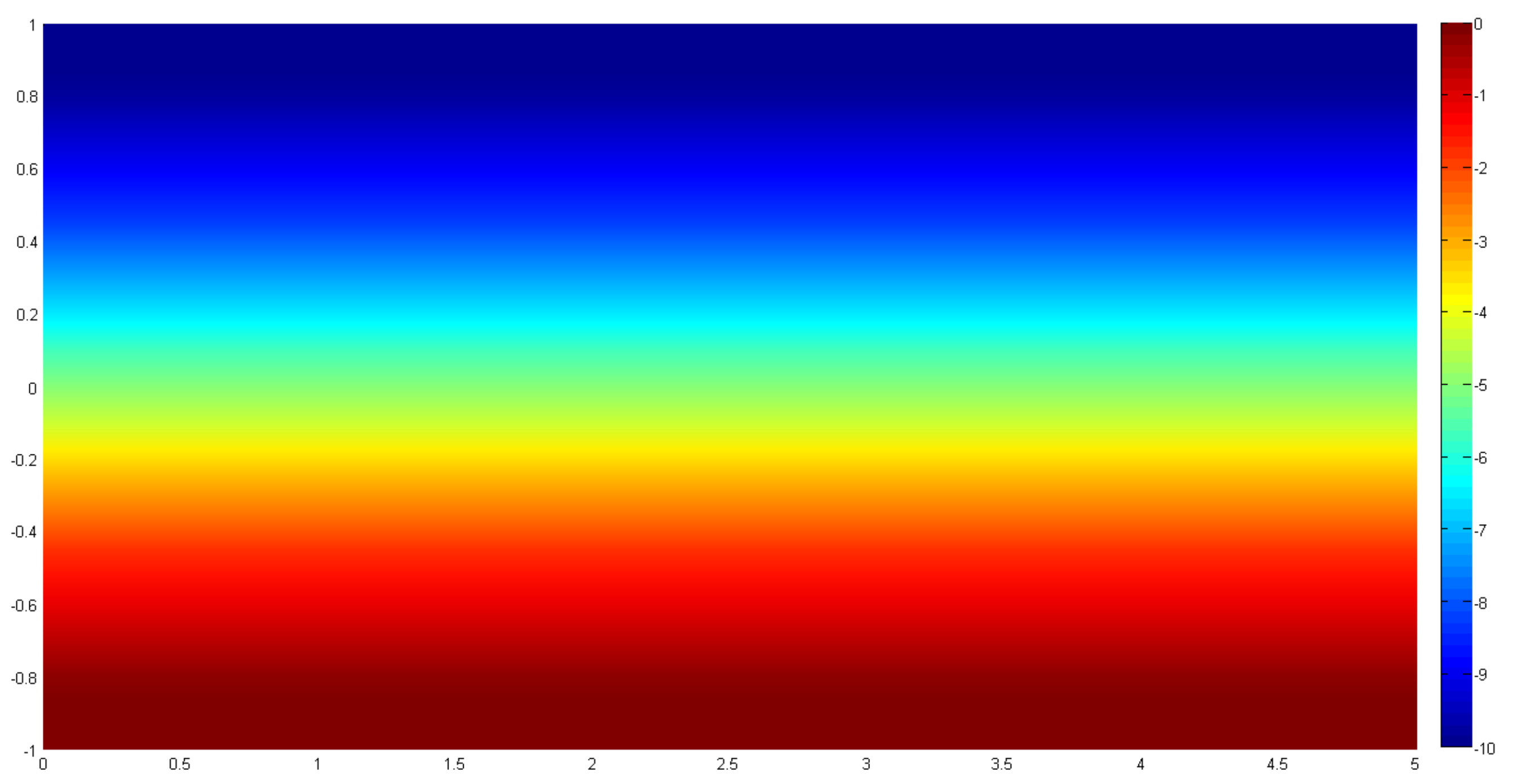}}
  \subfigure[]{
    \label{fig:ana_tau_xy}
    \includegraphics[width = 0.30\textwidth, bb = 0 0 1291 677]
    {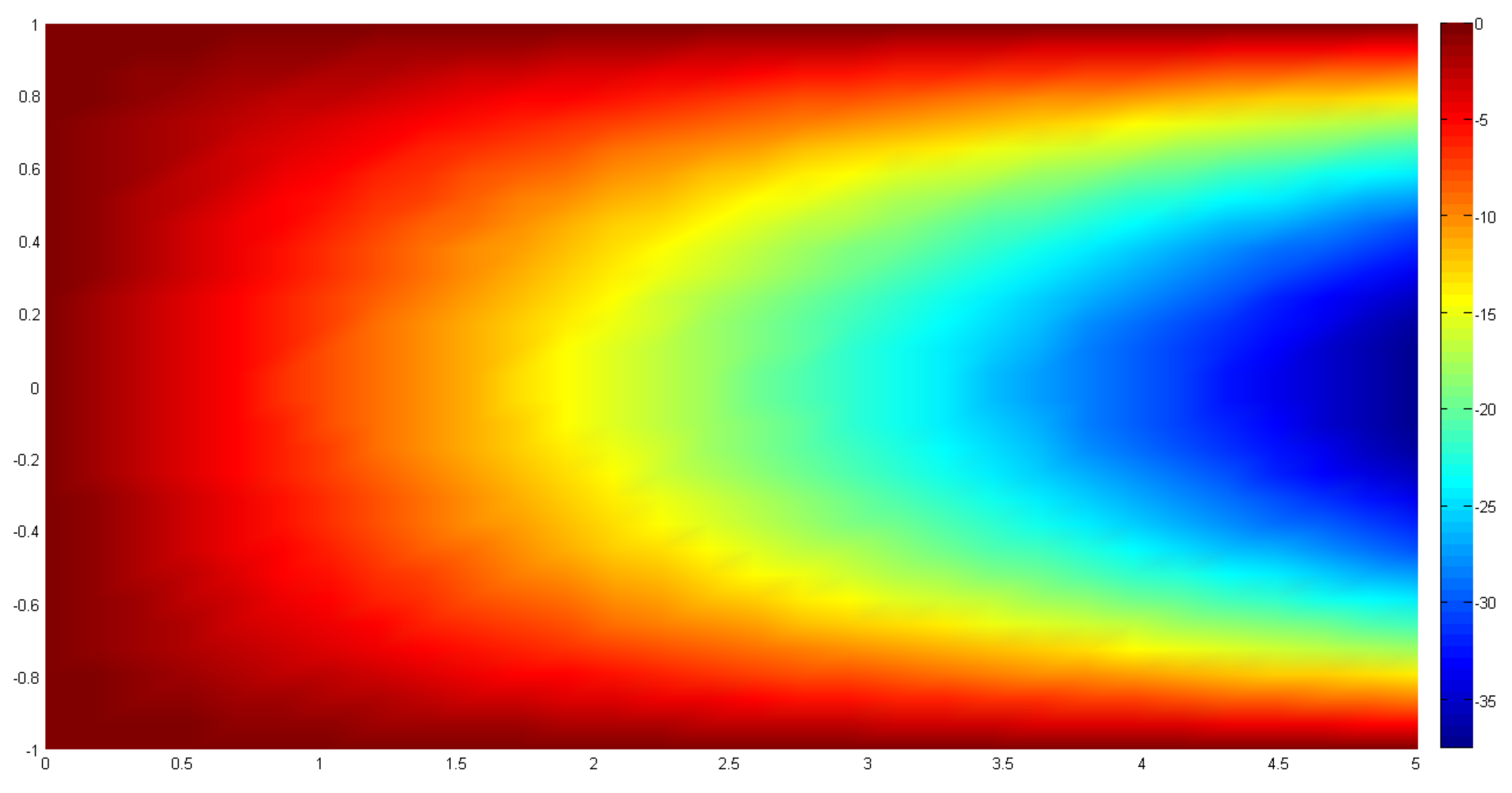}}
  \subfigure[]{
    \label{fig:num_delta_x}
    \includegraphics[width = 0.30\textwidth, bb = 0 0 1301 669]
    {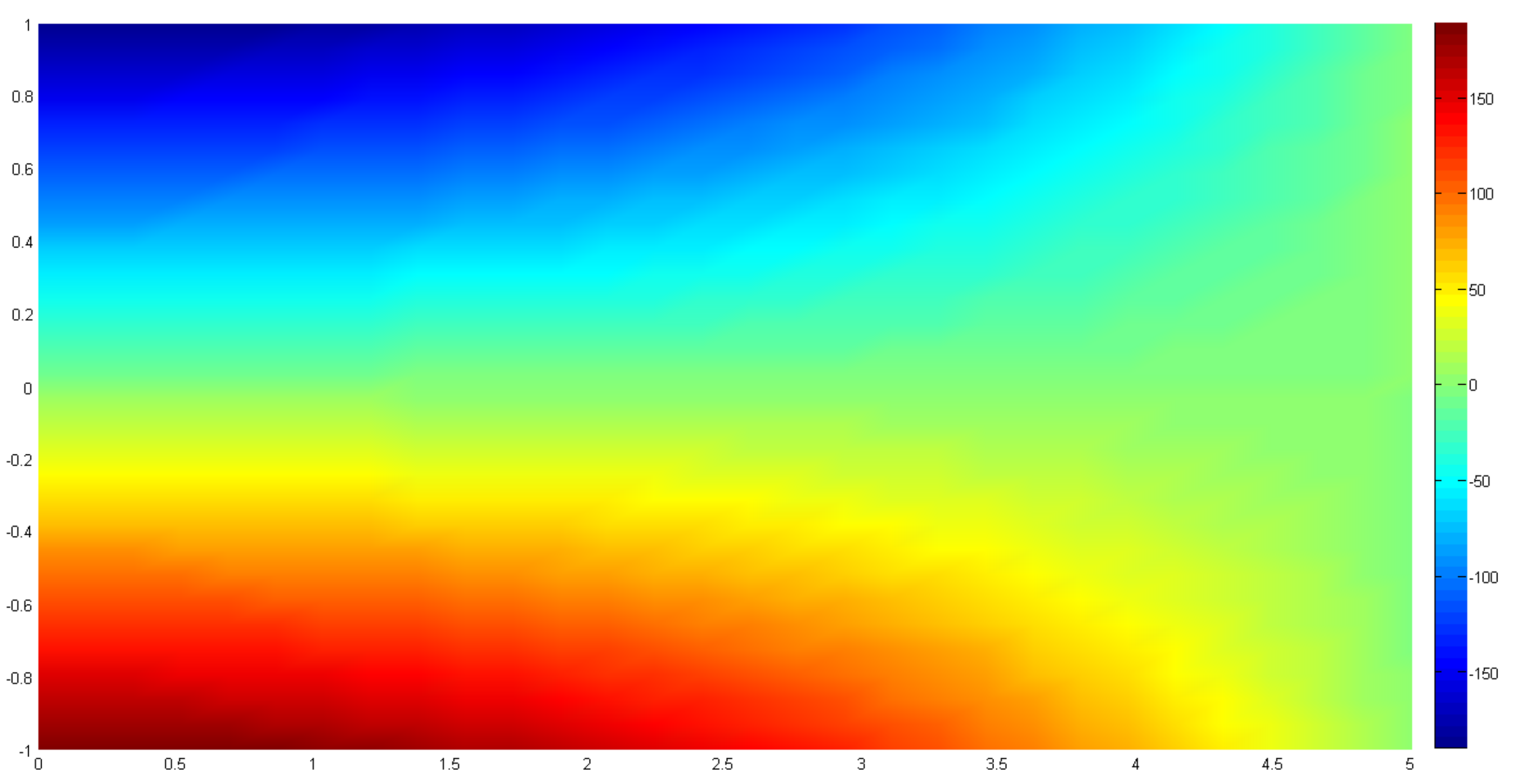}}
  \subfigure[]{
    \label{fig:num_delta_y}
    \includegraphics[width = 0.30\textwidth, bb = 0 0 1293 683]
    {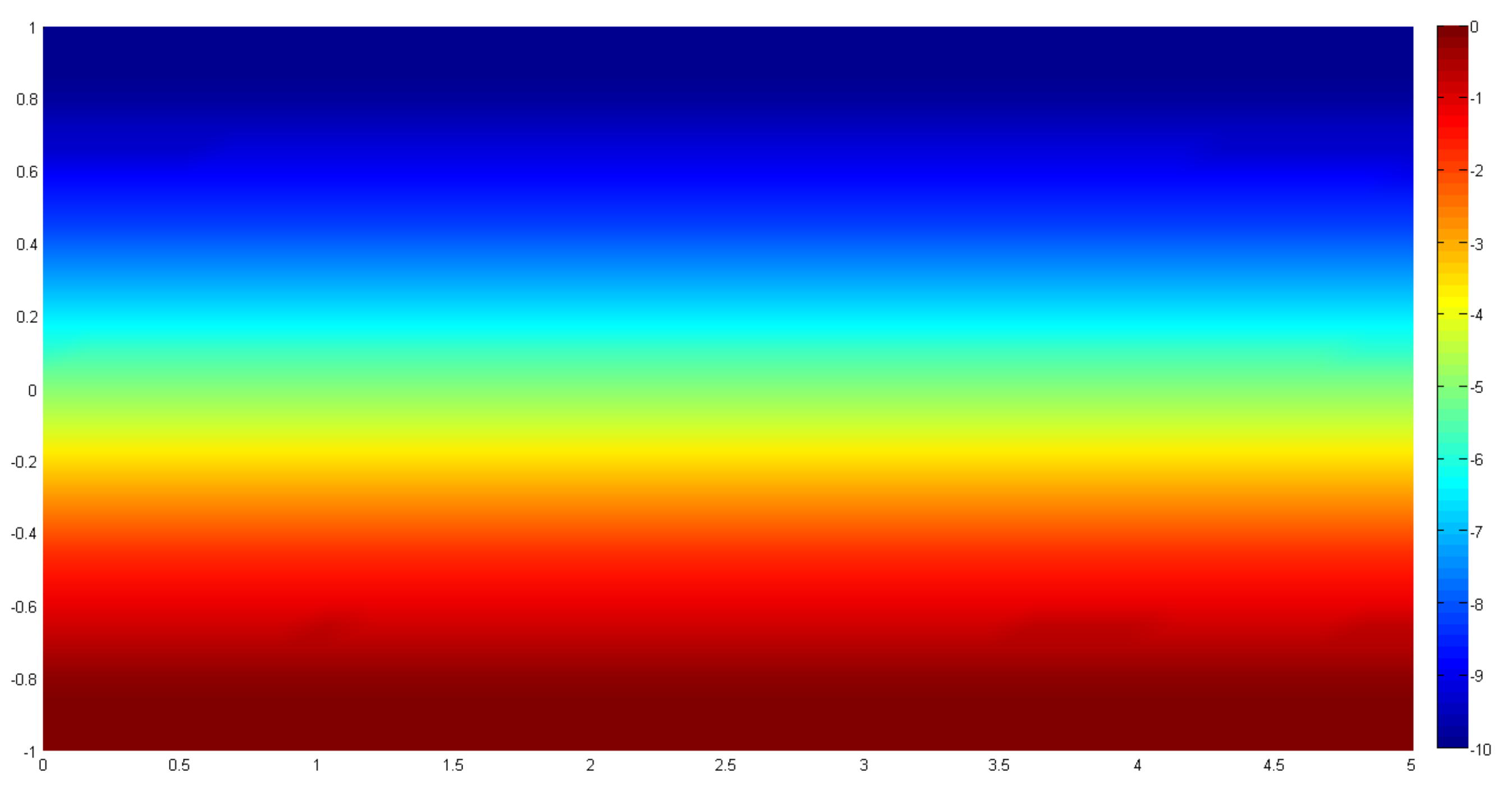}}
  \subfigure[]{
    \label{fig:num_tau_xy}
    \includegraphics[width = 0.30\textwidth, bb = 0 0 1287 675]
    {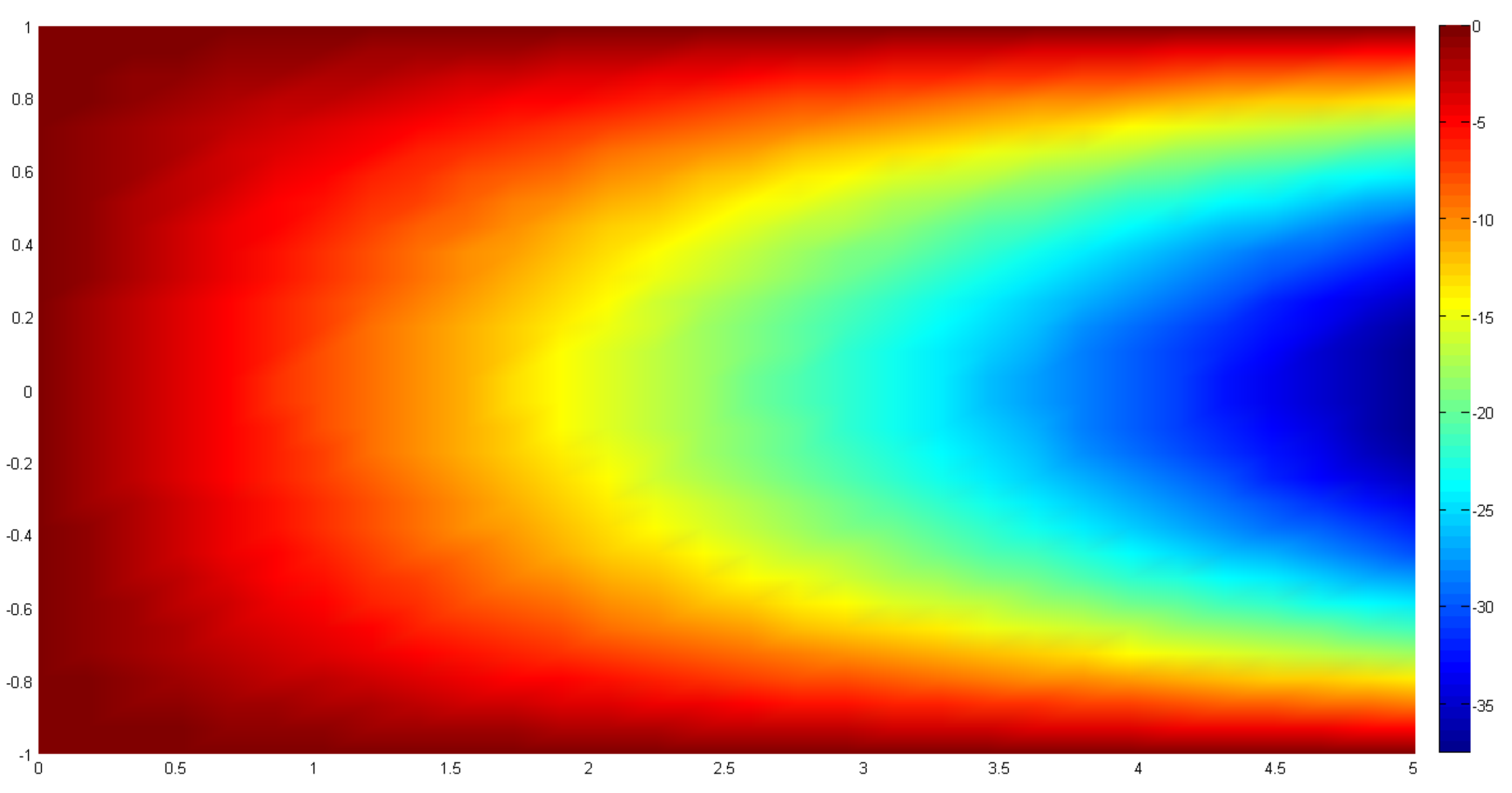}}
  \caption{Analytical solution for $\sigma_x$ (a), $\sigma_y$ (b),
    and, $\tau_{xy}$ (c) and numerical solution for $\sigma_x$ (d), $\sigma_y$ (e),
    and, $\tau_{xy}$ (f) generated by the IGA-L method with $11 \times 11$ control points
    and $16 \times 16$ Greville collocation points.}
  \label{fig:beam_solution}
\end{figure}

\begin{figure}[!htb]
\centering
 \subfigure[]{\label{fig:err_delta_x_igl}
    \includegraphics[width = 0.30\textwidth]
                {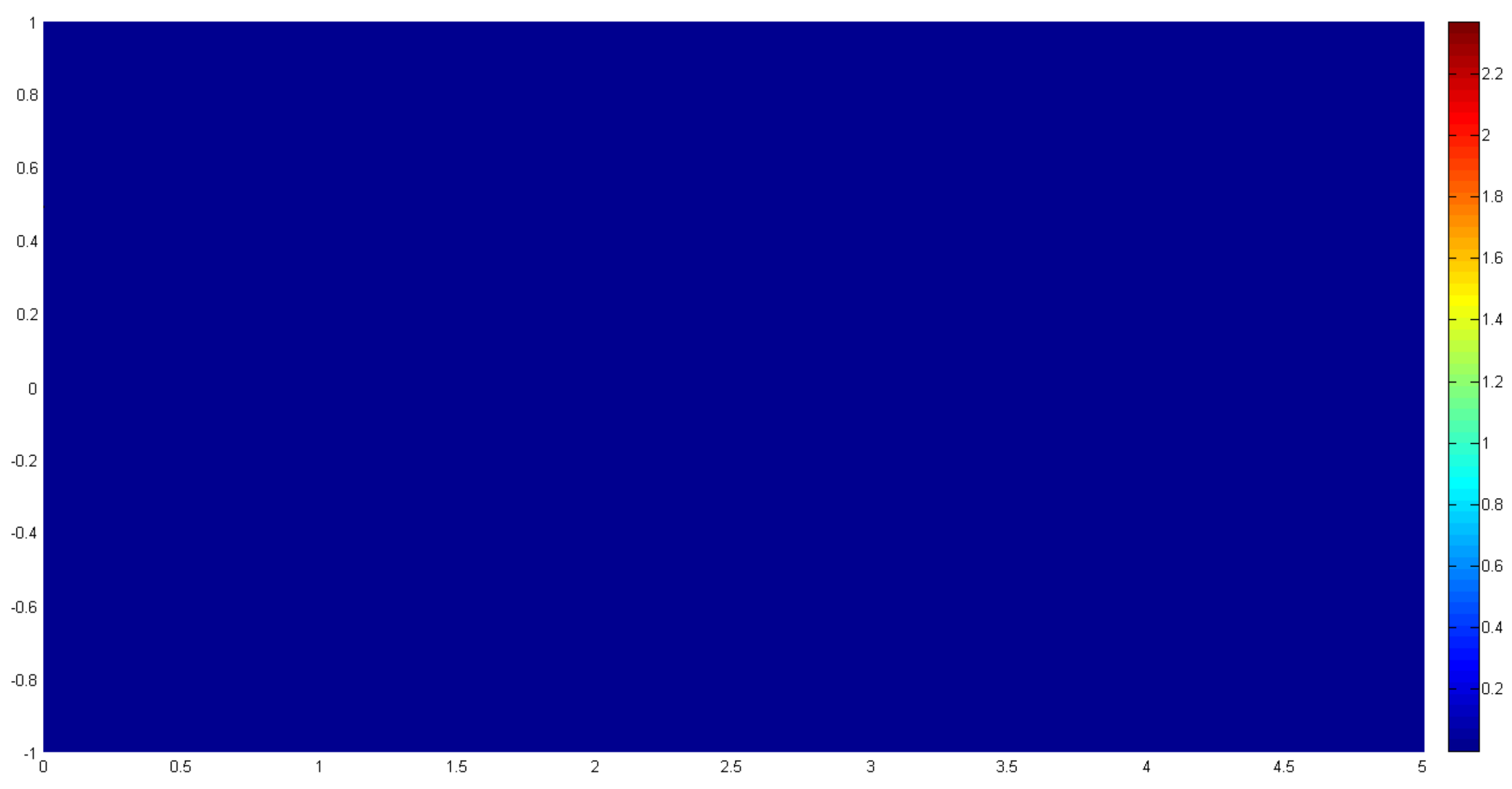}}
  \subfigure[]{\label{fig:err_delta_y_igl}
    \includegraphics[width = 0.30\textwidth]
                {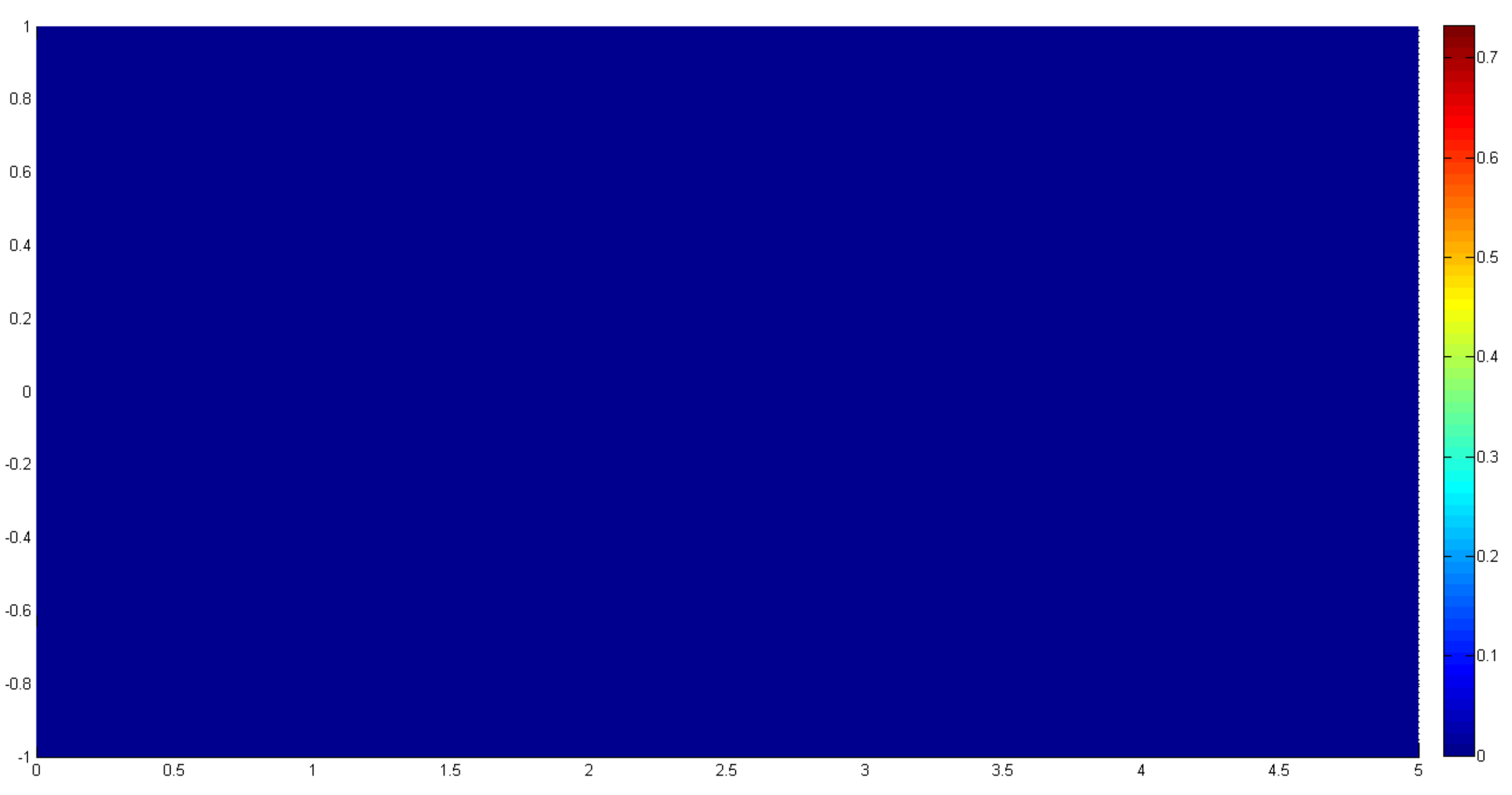}}
  \subfigure[]{\label{fig:err_tau_xy_igl}
    \includegraphics[width = 0.30\textwidth]
                {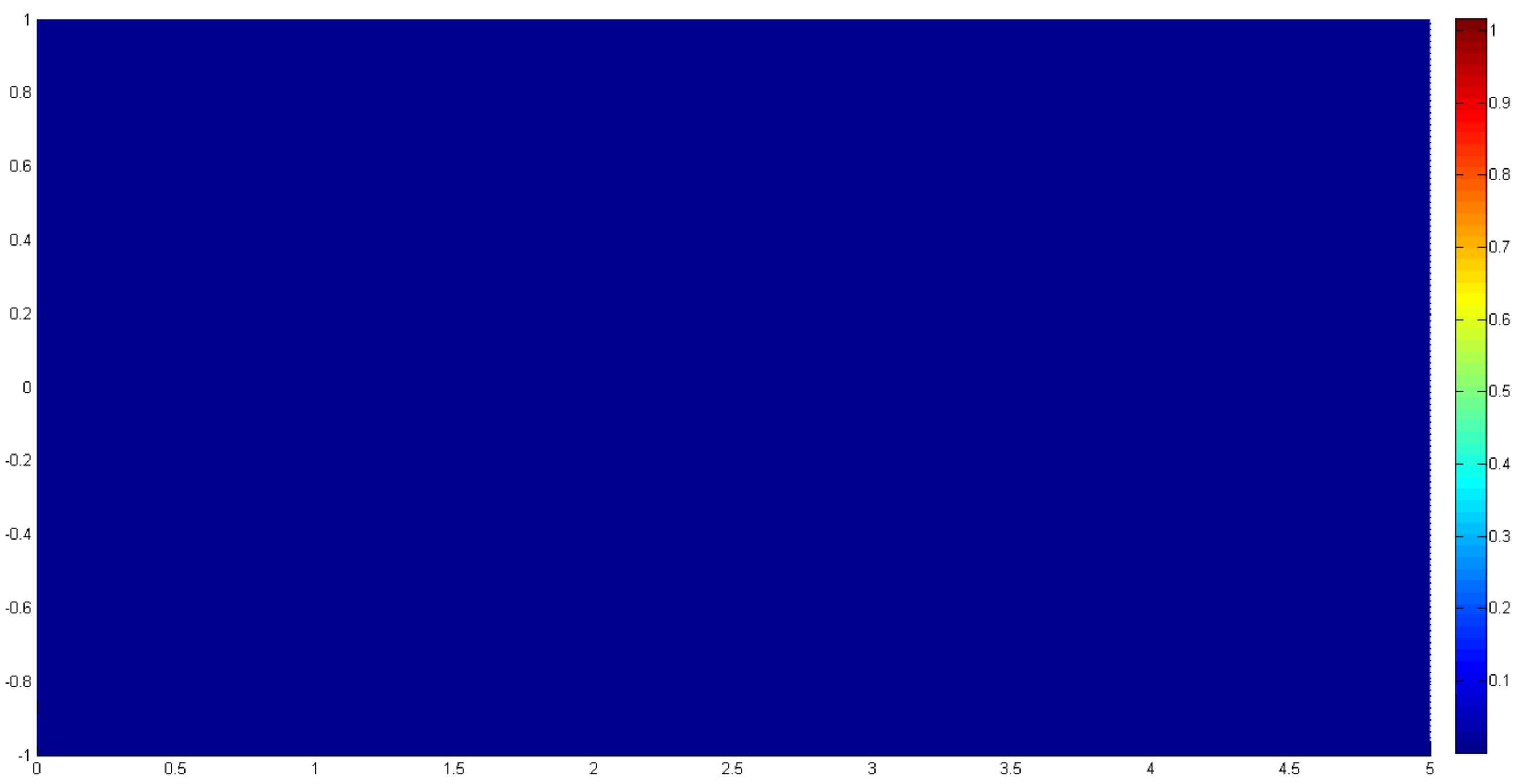}}
  \subfigure[]{\label{fig:err_delta_x_igc}
    \includegraphics[width = 0.30\textwidth]
                {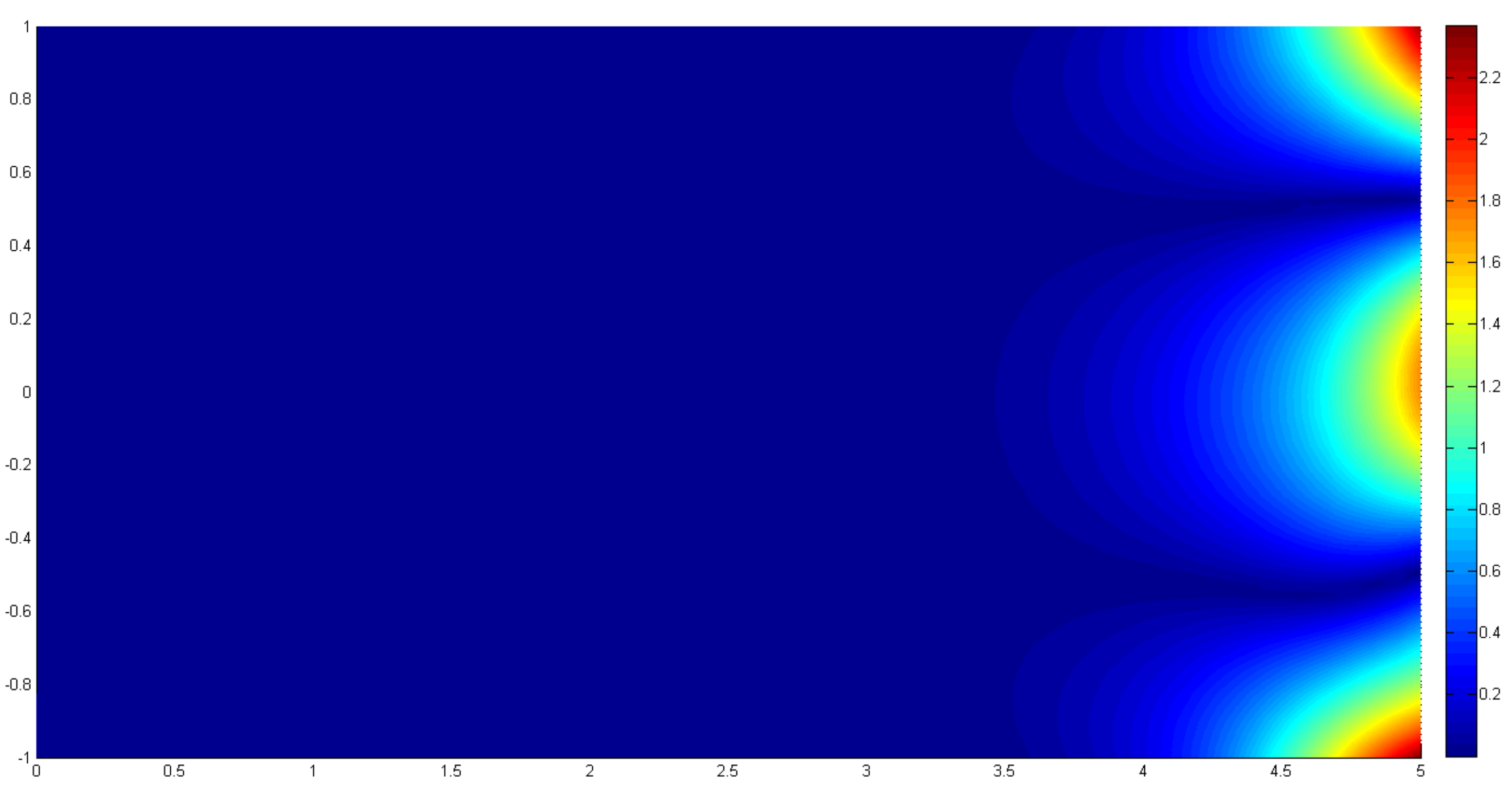}}
  \subfigure[]{\label{fig:err_delta_y_igc}
    \includegraphics[width = 0.30\textwidth]
                {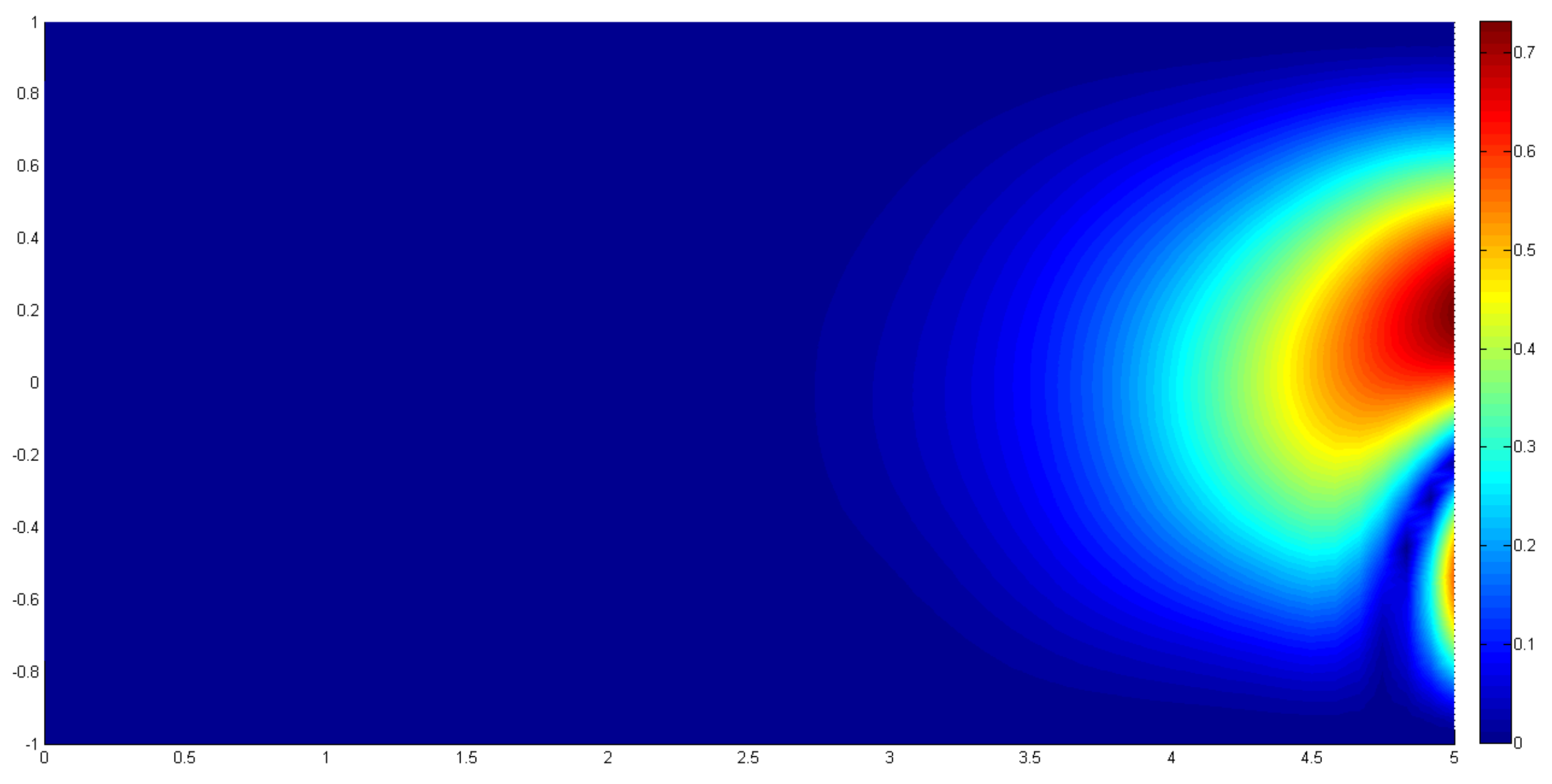}}
  \subfigure[]{\label{fig:err_tau_xy_igc}
    \includegraphics[width = 0.30\textwidth]
                {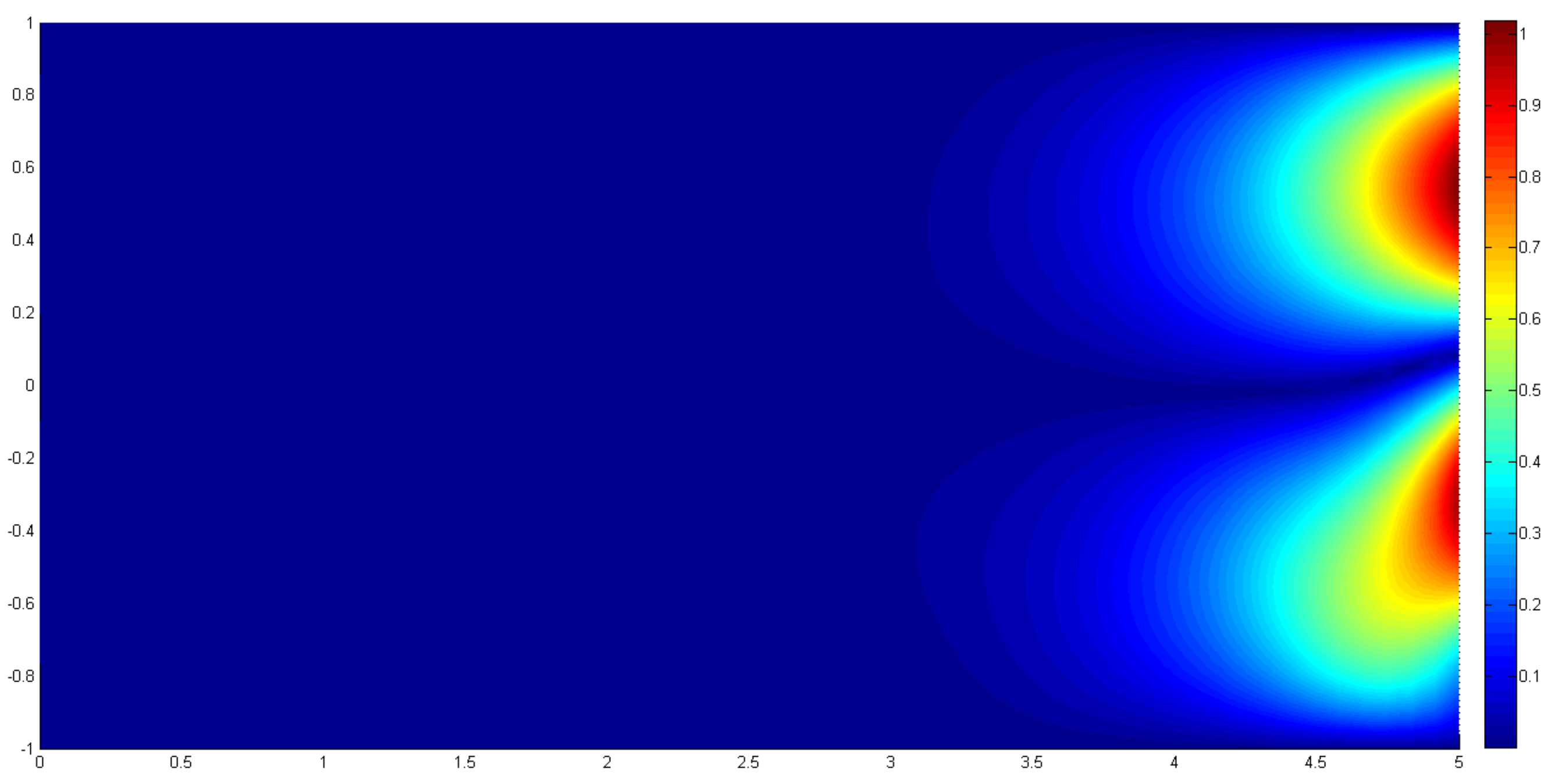}}
  \subfigure[]{\label{fig:err_delta_x_igsc}
    \includegraphics[width = 0.30\textwidth, height = 0.17\textwidth]
                {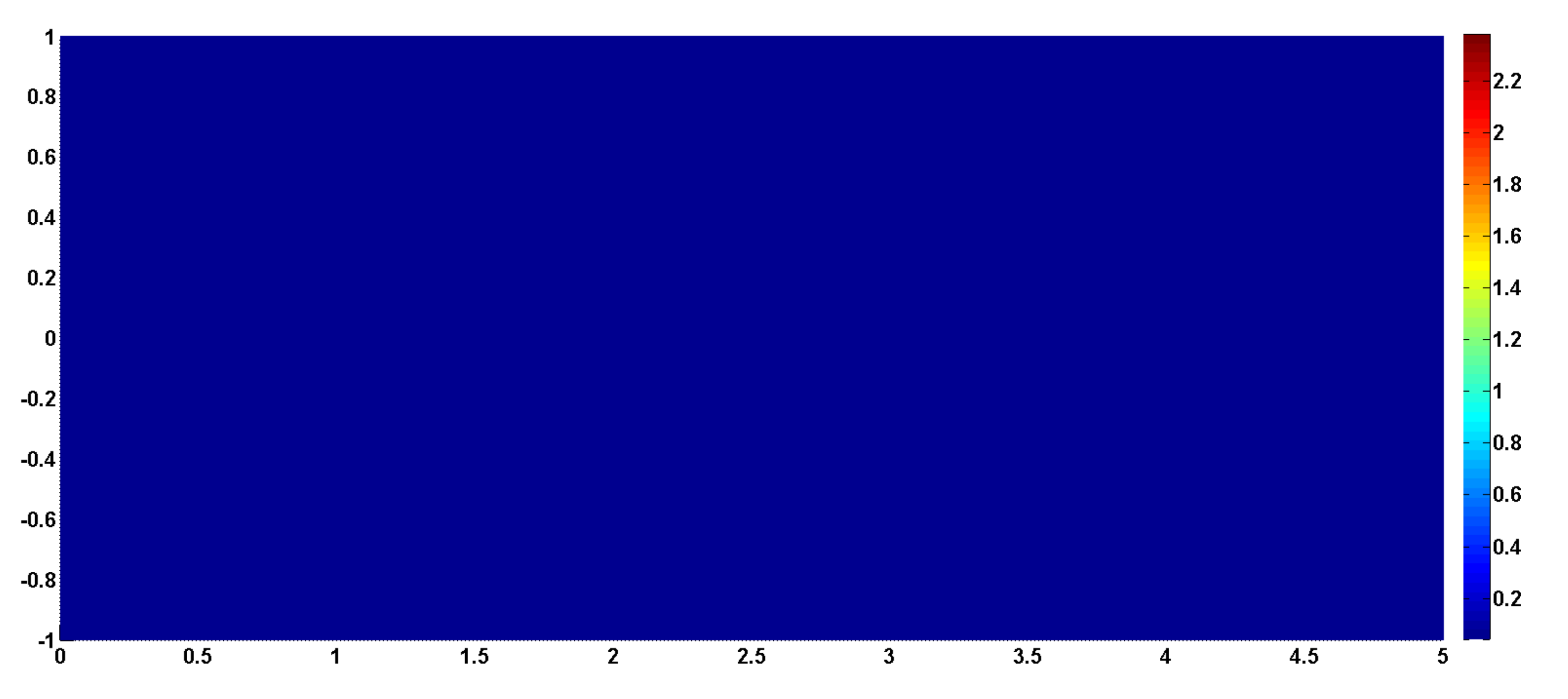}}
  \subfigure[]{\label{fig:err_delta_y_igsc}
    \includegraphics[width = 0.30\textwidth, height = 0.17\textwidth]
                {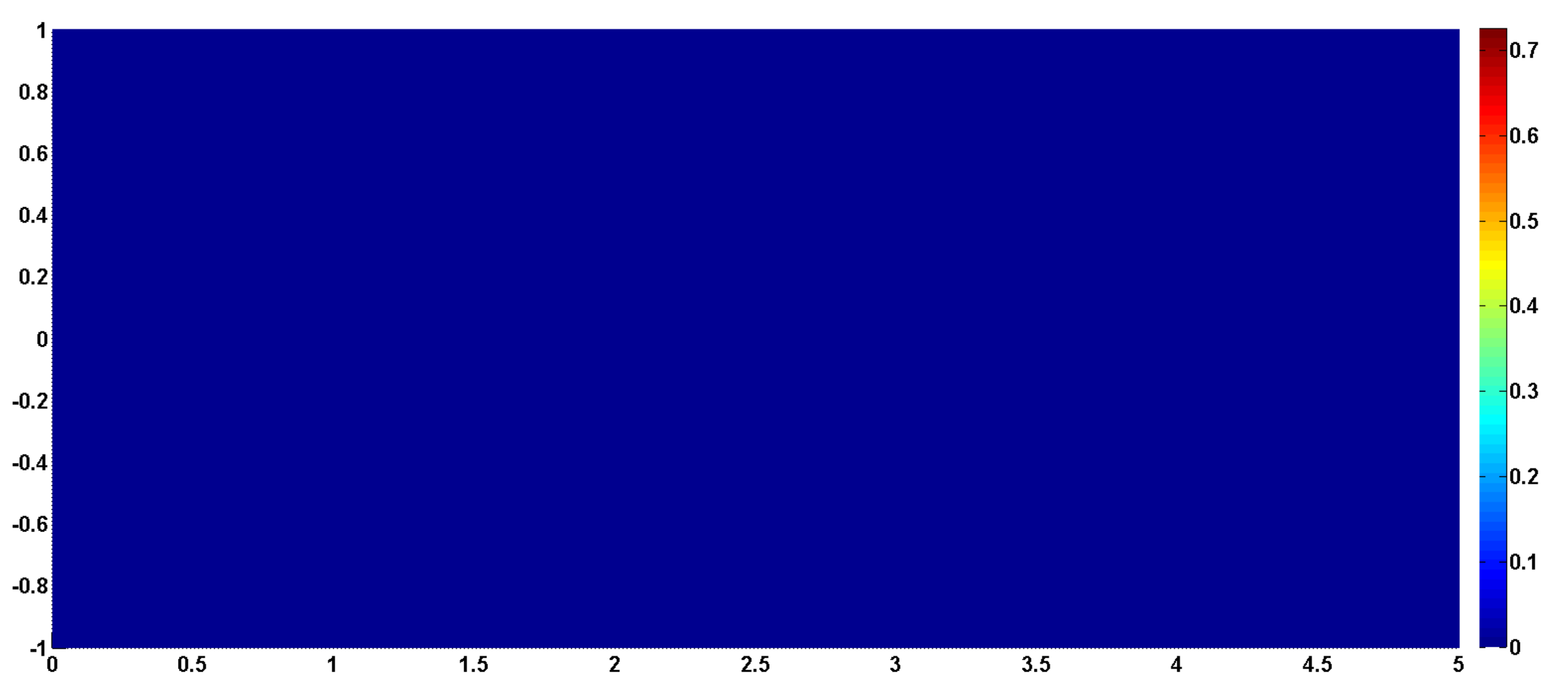}}
  \subfigure[]{\label{fig:err_tau_xy_igsc}
    \includegraphics[width = 0.30\textwidth, height = 0.17\textwidth]
                {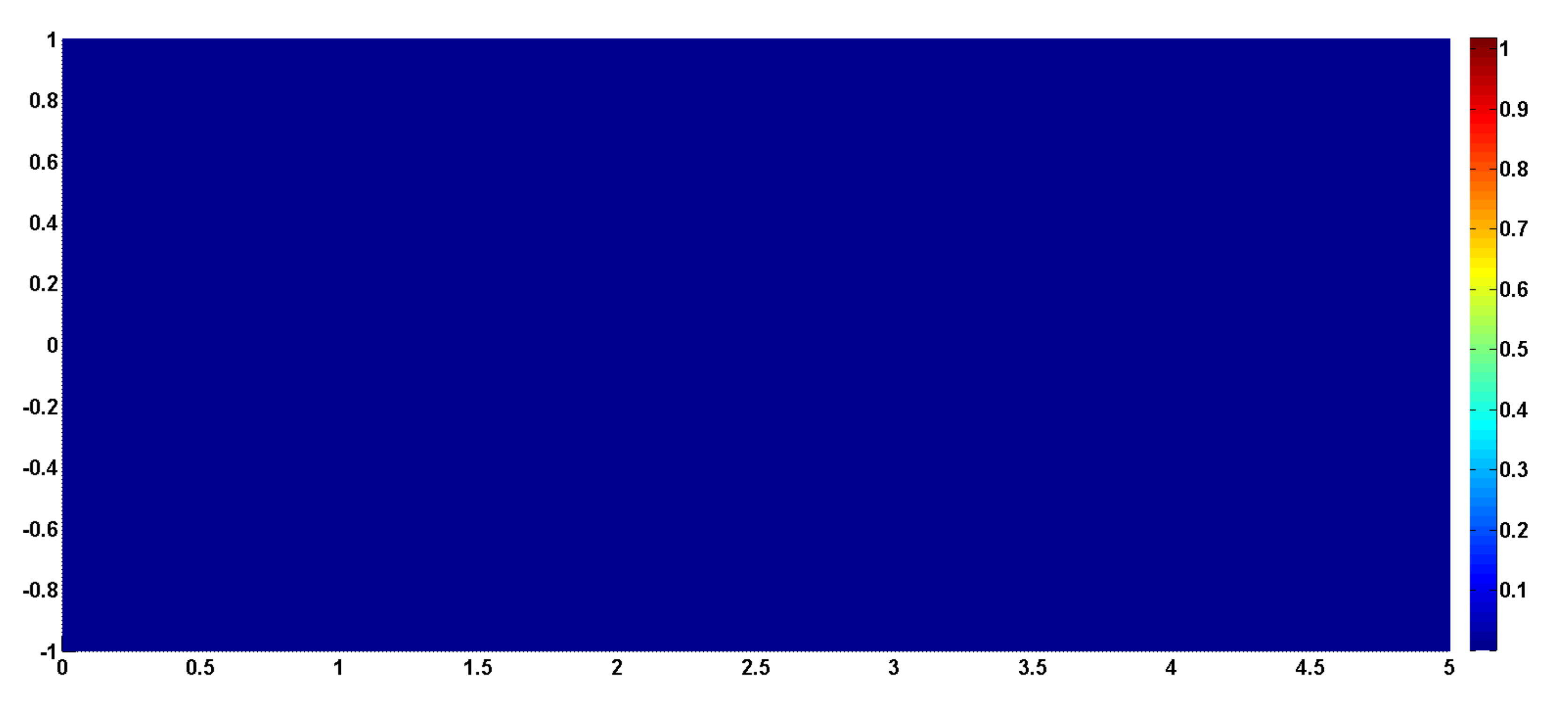}}
  \caption{Absolute error~\pref{eq:abs_err} distribution for the IGA-L
    solution $\sigma_x$ (a), $\sigma_y$ (b), and $\tau_{xy}$ (c),
    the IGA-C solution $\sigma_x$ (d), $\sigma_y$ (e), and $\tau_{xy}$ (f), and the IGA-SC solution $\sigma_x$ (g), $\sigma_y$ (h), and $\tau_{xy}$ (i).}
  \label{fig:beam_error}
\end{figure}

 \textbf{Example IV:}
 elastic problem, i.e., the simply supported
    beam (see Fig.~\ref{fig:beam}).
 As illustrated in Fig.~\ref{fig:beam},
    the simply supported beam with a rectangular cross section has
    depth $h$ and length $2l$.
 Uniformly distributed loading $q$ was applied on the upper
    surface,
    and equilibrium was maintained by reaction force $ql$ at both ends.
 Here, the body force need not be considered.
 The analytical solution of the simply supported beam problem is
 \begin{align*}
        \sigma_x & = \frac{6q}{h^3} (l^2-x^2) y + q \frac{y}{h}
            \left(4 \frac{y^2}{h^2} - \frac{3}{5} \right), \\
        \sigma_y & = - \frac{q}{2} \left(1+\frac{y}{h}\right)
            \left( 1 - \frac{2y}{h} \right)^2 , \\
        \tau_{xy} & = - \frac{6q}{h^3} x \left(\frac{h^2}{4} - y^2
            \right).
 \end{align*}

 We calculated the simply supported beam problem using the IGA-L, IGA-C,
    and IGA-SC methods,
    with $q = 10,\ h = 2,\ l = 5$ (see Fig.~\ref{fig:beam}).
 The physical domain was represented by a cubic B-spline patch
    presented in Appendix A4.
 Fig.~\ref{fig:beam_solution} illustrates the analytical
    and numerical solutions for $\sigma_x,\ \sigma_y$, and
    $\tau_{xy}$,
    generated by the IGA-L method with $11 \times 11$ control points
    and $16 \times 16$ Greville collocation points.
 The relative errors for $\sigma_x,\ \sigma_y$, and $\tau_{xy}$  of the IGA-L
    solutions are $1.10 \times 10^{-5}$, $3.29 \times 10^{-4}$,
    and $7.20 \times 10^{-5}$, respectively.
 For comparison, the relative errors for $\sigma_x,\ \sigma_y$,
    and $\tau_{xy}$ of the IGA-SC solutions with $11 \times 11$ control points and $16 \times 16$ collocation points are, $3.20 \times 10^{-5}$, $4.76 \times 10^{-4}$, and $8.70 \times 10^{-5}$, respectively;
    those of the IGA-C solutions with $11 \times 11$ control points and Greville collocation points are $4.0 \times 10^{-3}, 6.14 \times 10^{-2}$, and $2.15 \times 10^{-2}$.
 Additionally, Fig.~\ref{fig:beam_error} demonstrates the absolute error
    distribution for the IGA-L, IGA-C, and IGA-SC solutions.

\begin{figure}[!htb]
\centering
  \subfigure[]{
    \label{subfig:iga_l_delta_x}
    \includegraphics[width = 0.32\textwidth]
    {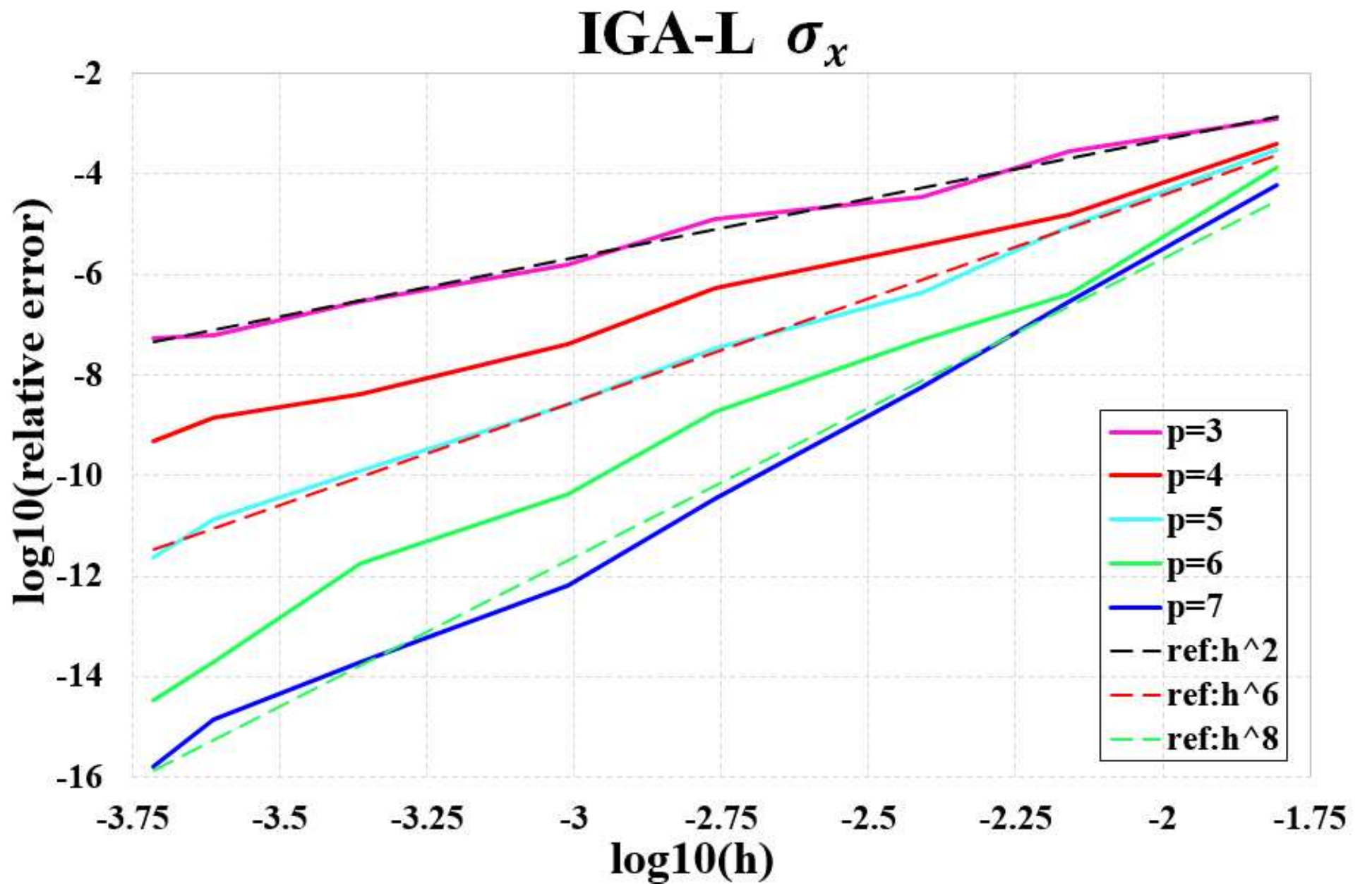}}
  \subfigure[]{
    \label{subfig:iga_l_delta_y}
    \includegraphics[width = 0.32\textwidth]
    {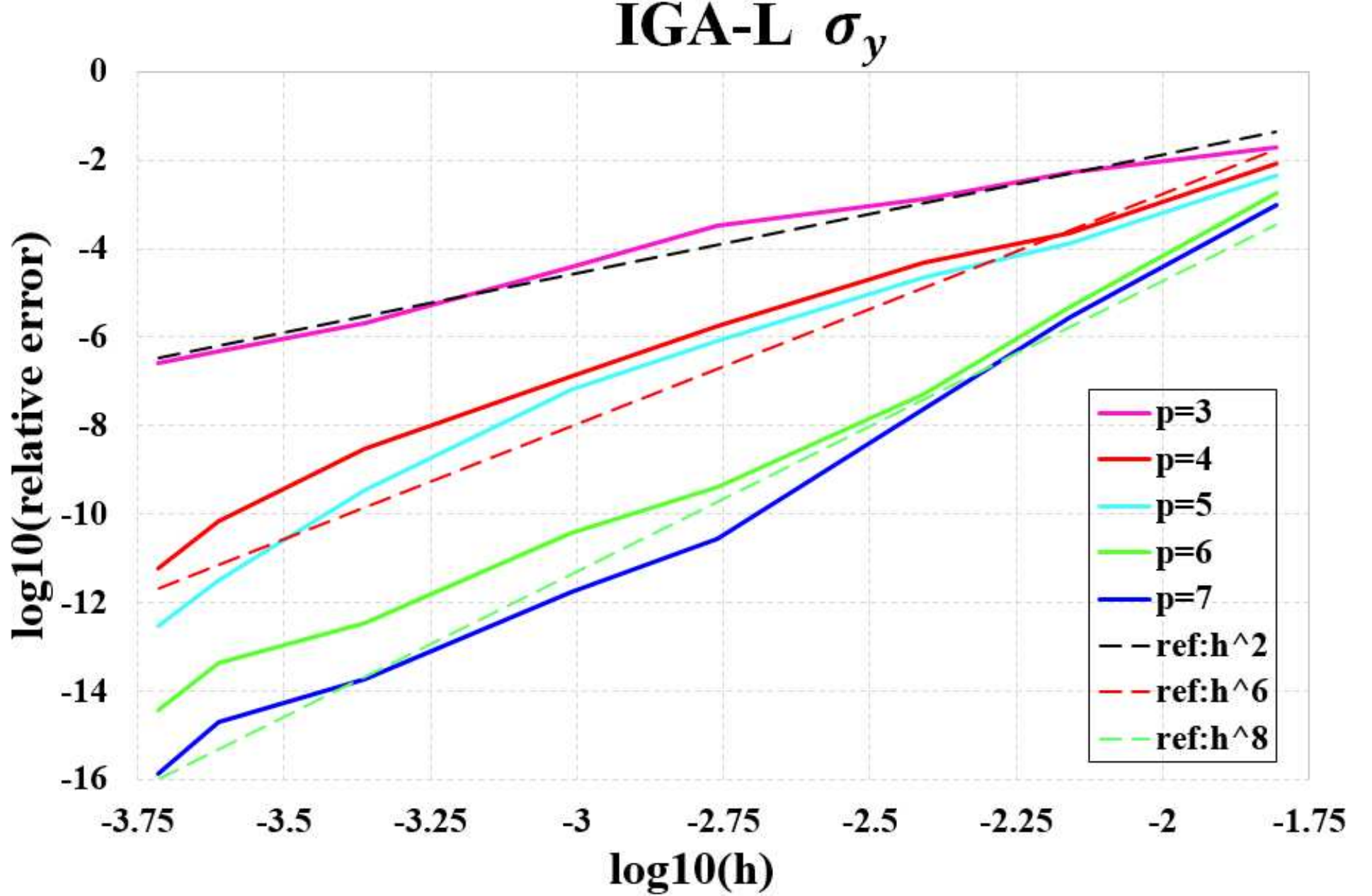}}
  \subfigure[]{
    \label{subfig:iga_l_tau_xy}
    \includegraphics[width = 0.32\textwidth]
    {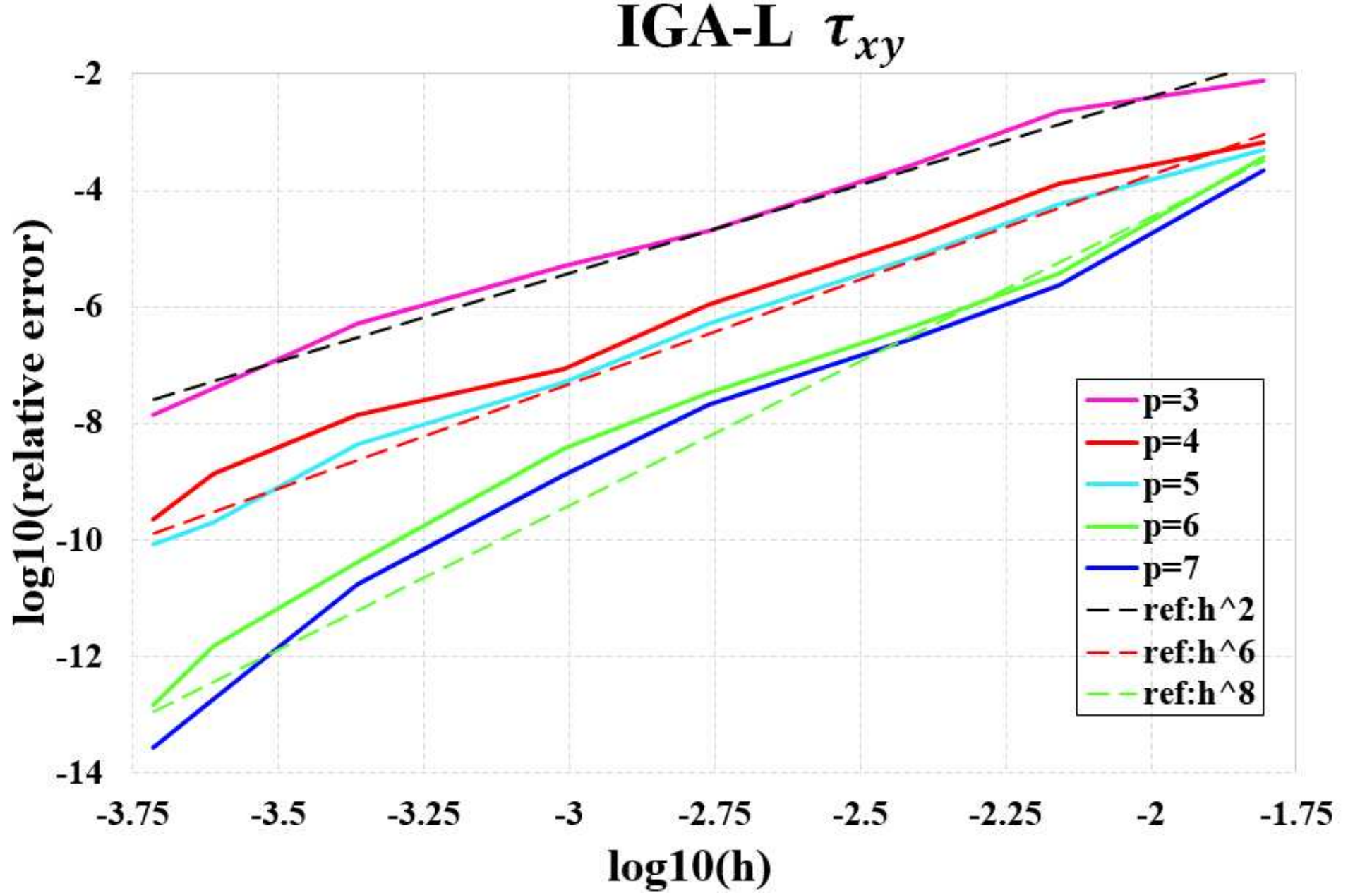}}
  \subfigure[]{
    \label{subfig:iga_c_delta_x}
    \includegraphics[width = 0.32\textwidth]
    {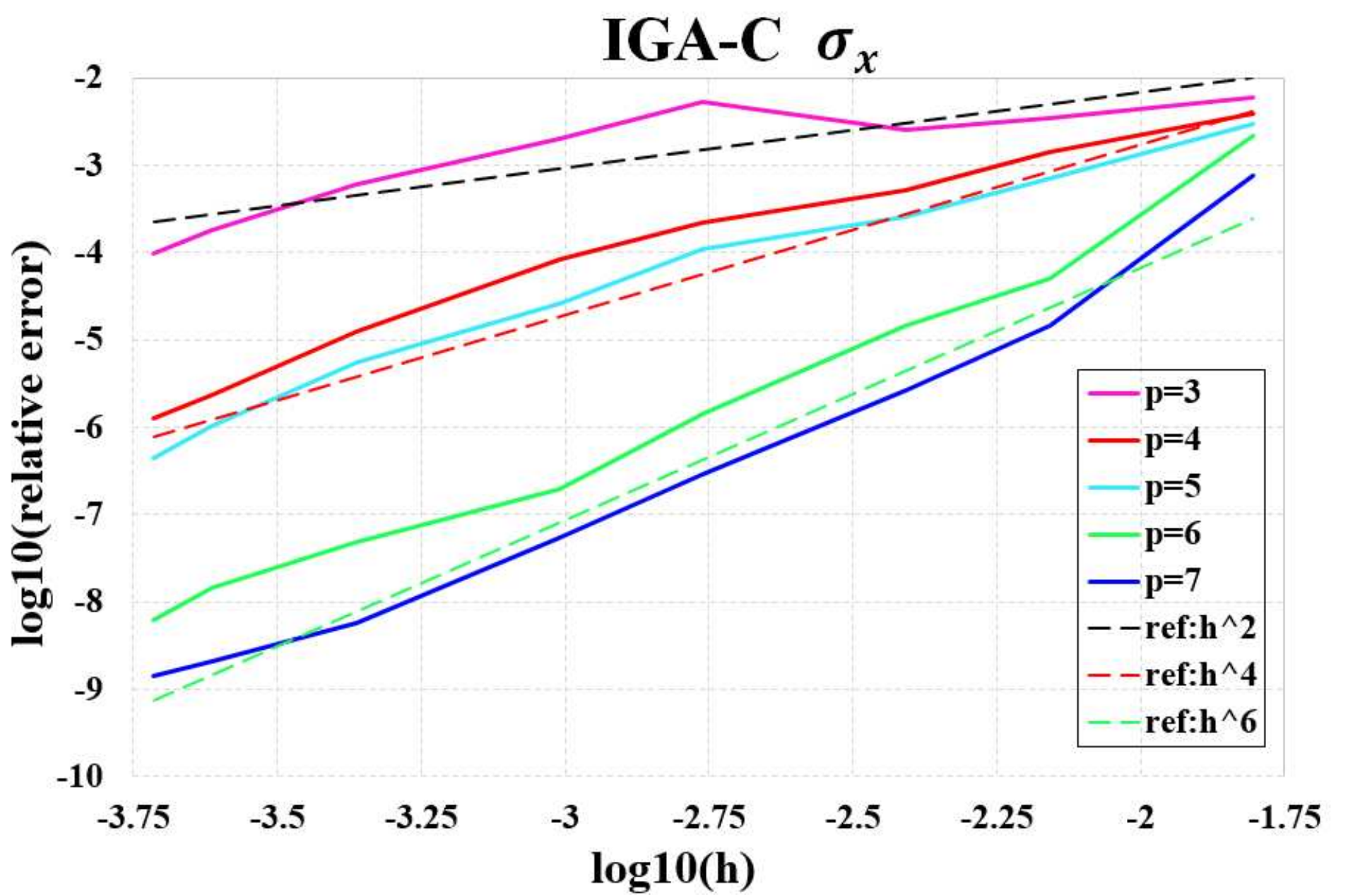}}
  \subfigure[]{
    \label{subfig:iga_c_delta_y}
    \includegraphics[width = 0.32\textwidth]
    {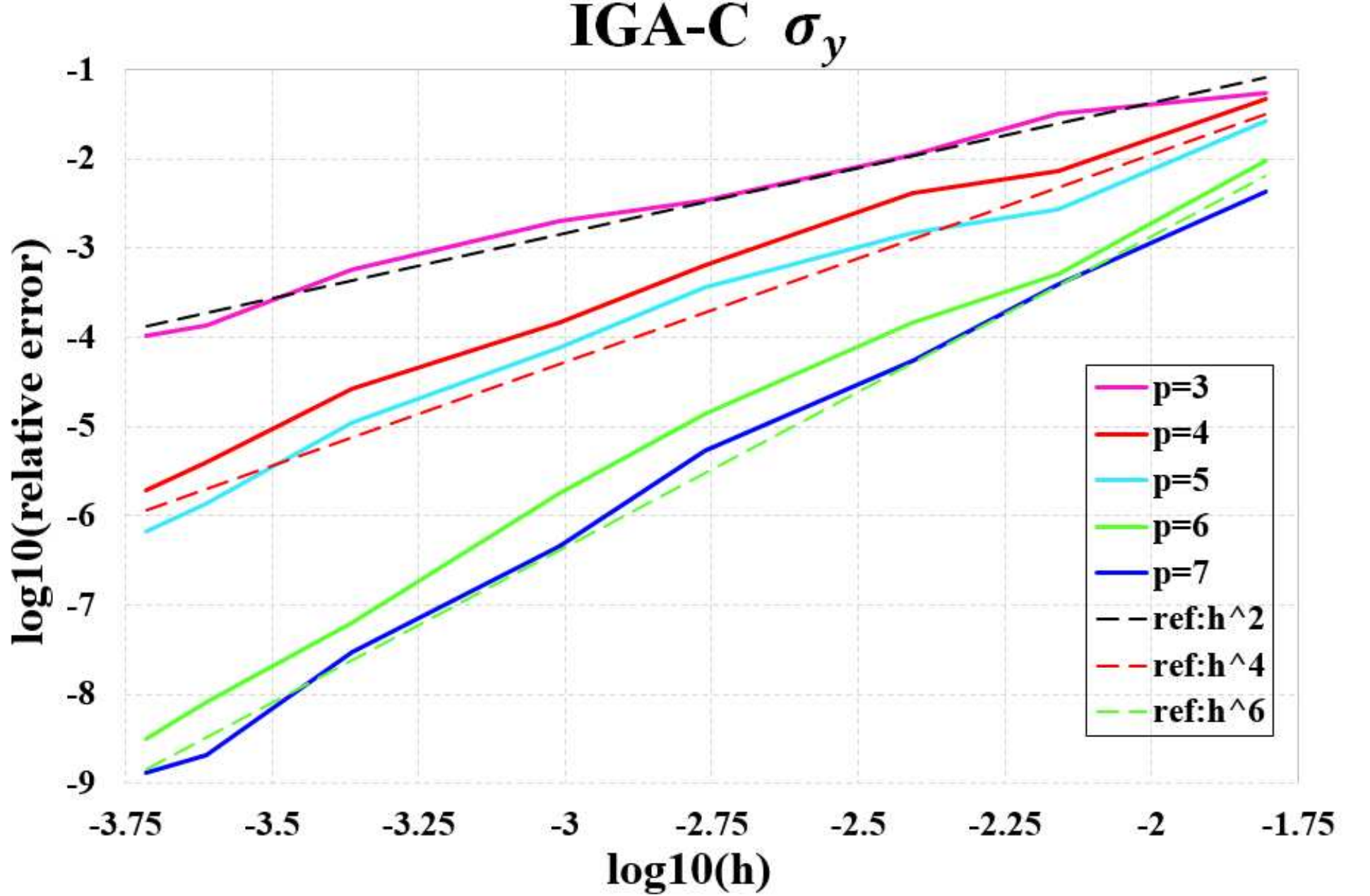}}
  \subfigure[]{
    \label{subfig:iga_c_tau_xy}
    \includegraphics[width = 0.32\textwidth]
    {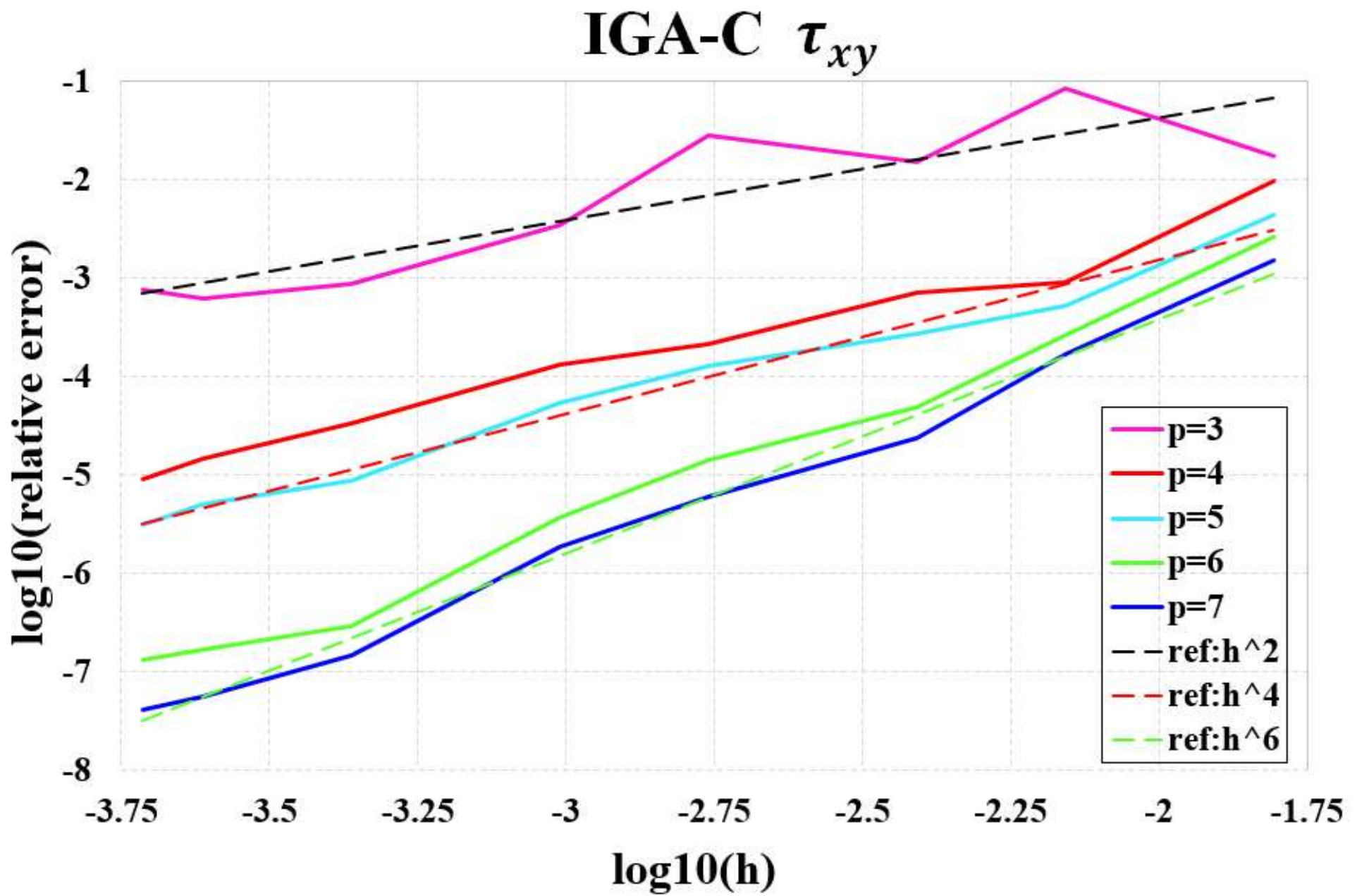}}
  \subfigure[]{
    \label{subfig:iga_sc_delta_x}
    \includegraphics[width = 0.32\textwidth]
    {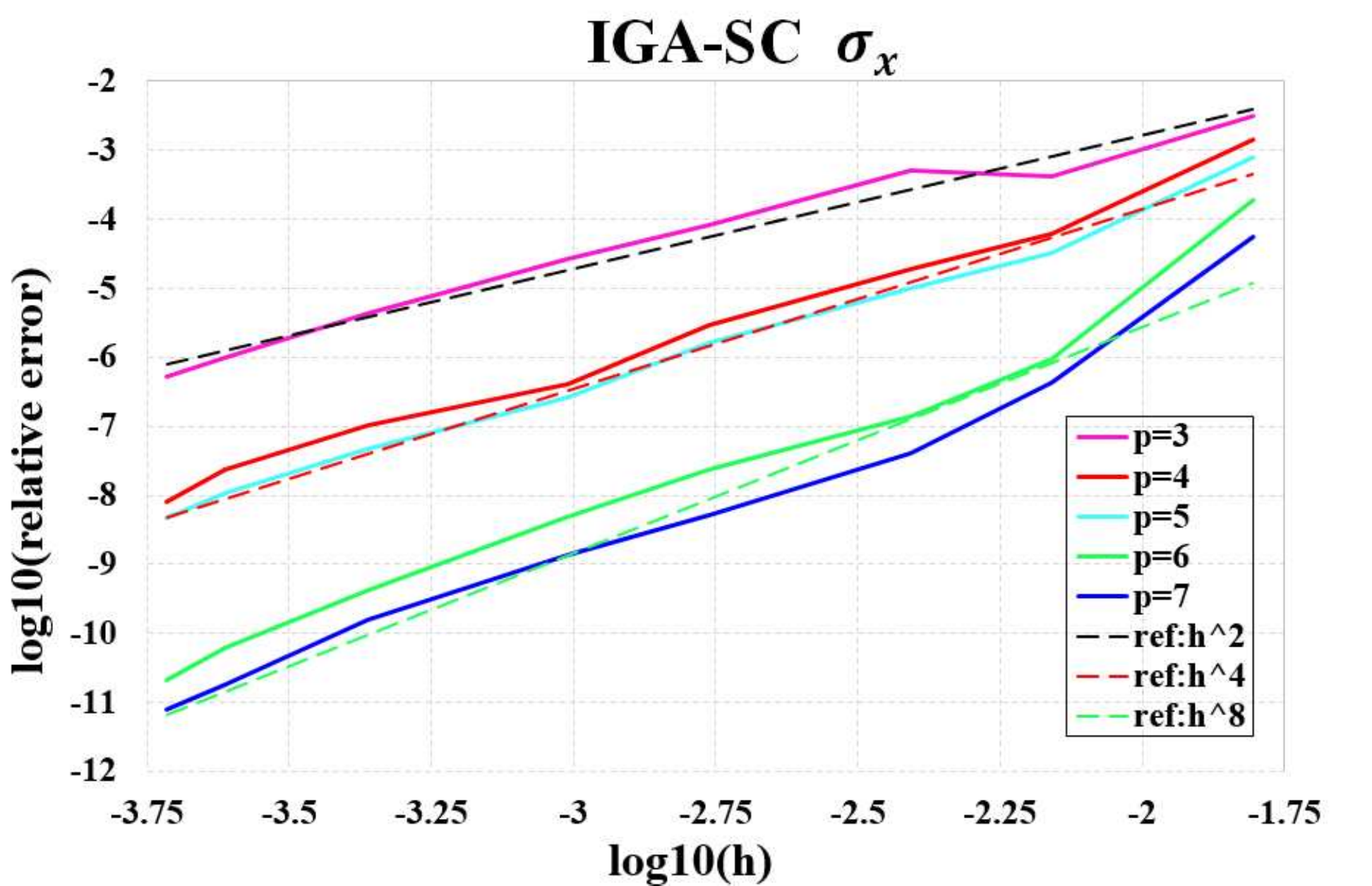}}
  \subfigure[]{
    \label{subfig:iga_sc_delta_y}
    \includegraphics[width = 0.32\textwidth]
    {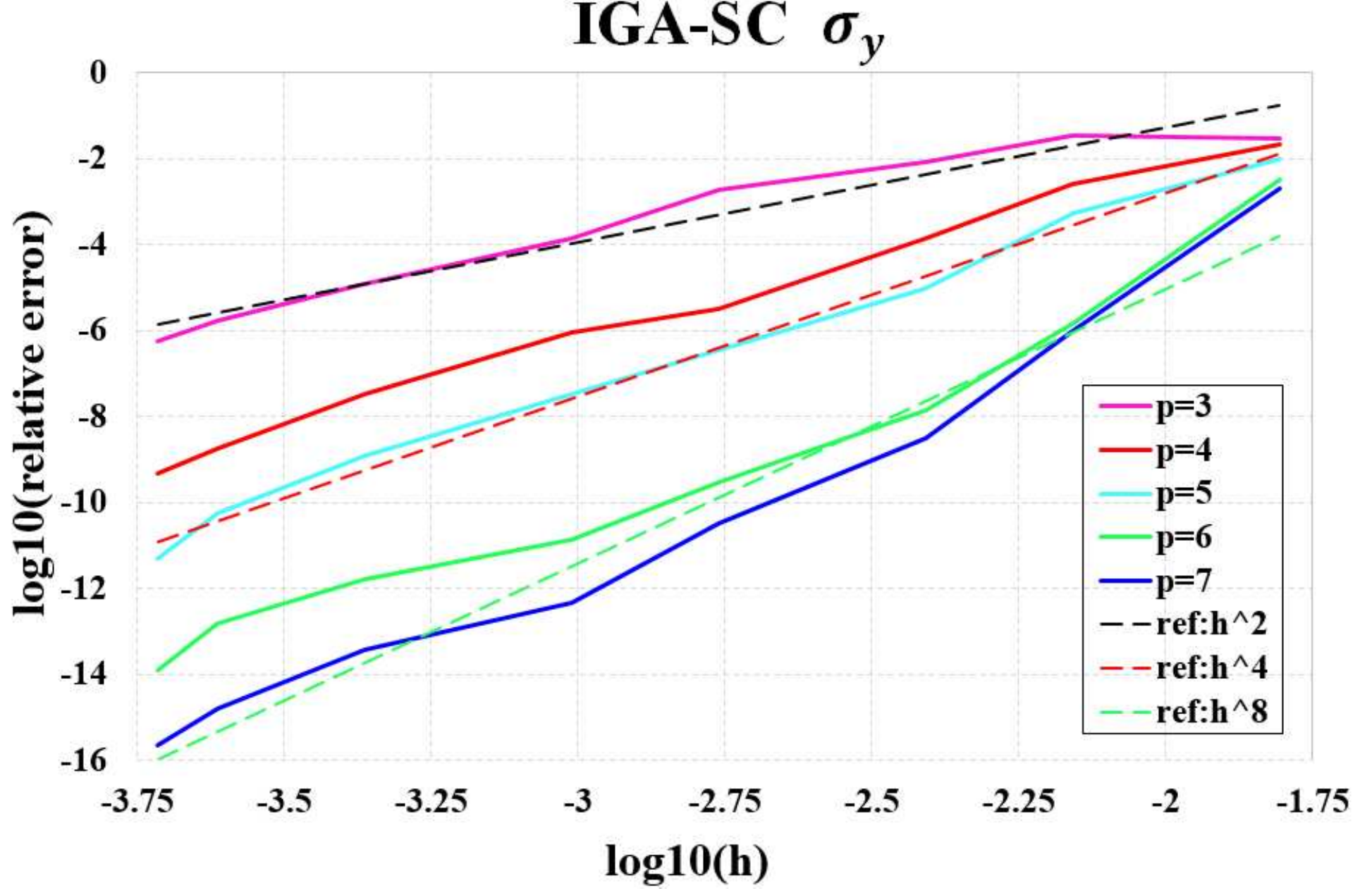}}
  \subfigure[]{
    \label{subfig:iga_sc_tau_xy}
    \includegraphics[width = 0.32\textwidth]
    {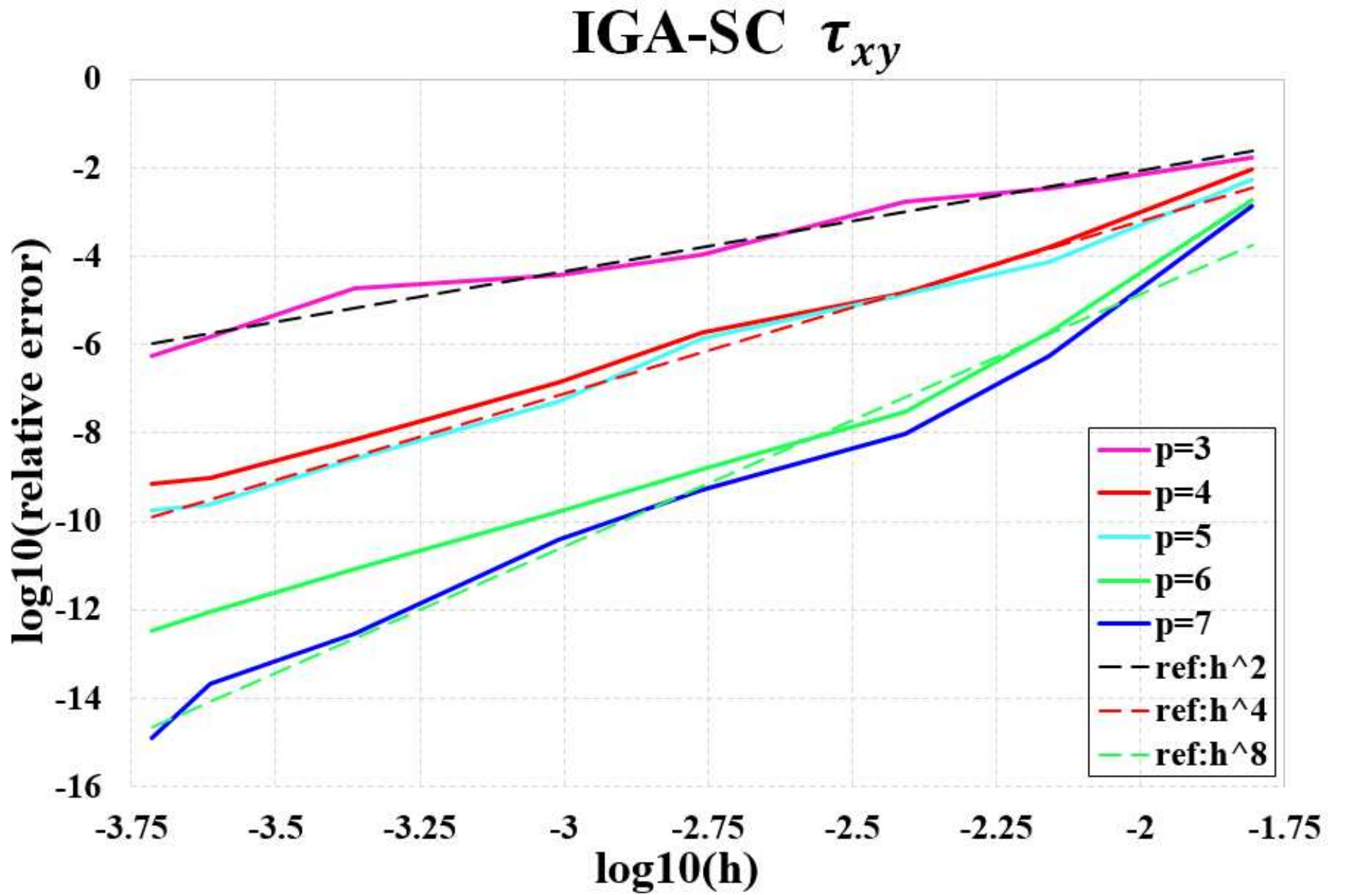}}
  \caption{Numerical results for the simply supported beam
    (Fig.~\ref{fig:beam}).
    (a, b, c) Diagrams of $\log_{10}$(Relative error) vs. $\log_{10}(h)$ for
        $\sigma_x$, $\sigma_y$, and $\tau_{xy}$, respectively, using IGA-L method.
    (d, e, f) Diagrams of $\log_{10}$(Relative error) vs. $\log_{10}(h)$ for
        $\sigma_x$, $\sigma_y$, and $\tau_{xy}$, respectively, using IGA-C method.
    (h, i, j) Diagrams of $\log_{10}$(Relative error) vs. $\log_{10}(h)$ for
        $\sigma_x$, $\sigma_y$, and $\tau_{xy}$, respectively, using IGA-SC method.}
  \label{fig:beam_comparison}
\end{figure}

\begin{figure}[!htb]
\centering
  \subfigure[]{
    \label{subfig:delta_x_time_error}
    \includegraphics[width = 0.32\textwidth]
    {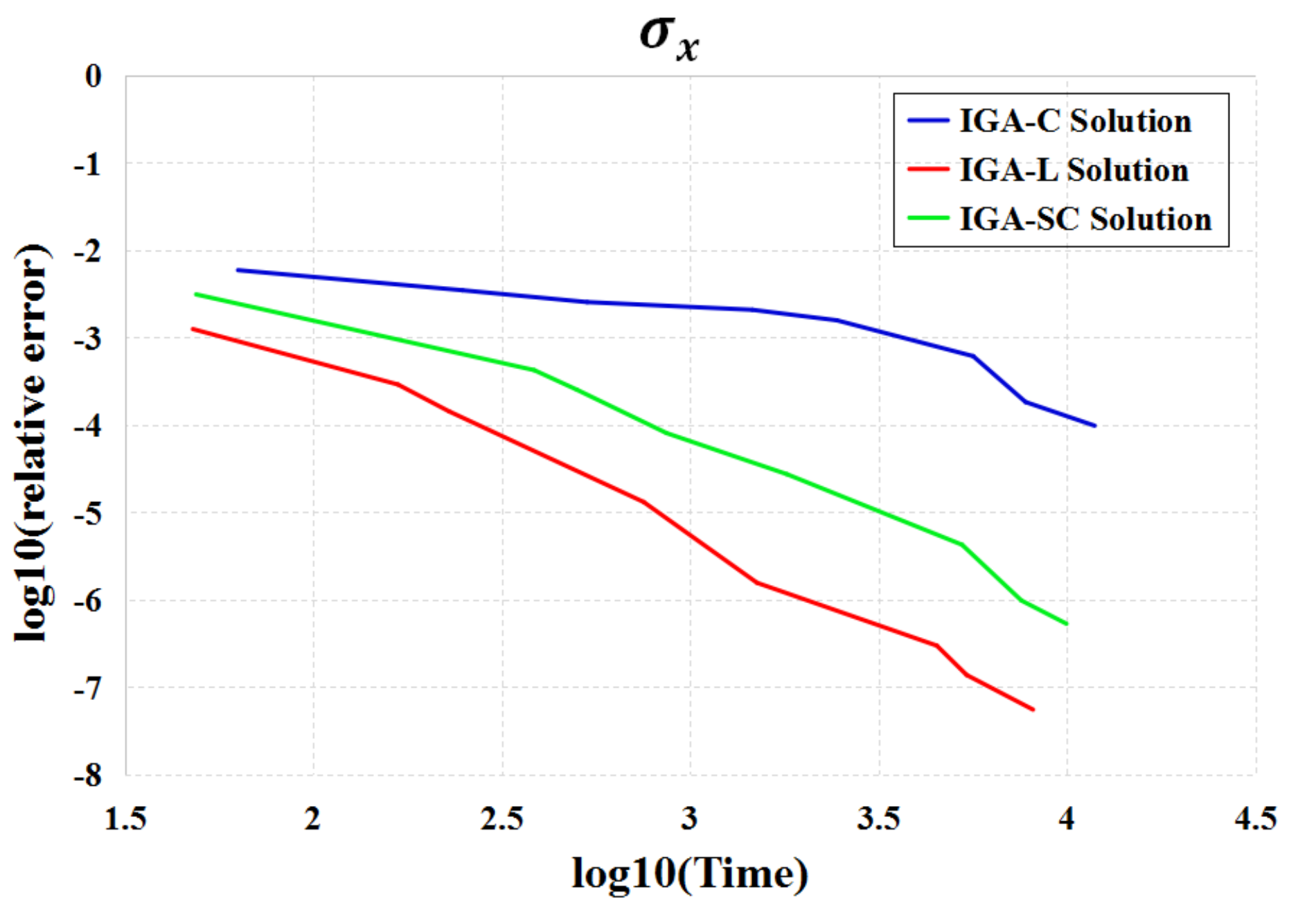}}
  \subfigure[]{
    \label{subfig:delta_y_time_error}
    \includegraphics[width = 0.32\textwidth]
    {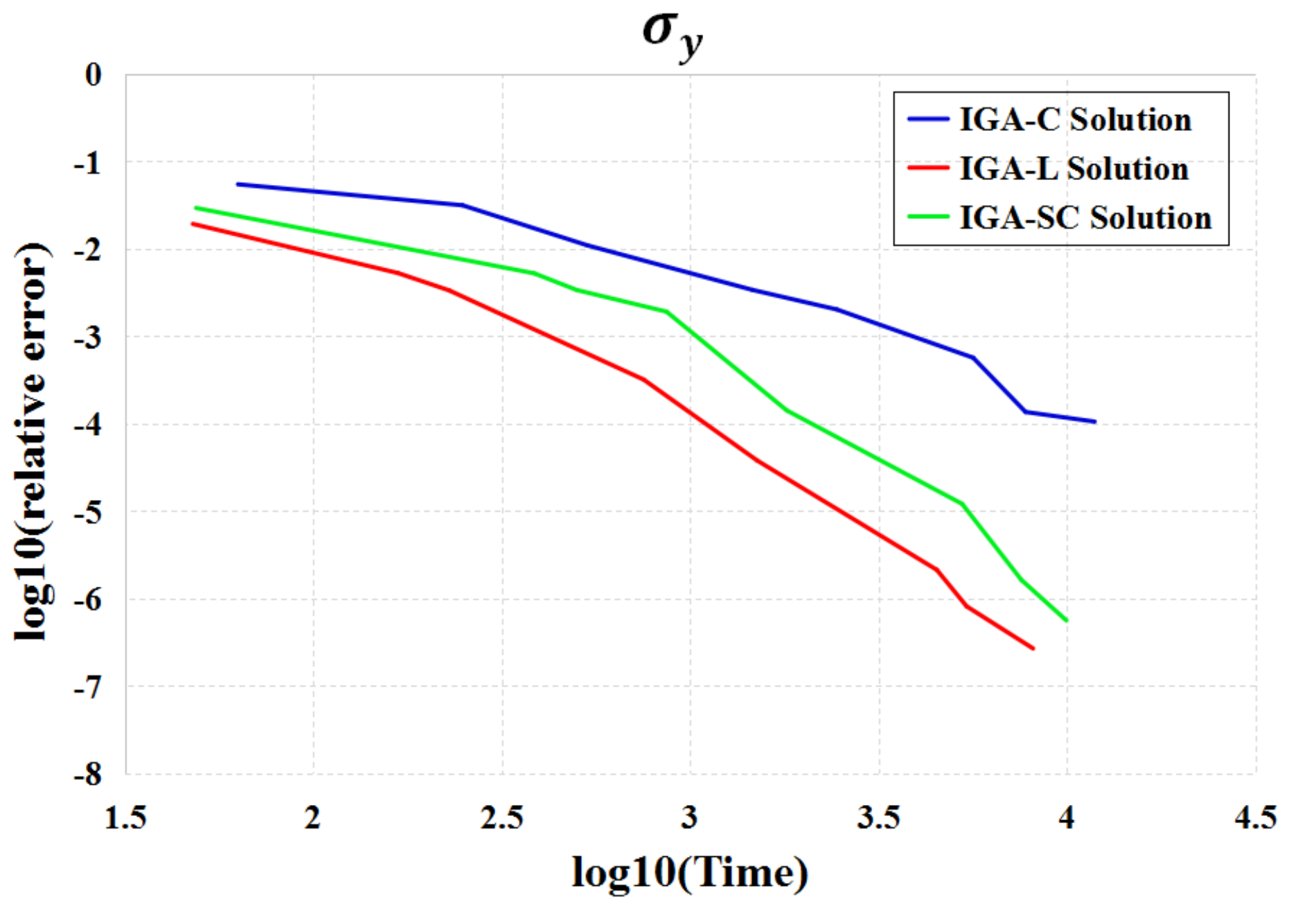}}
  \subfigure[]{
    \label{subfig:tau_xy_time_error}
    \includegraphics[width = 0.32\textwidth]
    {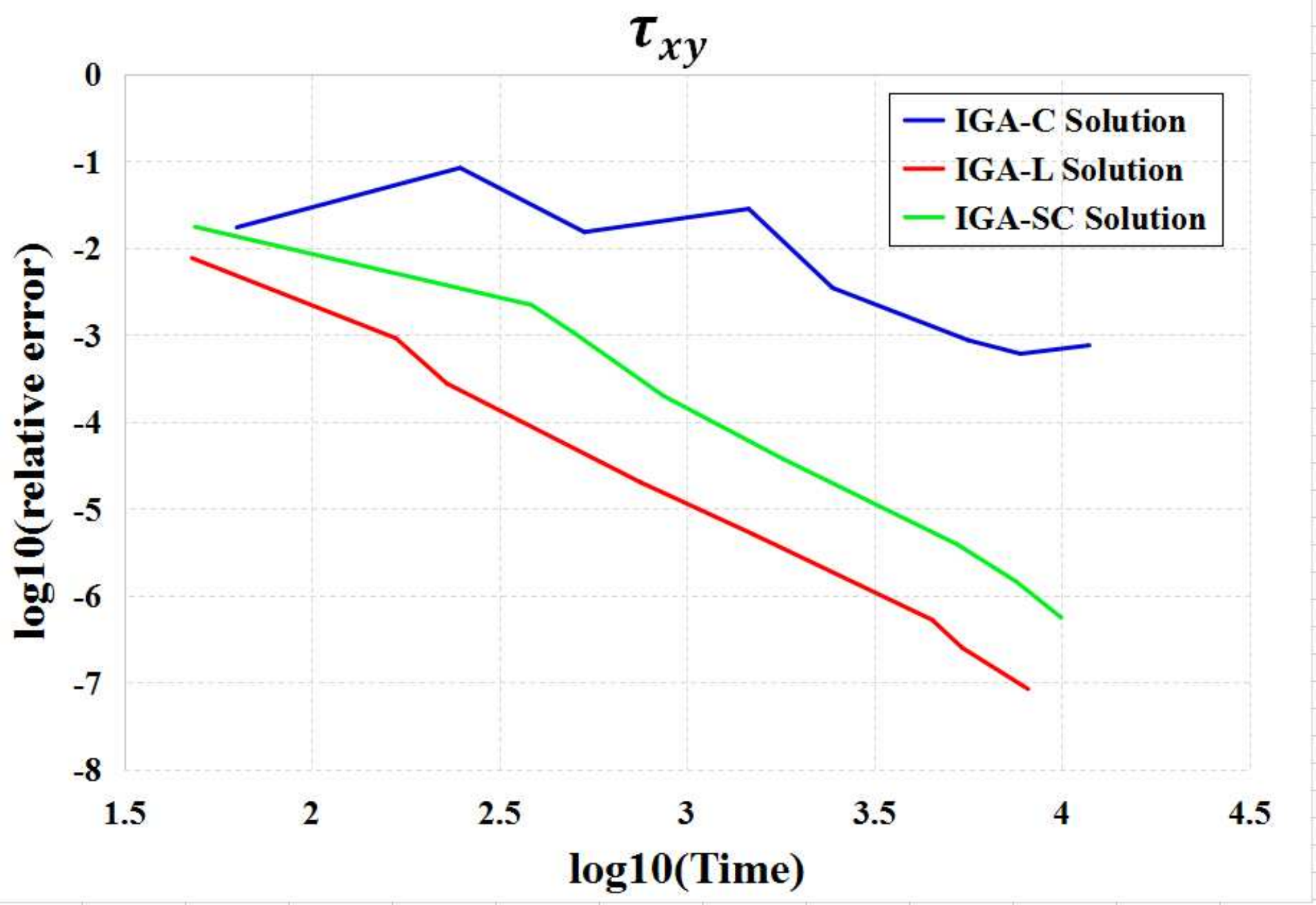}}
  \caption{Diagrams of $\log_{10}$(Time) v.s. $\log_{10}$(relative error) for the simply supported beam (Fig.~\ref{fig:beam}),
  for computing $\delta_x$ (a), $\delta_y$ (b), and $\tau_{xy}$ (c), respectively.}
  \label{fig:beam_time_error}
\end{figure}

 Furthermore, in Fig.~\ref{fig:beam_comparison},
    we demonstrate the diagrams of $\log_{10}(h)$ v.s. $\log_{10}$(relative error).
 We can see that, on one hand, while the convergence rate of IGA-C and IGA-SC
    solutions is $O(h^4)$ for degree $p = 4, 5$,
    that of IGA-L solutions is $O(h^6)$ for degree $p=4, 5$.
 On the other hand, while the convergence rate of IGA-C solutions is $O(h^6)$
    for degree $p = 6, 7$,
    that of IGA-L solutions is $O(h^8)$ for degree $p = 6, 7$.
 In addition, diagrams of $\log_{10}$(time) v.s. $\log_{10}$(relative error)
    for the three methods IGA-L, IGA-C, and IGA-SC are presented in Fig.~\ref{fig:beam_time_error},
    and diagrams for IGA-L have the best performance.

\begin{figure}[!htb]
\centering
  \subfigure[Analytical solution, IGA-L solution
            with $14$ Greville collocation points, and IGA-SC solution with $14$ collocation points, of Eq.~\pref{eq:exmp_one_dim_nuemann}.]{
    \label{fig:one-dim-igl-solution-g}
    \includegraphics[width=0.40\textwidth]{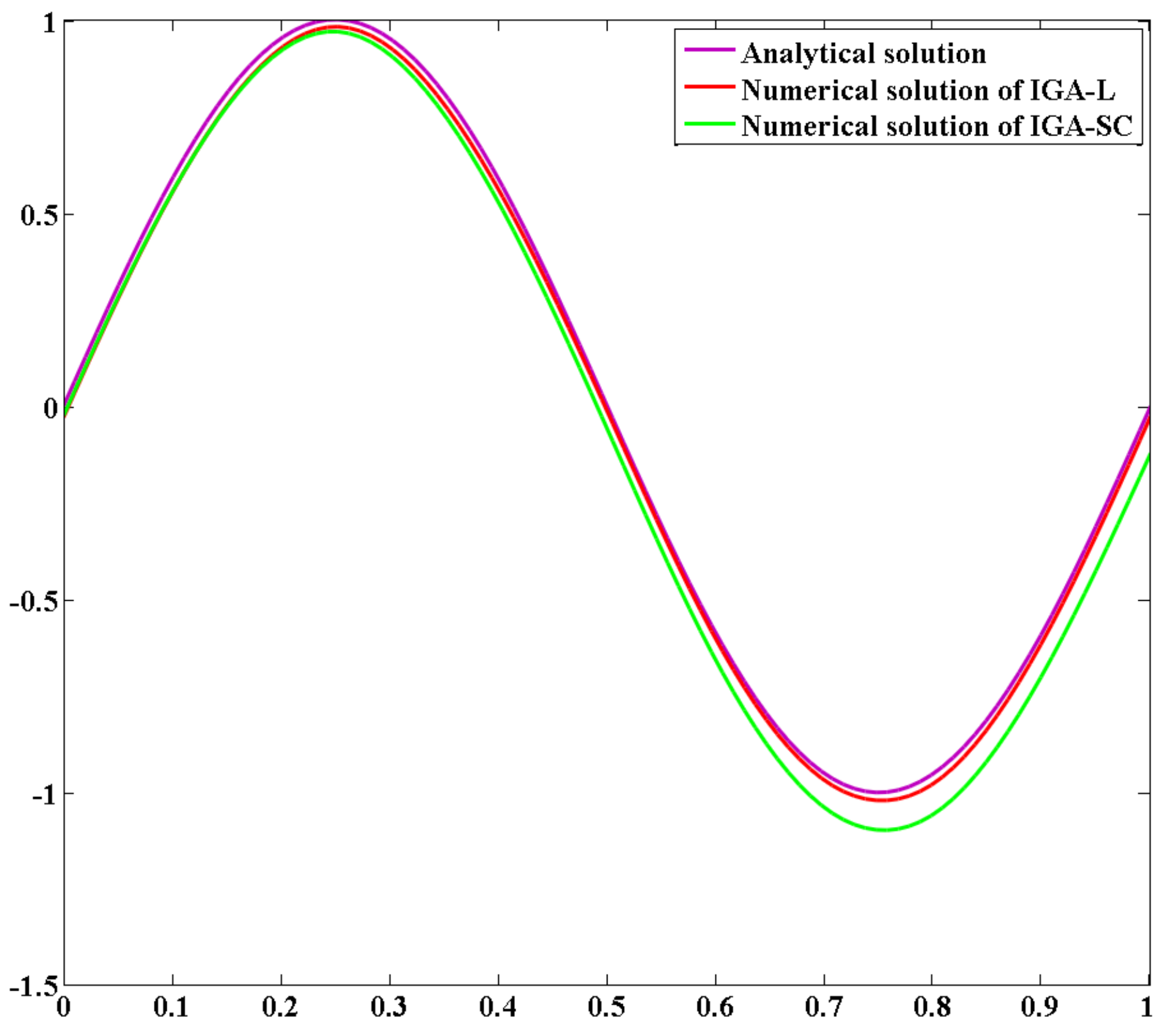}}
  \hspace{0.2cm}
  \subfigure[Absolute error curve~\pref{eq:abs_err} for the IGA-C solution of
            Eq.~\pref{eq:exmp_one_dim_nuemann} with Greville
            collocation points]{
    \label{fig:one-dim-igc-solution-g}
  \includegraphics[width = 0.46\textwidth]{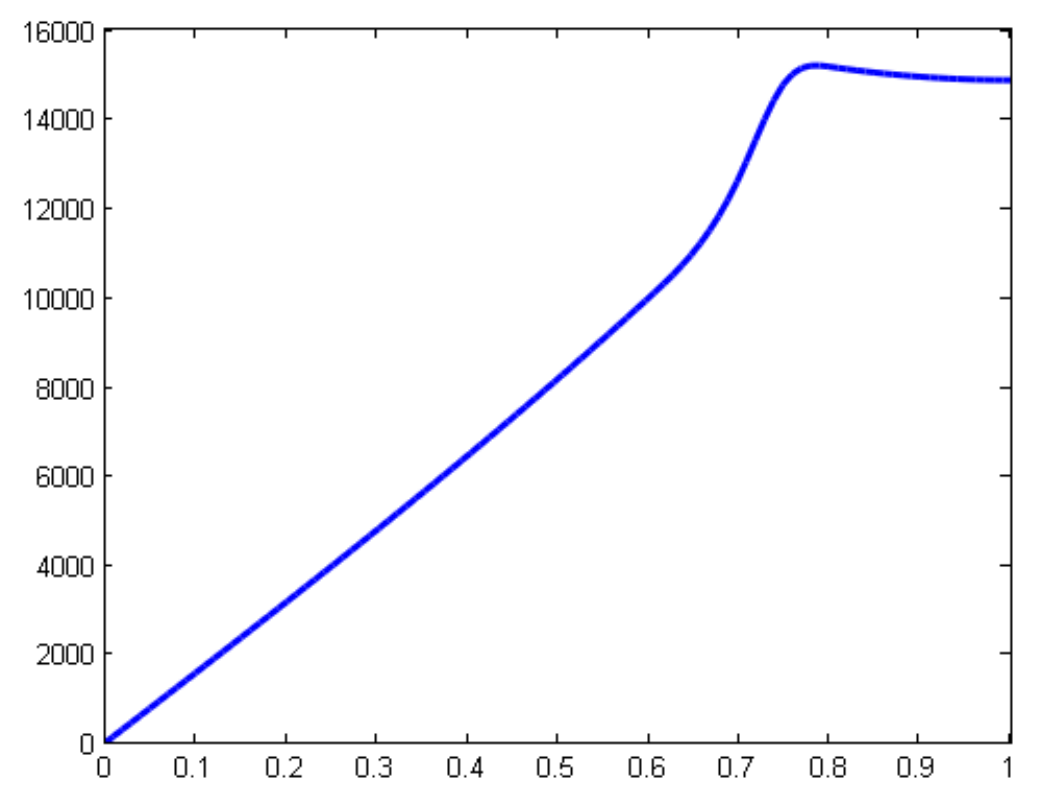}
  }
  \caption{IGA-L can be made more stable than IGA-C by changing the number of collocation points
    when solving Eq.~\pref{eq:exmp_one_dim_nuemann}.
    (a) The IGA-L solution to Eq.~\pref{eq:exmp_one_dim_nuemann} with
        $10$ control points and $14$ Greville collocation points is stable, with relative error $0.0343$.
        Meanwhile, the relative error of the IGA-SC solution to Eq.~\pref{eq:exmp_one_dim_nuemann} with $10$ control points and $14$ Greville collocation points is $0.0567$.
    (b) The IGA-C solution to Eq.~\pref{eq:exmp_one_dim_nuemann} with $10$
        control points and Greville collocation points is unstable, with relative error $2.6245 \times 10^3$.}
  \label{fig:igl-stable-than-igc}
\end{figure}

 \vspace{0.3cm}

 \textbf{Example V (Stability:):}
 Note that, in the IGA-C method the number of collocation points is
    fixed to be equal to the number of control points.
 However, in the IGA-L method,
    the number of collocation points is variable and larger than the
    number of control points.
 Therefore, the IGA-L method is more flexible than the IGA-C method.
 In this example, IGA-L, IGA-SC, and IGA-C are employed to solve a 1D source
    problem with Dirichlet boundary condition at the left end and Neumann boundary condition at the right end, i.e.,
    \begin{equation}
    \label{eq:exmp_one_dim_nuemann}
    \begin{cases}
    &-T'' + T = (1 + 4 \pi^2) \sin(2 \pi x),\ x \in \Omega = [0,1],\\
    &T(0) = 0, \\
    &T'(1) = 2 \pi. \\
    \end{cases}
    \end{equation}
 While the IGA-C method is unstable,
    the IGA-L method can be made stable by choosing appropriate number of collocation points.
 The analytical solution of Eq.~\pref{eq:exmp_one_dim_nuemann} is
    $T = \sin(2 \pi x)$.
 We still use the cubic B-spline curve, presented in Appendix A1,
    to represent the physical domain $\Omega = [0,1]$.

 Consider the case where the analytic solution of
    Eq.~\pref{eq:exmp_one_dim_nuemann} is approximated by a cubic
    B-spline function with $10$ control points,
    generated by inserting the following knots
    \begin{equation*}
        0.25,\ 0.5,\ 0.6,\ 0.7,\ 0.75,\ 0.8,
    \end{equation*}
    into the cubic B-spline curve in Appendix A1.
 When we use IGA-C to solve the source problem~\pref{eq:exmp_one_dim_nuemann}
    with Greville collocation points,
    the numerical solution is unstable,
    with relative error~\pref{eq:rel_err_func} $2.6245 \times 10^3$ (Fig.~\ref{fig:one-dim-igc-solution-g}).
 However, when IGA-L is employed to solve the source
    problem~\pref{eq:exmp_one_dim_nuemann} with $14$ Greville collocation points,
    the solution is stable, with relative error $0.0343$ (Fig.~\ref{fig:one-dim-igl-solution-g}).
 In addition, though the IGA-SC solution with $14$ Greville collocation
    points to Eq.\pref{eq:exmp_one_dim_nuemann} is also stable,
    its relative error is $0.0567$,
    larger than that of the IGA-L solution.

 \begin{figure}[!htb]
\centering
  \subfigure[Physical domain.]{ \label{subfig:fc_domain}
    \includegraphics[width = 0.28\textwidth]{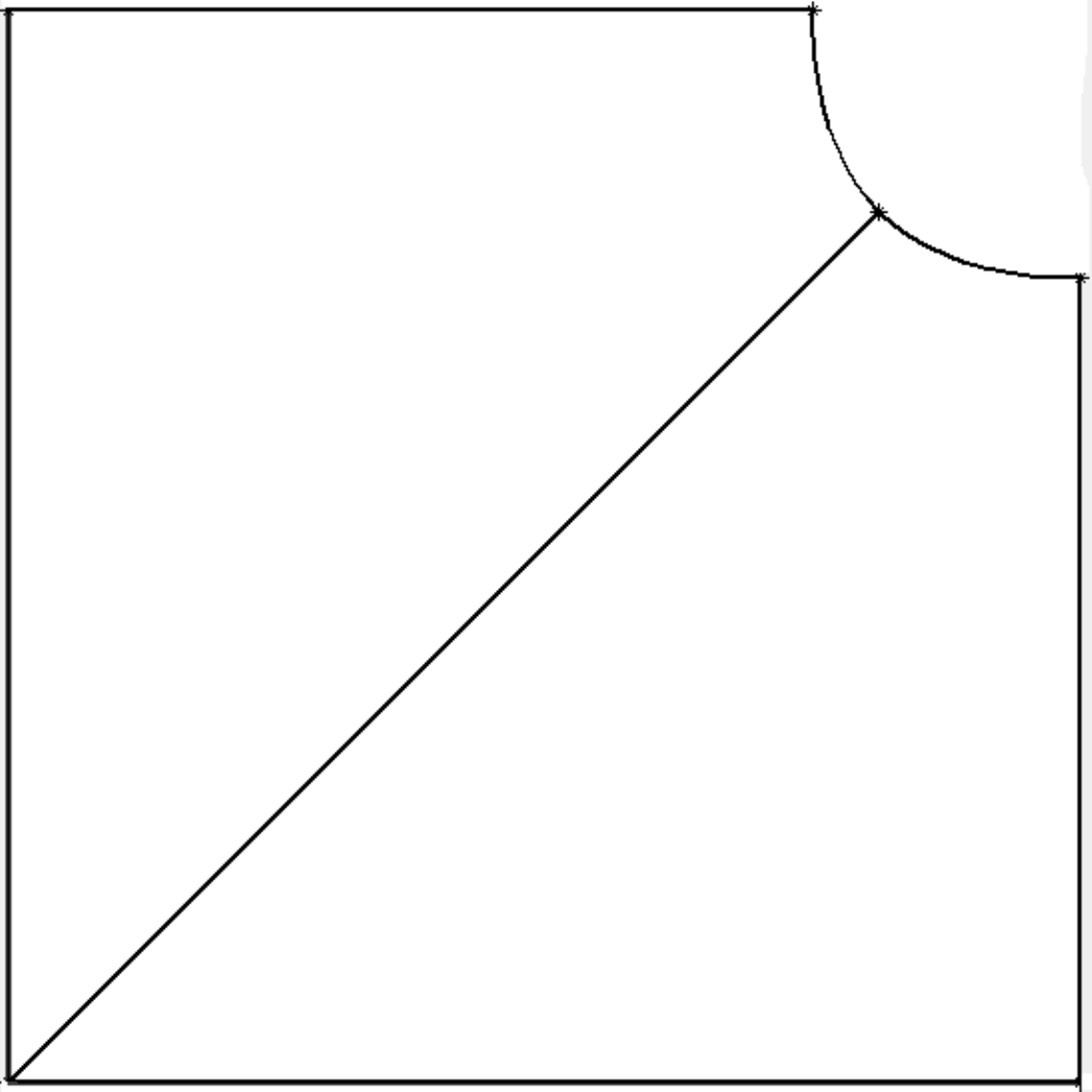}}
  \subfigure[Control net.]{ \label{subfig:fc_control_net}
    \includegraphics[width=0.28\textwidth]{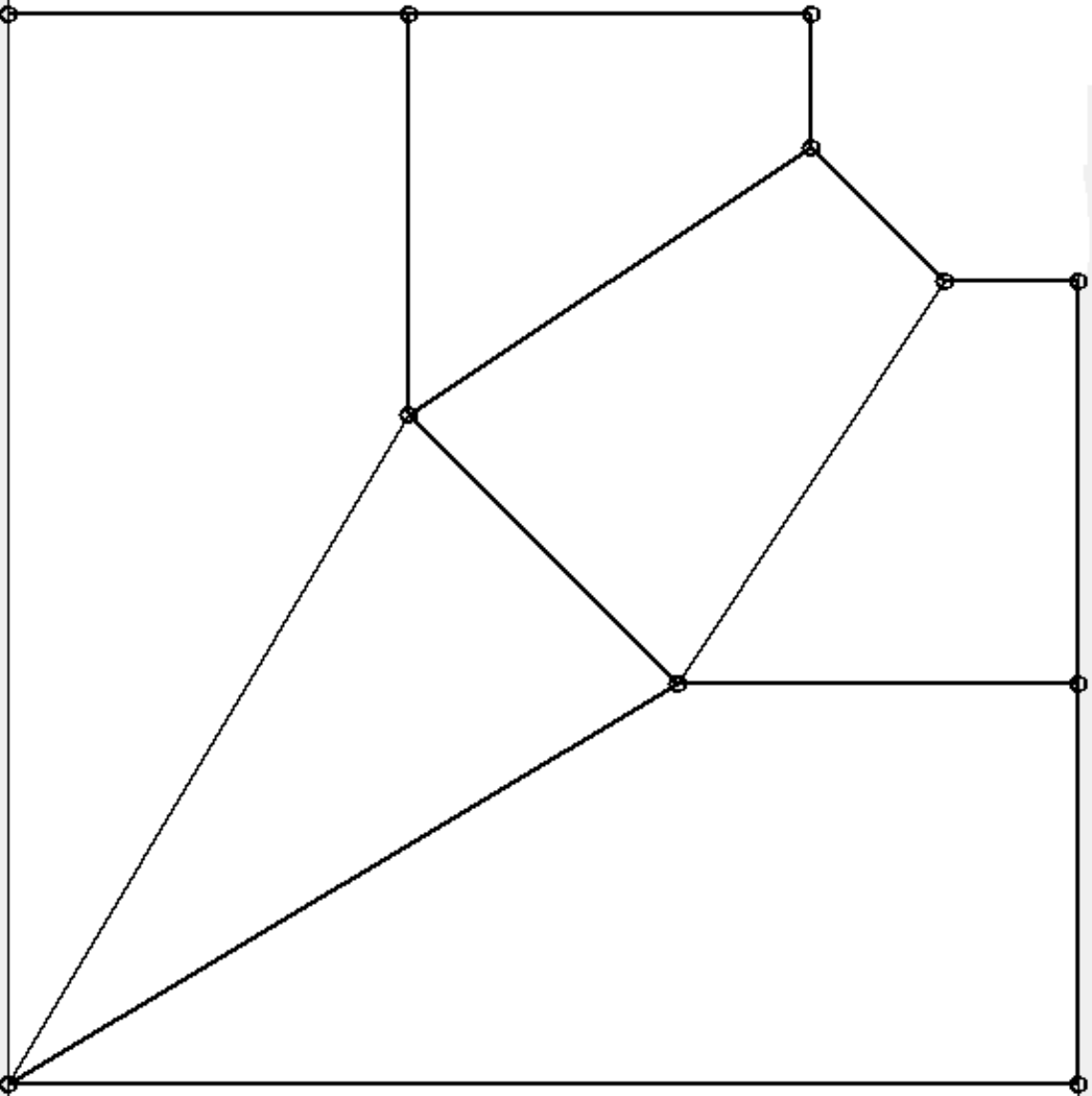}}
  \subfigure[Analytical solution.]{ \label{subfig:fc_ana_solution}
    \includegraphics[width = 0.36\textwidth]
        {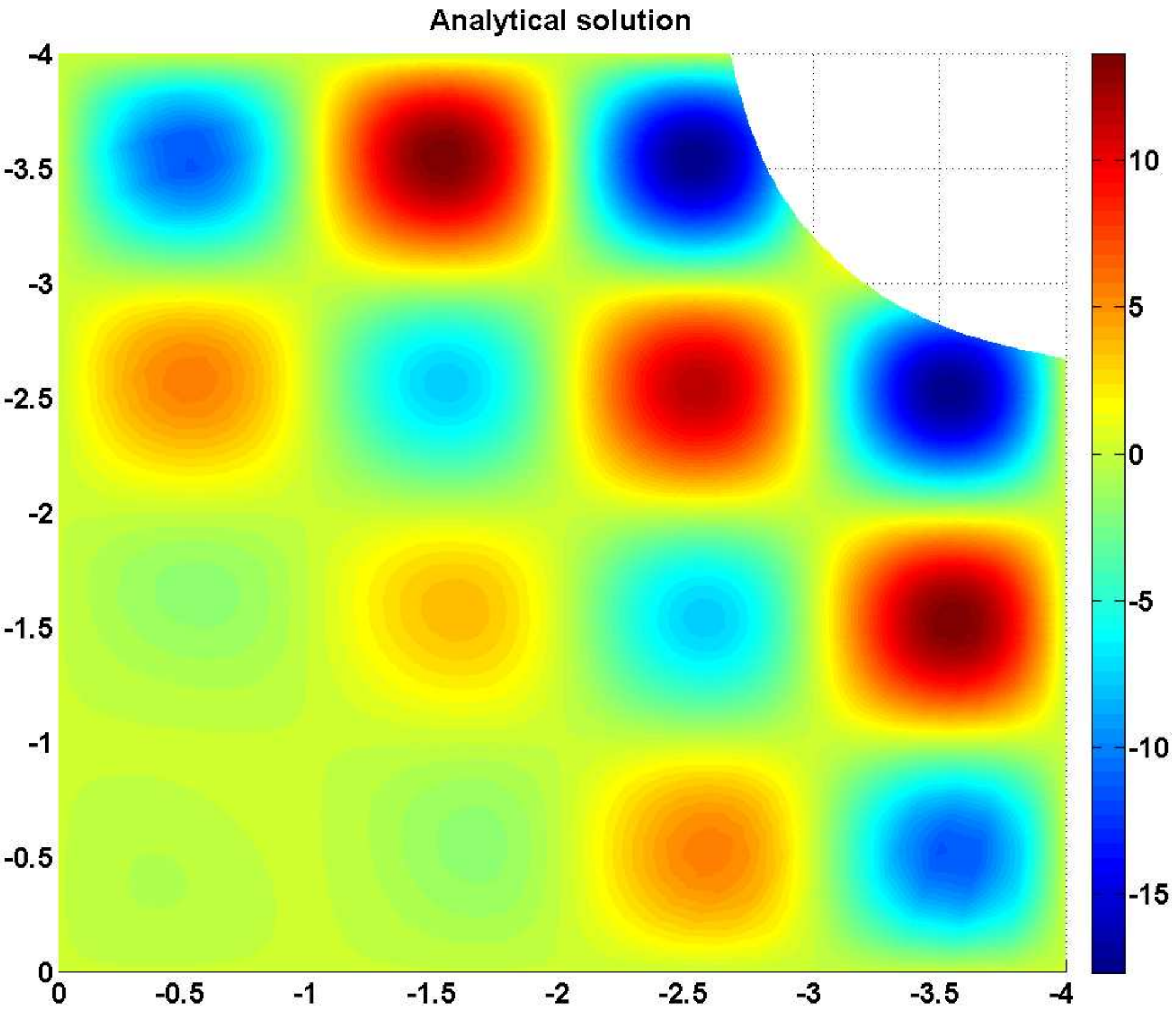}}
  \caption{Physical domain (a), control net (b), and the analytical solution (c) on the domain of frame corner.}
  \label{fig:frame_corner_domain}
\end{figure}

 \vspace{0.3cm}

 \textbf{Example VI (Frame Corner):}
 In this example, we solve a 2D source problem defined on the domain of frame
    corner (Fig.~\ref{subfig:fc_domain}), i.e.,
    \begin{equation}\label{eq:fc_2d_problem}
      \begin{cases}
        & - \triangle T + T = f, \quad (x,y) \in \Omega \\
        & T|_{\partial \Omega} = 0,
      \end{cases}
    \end{equation}
    where
    $$ f = (-4 + (2 \pi^2 + 1) x^2 + (2 \pi^2 + 1) y^2 - (2 \pi^2 + 1))
           \sin(\pi x) \sin(\pi y) - 4 \pi x \cos(\pi x) \sin(\pi y) -
           4 \pi y \cos(\pi y) \sin(\pi x). $$
 The analytical solution of the 2D source problem
    (Eq.~\pref{eq:fc_2d_problem}) is (Fig~\ref{subfig:fc_ana_solution}),
 $$ T = (x^2 + y^2 - 1) \sin(\pi x) \sin(\pi y).$$

 The physical domain of the 2D source problem~\pref{eq:fc_2d_problem} is a
    frame-corner-like shape (see Fig.~\ref{subfig:fc_domain}),
    which is modeled by a bi-quadratic NURBS surface with two patches (see Appendix A5).
 Fig.~\ref{subfig:fc_control_net} illustrates the control net of the
    bi-quadratic NURBS surface,
    where there are two overlapping control points at the lower left corner.
 So, the two patches are $C^0$ continuous across their common boundary inside
    the domain (Fig.~\ref{subfig:fc_domain}).
 In other words, the physical domain contains a $C^0$ line
    (Fig.~\ref{subfig:fc_domain}).

 \begin{figure}[!htb]
\centering
  \subfigure[IGA-L solution of Eq.~\pref{eq:fc_2d_problem} with relative error $7.84 \times 10^{-4}$.]{ \label{subfig:fc_iga_l_solution}
    \includegraphics[width=0.4\textwidth]{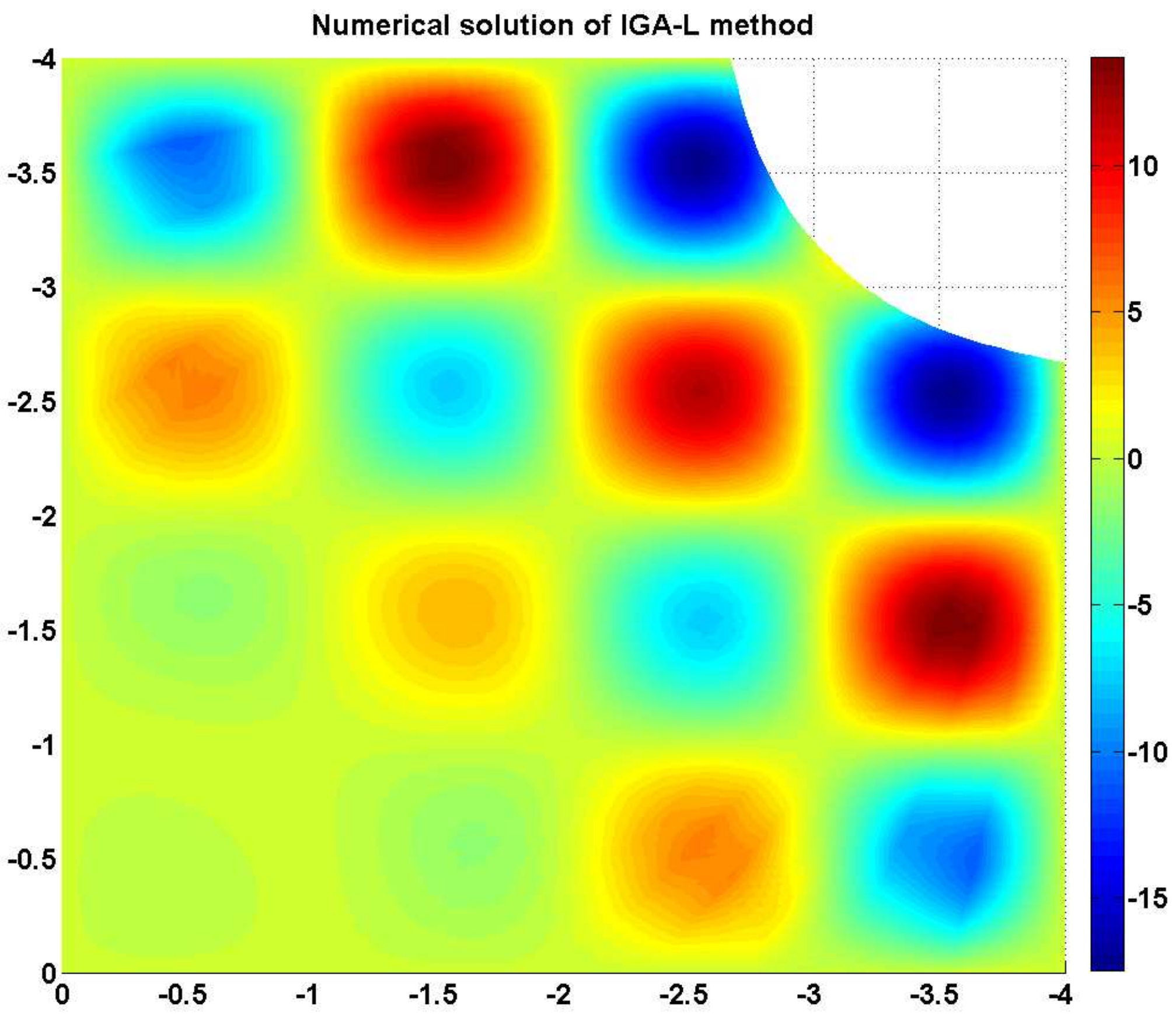}}
  \subfigure[IGA-G solution of Eq.~\pref{eq:fc_2d_problem} with relative error $4.25 \times 10^{-4}$.]{ \label{subfig:fc_iga_g_solution}
    \includegraphics[width = 0.4\textwidth]
        {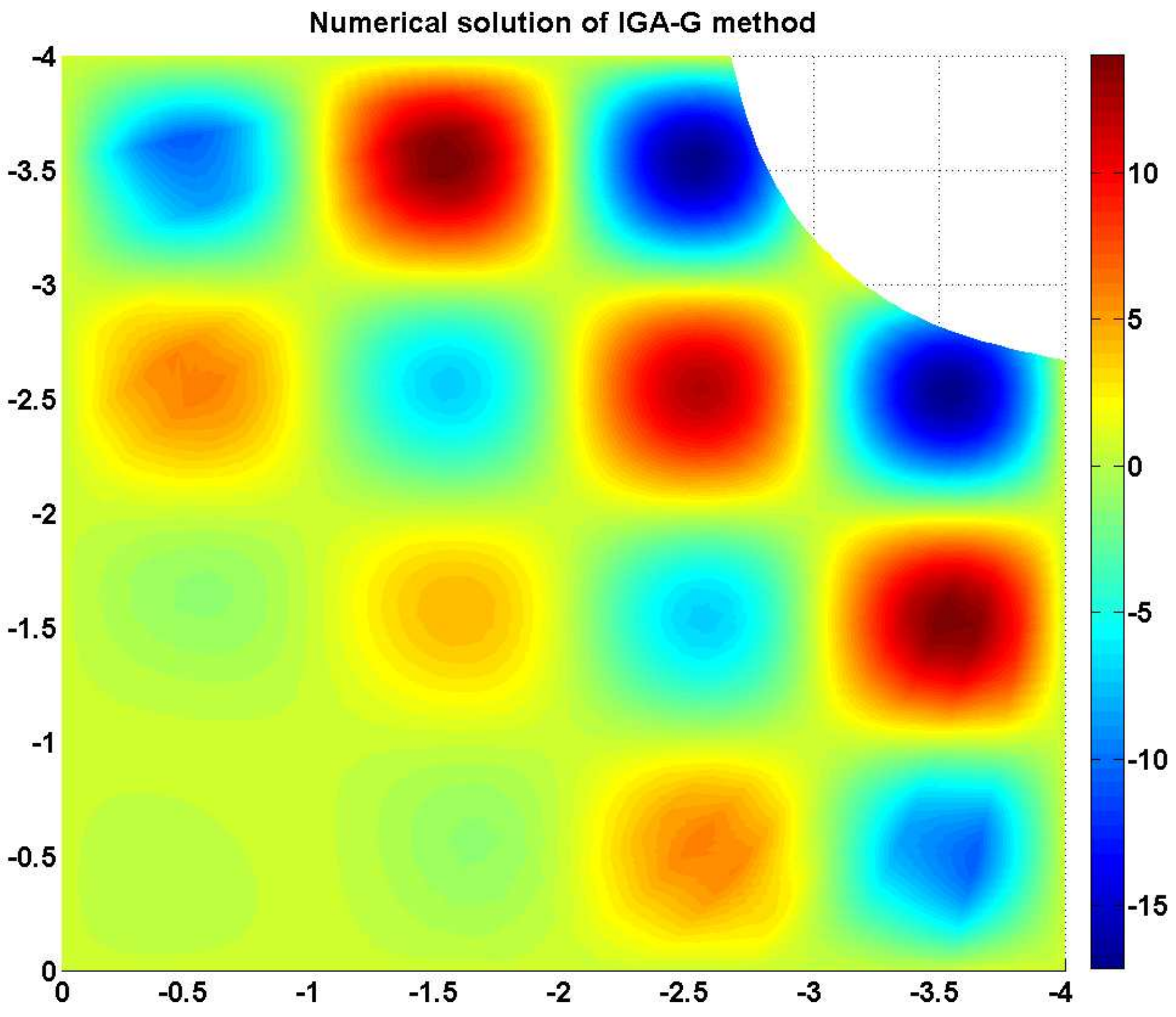}}
  \subfigure[IGA-C solution of Eq.~\pref{eq:fc_2d_problem} with relative error $5.42 \times 10^{-2}$.]{ \label{subfig:fc_iga_c_solution}
    \includegraphics[width = 0.4\textwidth]
        {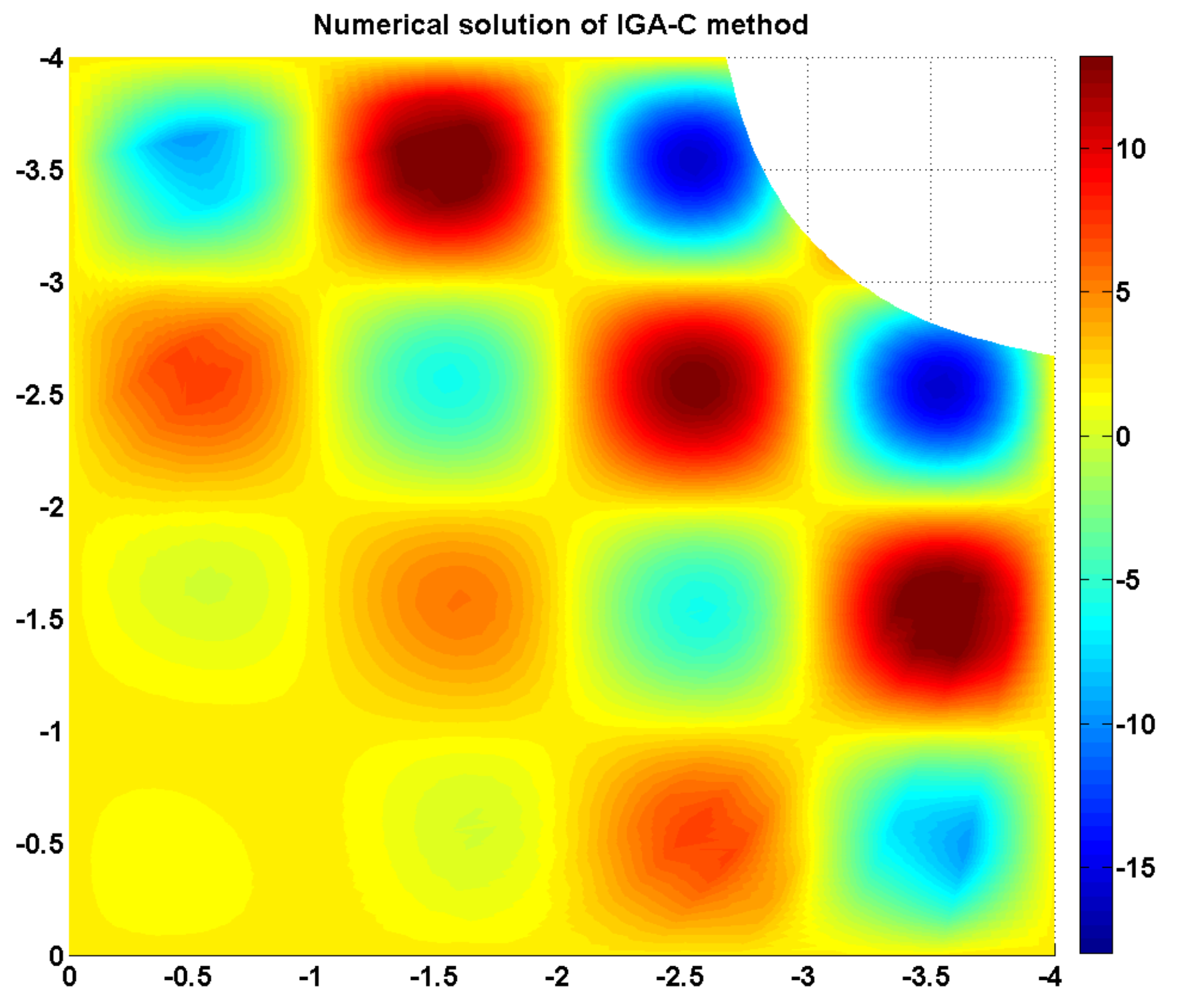}}
  \subfigure[IGA-SC solution of Eq.~\pref{eq:fc_2d_problem} with relative error $8.53 \times 10^{-4}$.]{ \label{subfig:fc_iga_sc_solution}
    \includegraphics[width = 0.4\textwidth]
        {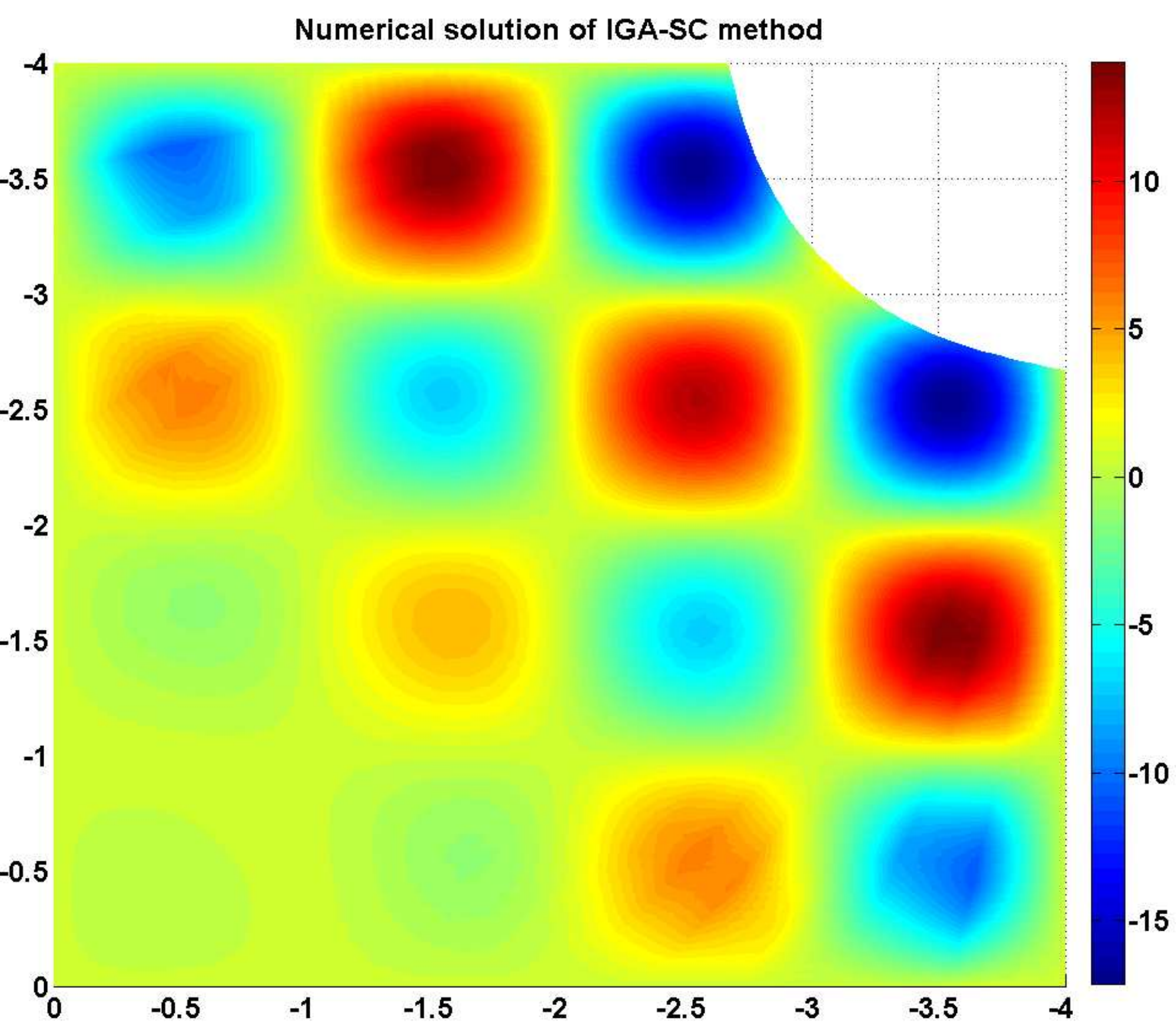}}
  \caption{The numerical solutions by IGA-L (a), IGA-G (b), IGA-C (c), and IGA-SC (d) methods, respectively.}
  \label{fig:frame_corner_numerical_solutions}
\end{figure}

 In this example, we compare the IGA-L method with IGA-C, IGA-SC,
    and isogeometric Galerkin method (IGA-G).
 Fig.~\ref{fig:frame_corner_numerical_solutions} illustrates the numerical
    solutions generated by IGA-L, IGA-G, IGA-C, and IGA-SC methods.
 All of the numerical solutions are bi-cubic NURBS functions with
    $17 \times 18$ control points.
 The relative error of IGA-C solution is $5.42 \times 10^{-2}$
    (Fig.\ref{subfig:fc_iga_c_solution}).
 Using $64 \times 65$ quadrature points,
    the relative error of IGA-G solution is $4.25 \times 10^{-4}$ (Fig.\ref{subfig:fc_iga_g_solution}).
 With $28 \times 30$ collocation points,
    the relative error of IGA-SC method is $8.53 \times 10^{-4}$.
 However, the relative error of our IGA-L method reaches
    $7.84 \times 10^{-4}$ using $24 \times 26$ collocation points.
 Note that, the number of collocation points of our IGA-L method is smaller
    than that of IGA-SC method,
    but the precision of IGA-L method is better than that of IGA-SC method.

\begin{figure}[!htb]
\centering
  \subfigure[]{ \label{subfig:fc_iga_l_error}
    \includegraphics[width=0.4\textwidth]{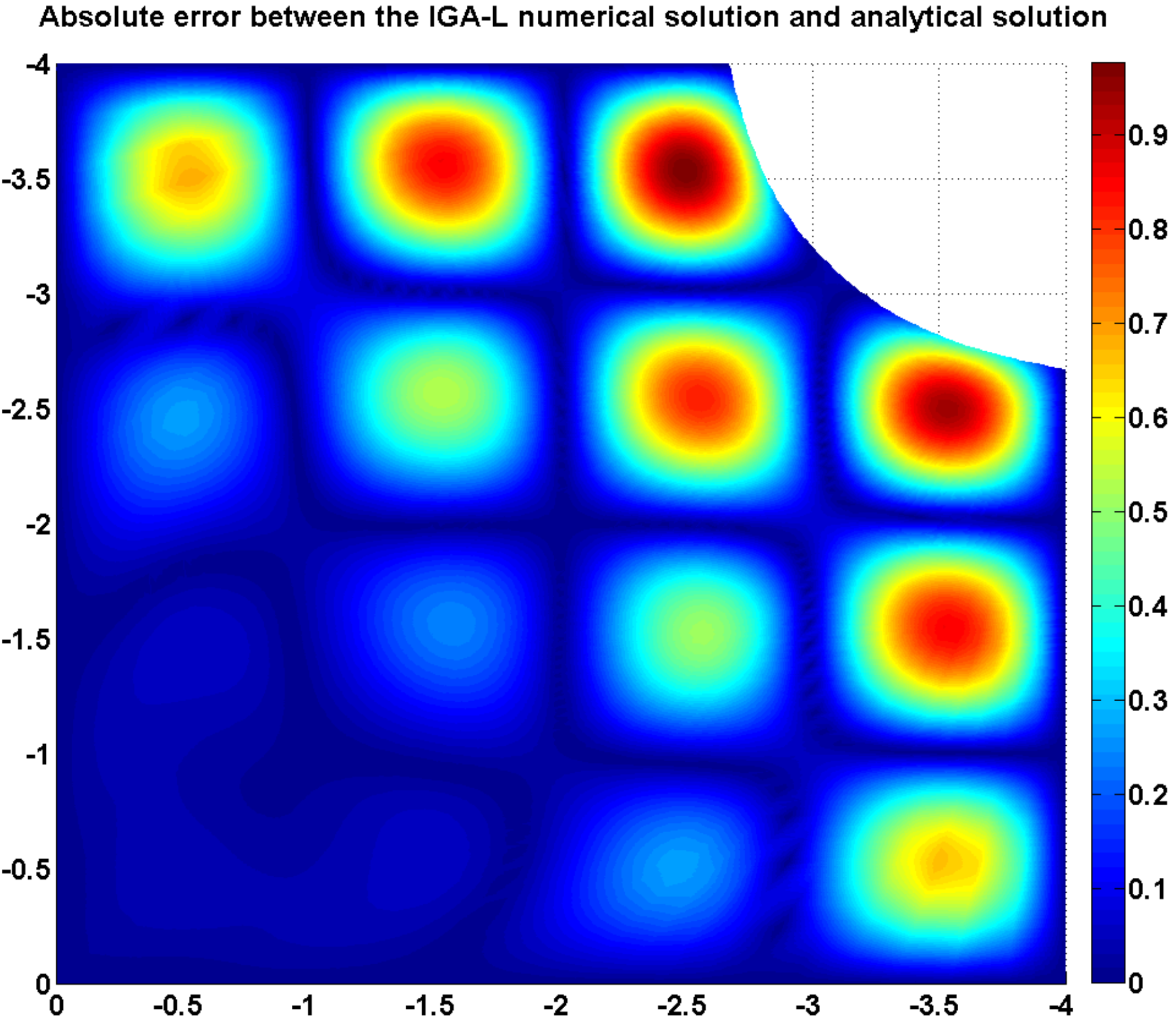}}
  \subfigure[]{ \label{subfig:fc_iga_g_error}
    \includegraphics[width = 0.4\textwidth]{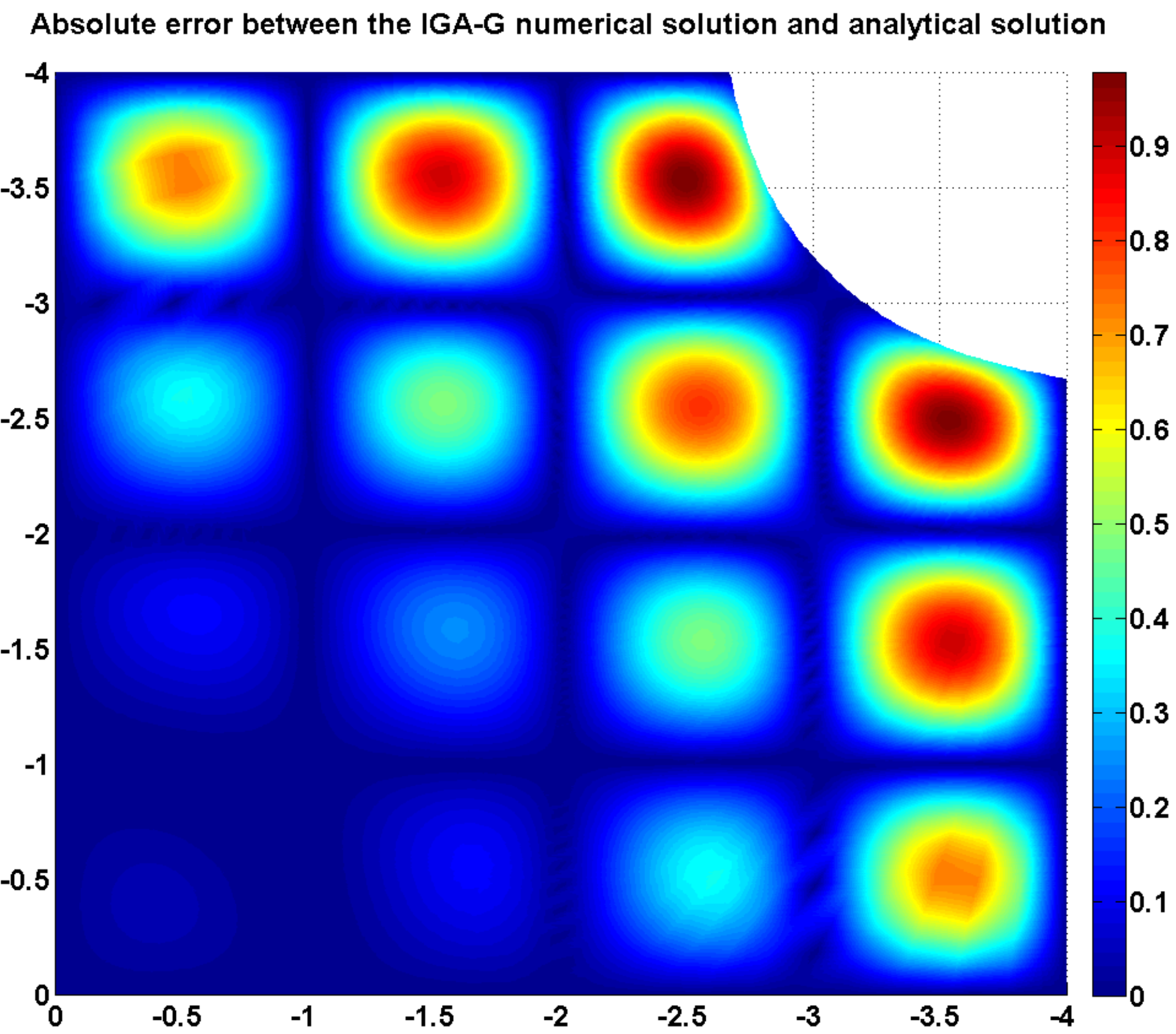}}
  \subfigure[]{ \label{subfig:fc_iga_c_error}
    \includegraphics[width = 0.4\textwidth]{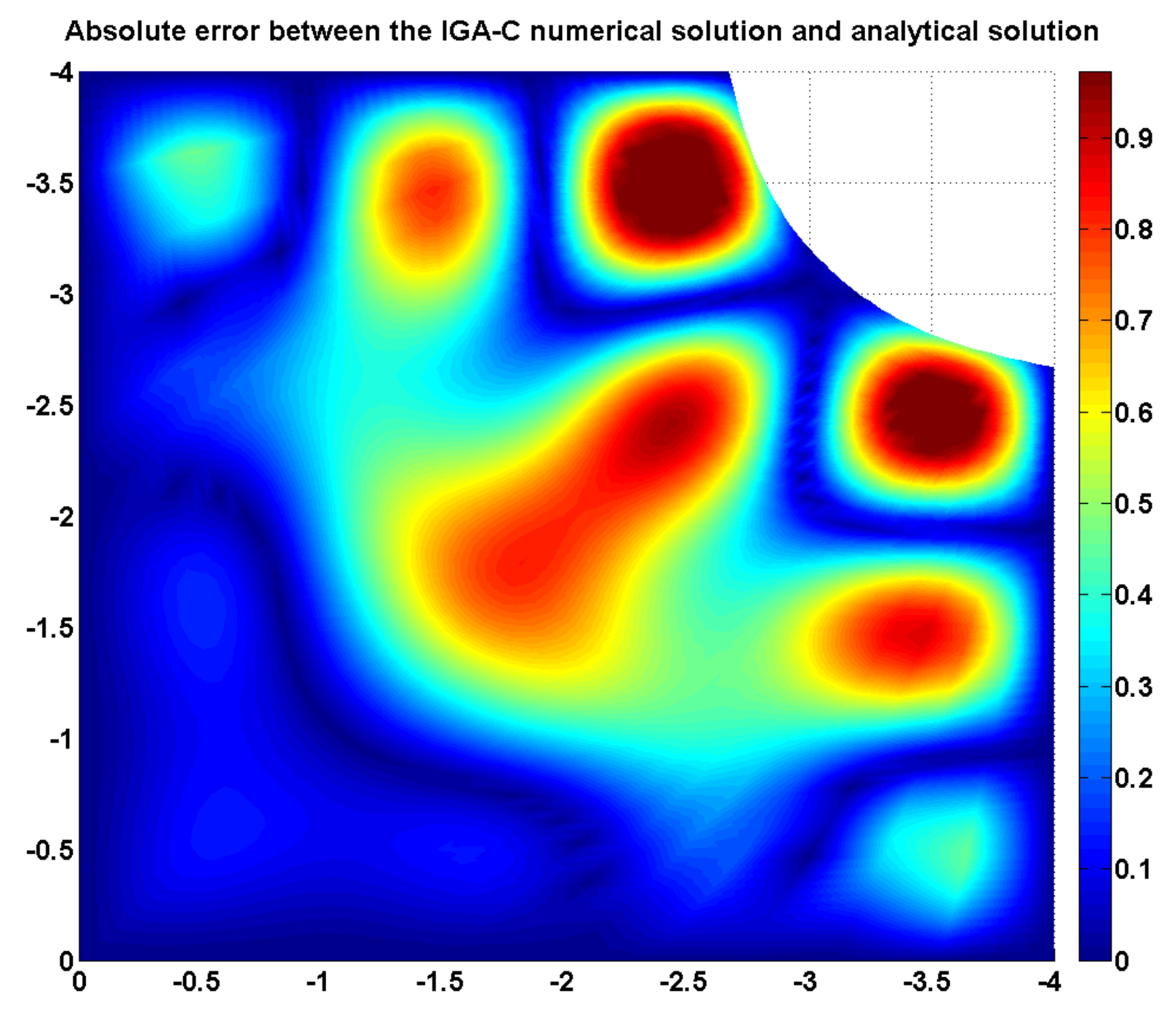}}
  \subfigure[]{ \label{subfig:fc_iga_sc_error}
    \includegraphics[width = 0.4\textwidth]{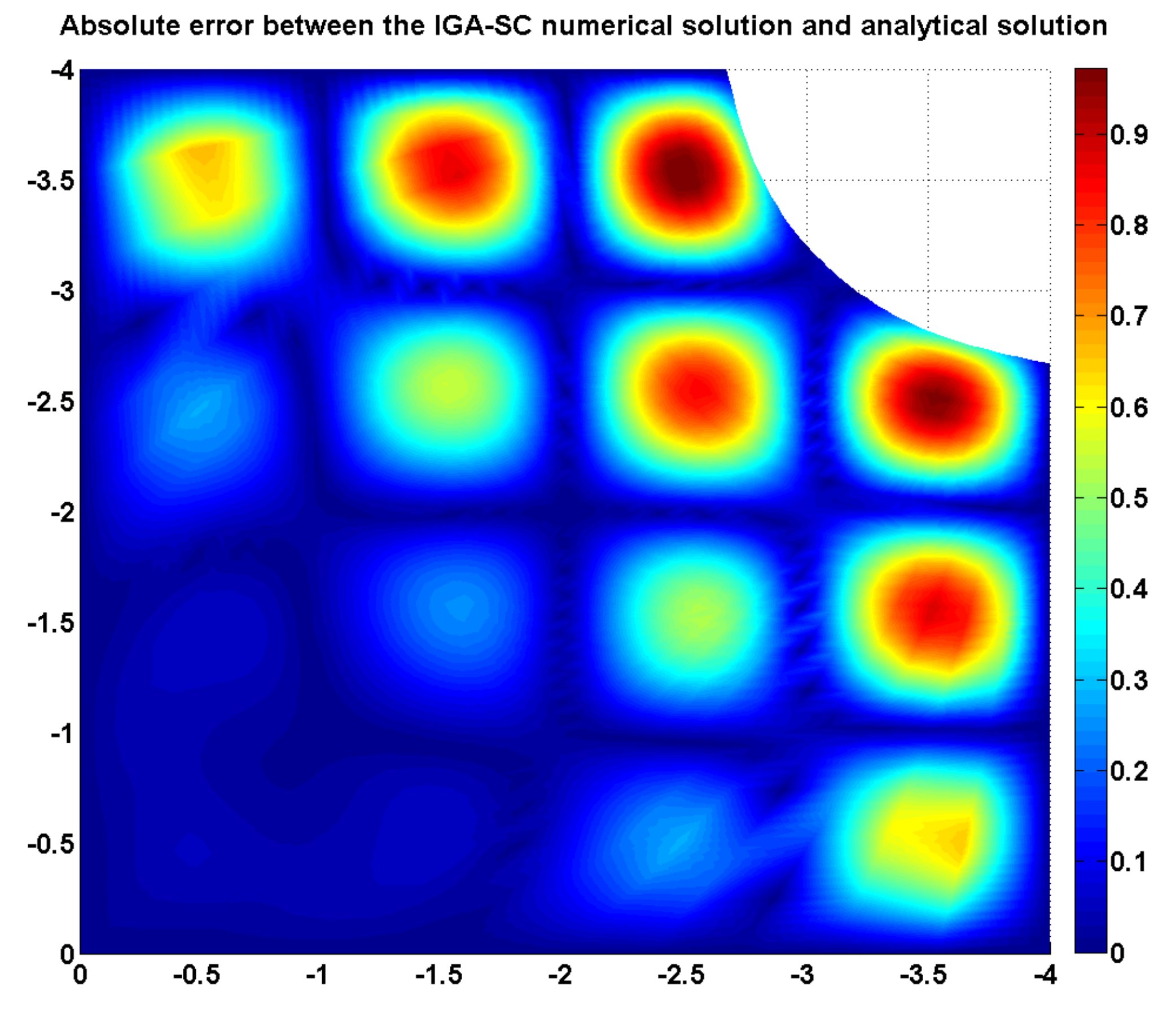}}
  \caption{Absolute error distributions for IGA-L solution (a), IGA-G solution (b), IGA-C solution (c), and IGA-SC solution (d), respectively.}
  \label{fig:frame_corner_error_distribution}
\end{figure}

 Moreover, Fig.~\ref{fig:frame_corner_error_distribution} demonstrates the
    absolute error distributions of IGA-L, IGA-G, IGA-C, and IGA-SC solutions.
 It can be noticed that,
    while the absolute error distribution of IGA-C method is heavily influenced by the $C^0$ line of the the physical domain (Fig.~\ref{subfig:fc_iga_c_error}),
    the $C^0$ line almost does not affect the absolute error distribution of IGA-L method (Fig.~\ref{subfig:fc_iga_l_error}).
 Moreover, compared with the IGA-SC method
    (Fig.~\ref{subfig:fc_iga_sc_error}),
     the absolute error distribution of IGA-L method (Fig.~\ref{subfig:fc_iga_l_error}) is closer to that of IGA-G method (Fig.~\ref{subfig:fc_iga_g_error}).

\begin{figure}[!htb]
\centering
  \subfigure[]{ \label{subfig:fc_diagram_iga_l}
    \includegraphics[width=0.42\textwidth]{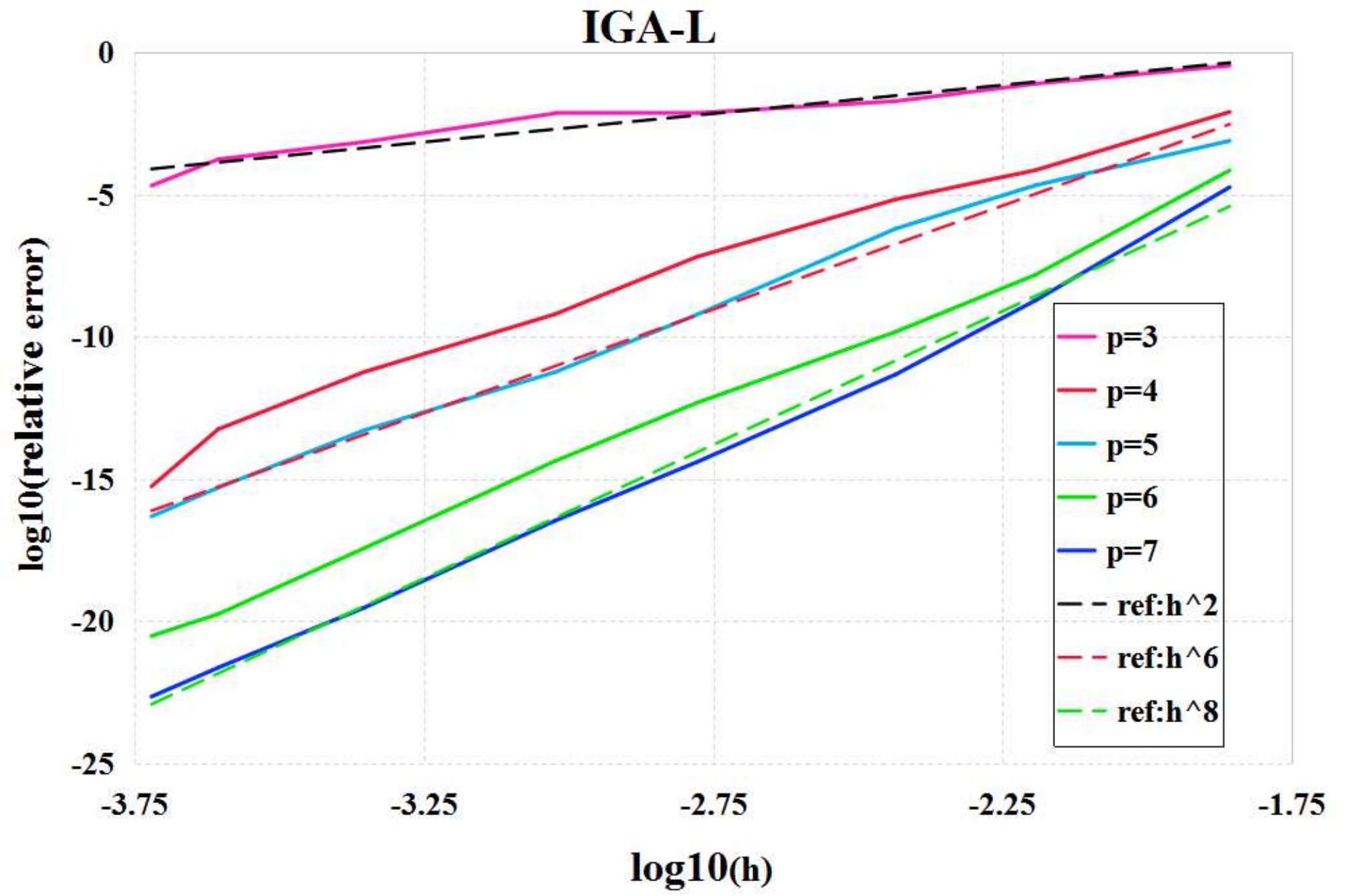}}
  \subfigure[]{ \label{subfig:fc_diagram_iga_g}
    \includegraphics[width = 0.4\textwidth]{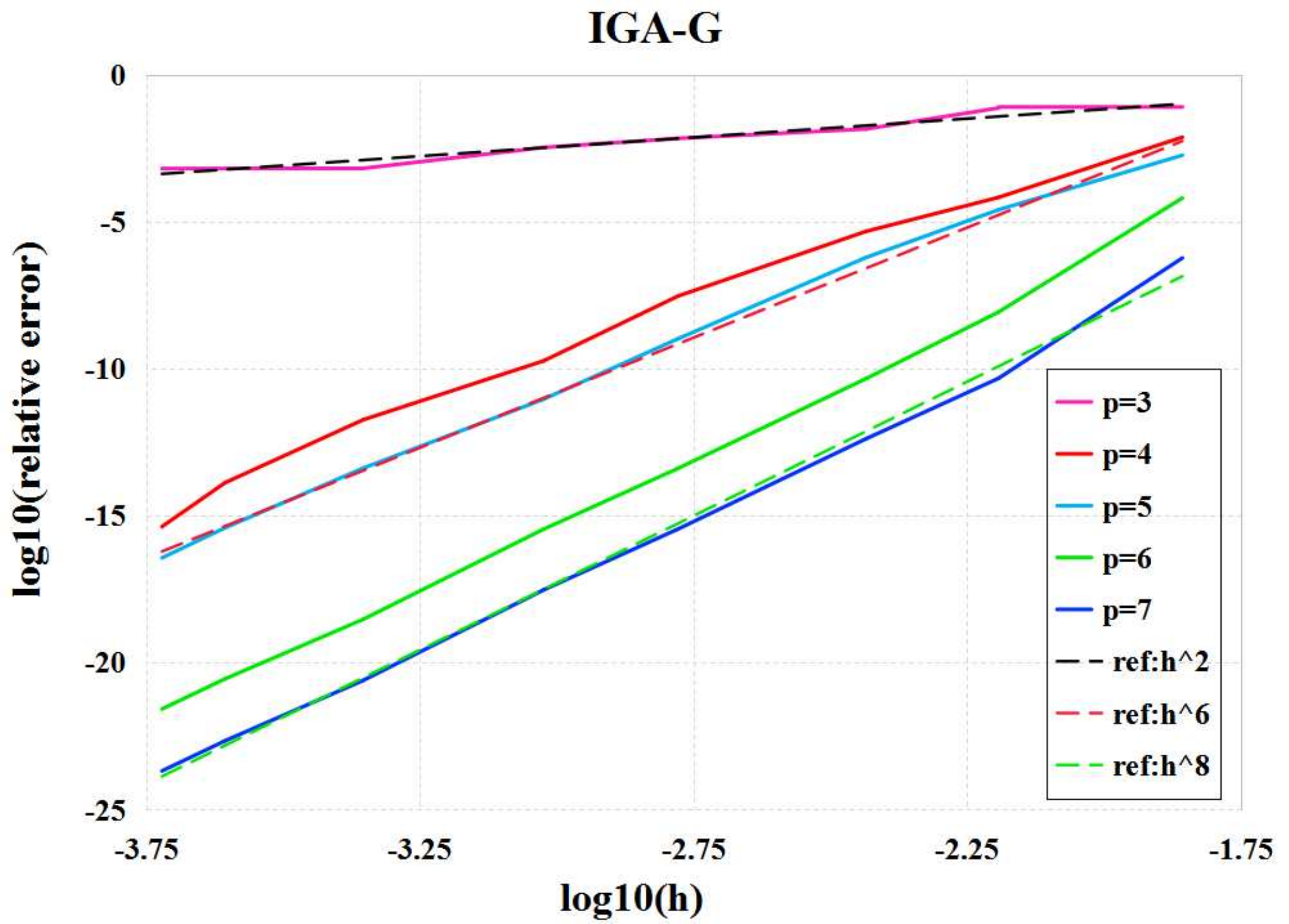}}
  \subfigure[]{ \label{subfig:fc_diagram_iga_c}
    \includegraphics[width = 0.4\textwidth]{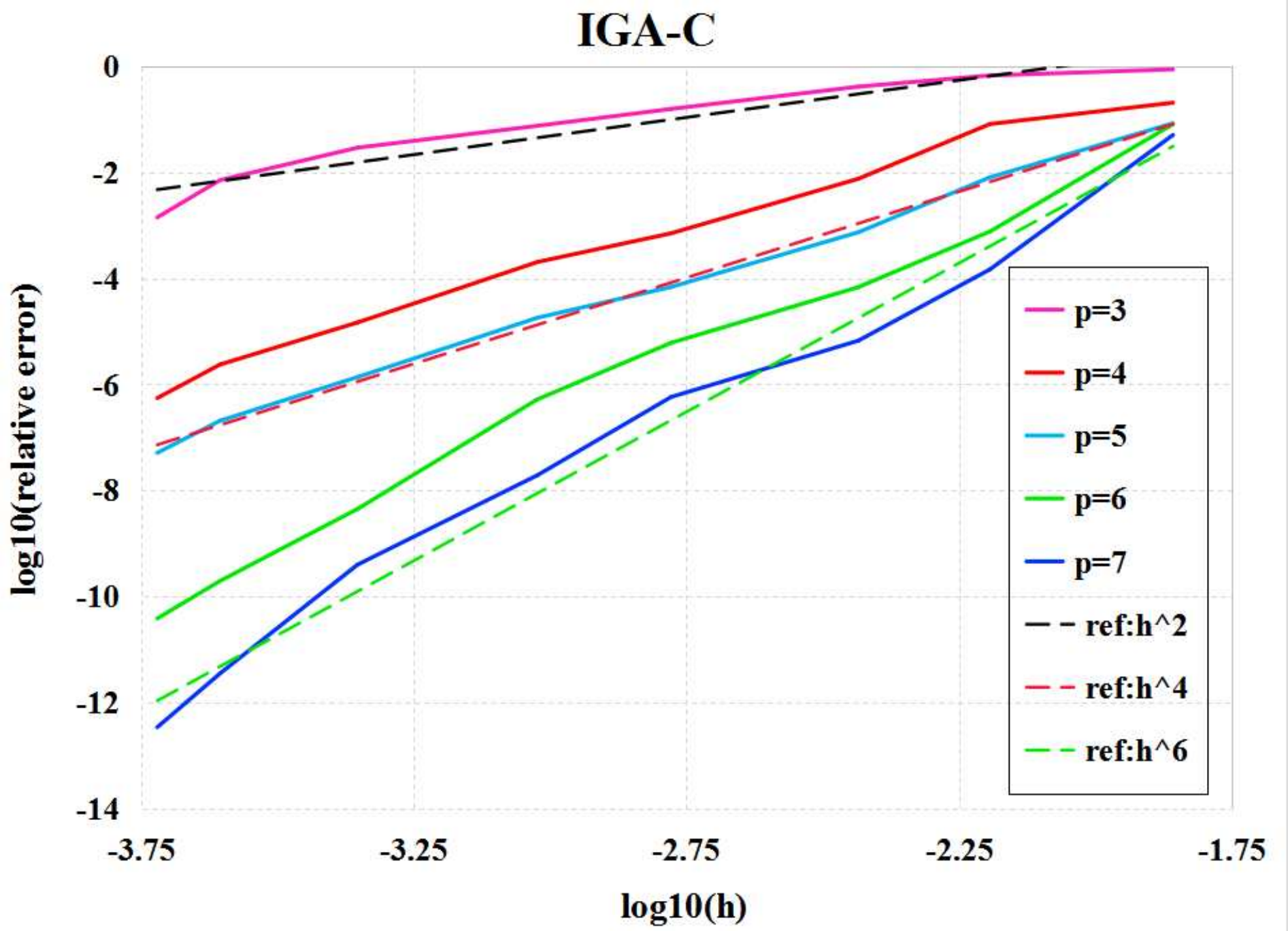}}
  \subfigure[]{ \label{subfig:fc_diagram_iga_sc}
    \includegraphics[width = 0.4\textwidth]{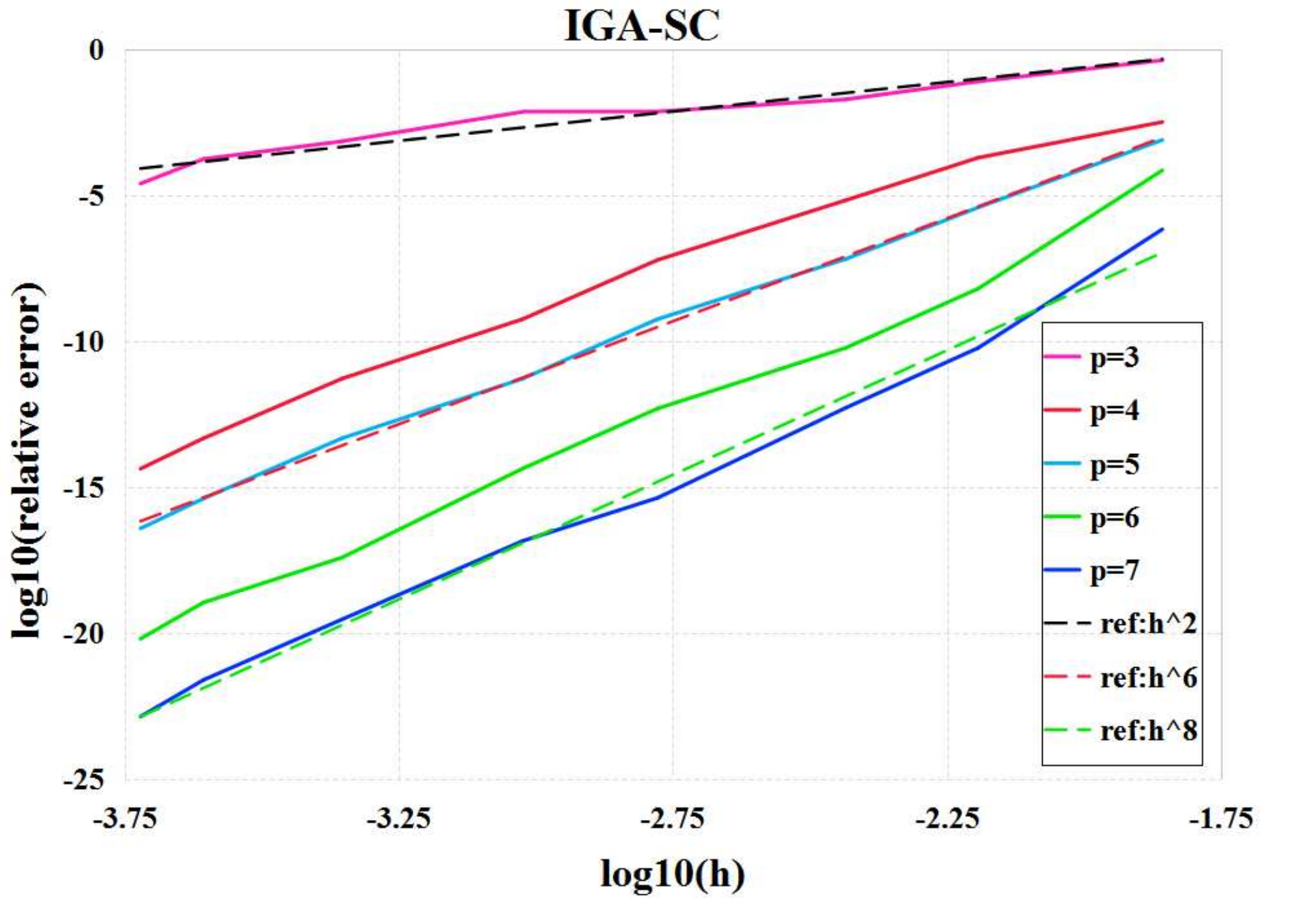}}
  \subfigure[]{ \label{subfig:fc_time_error}
    \includegraphics[width = 0.4\textwidth]{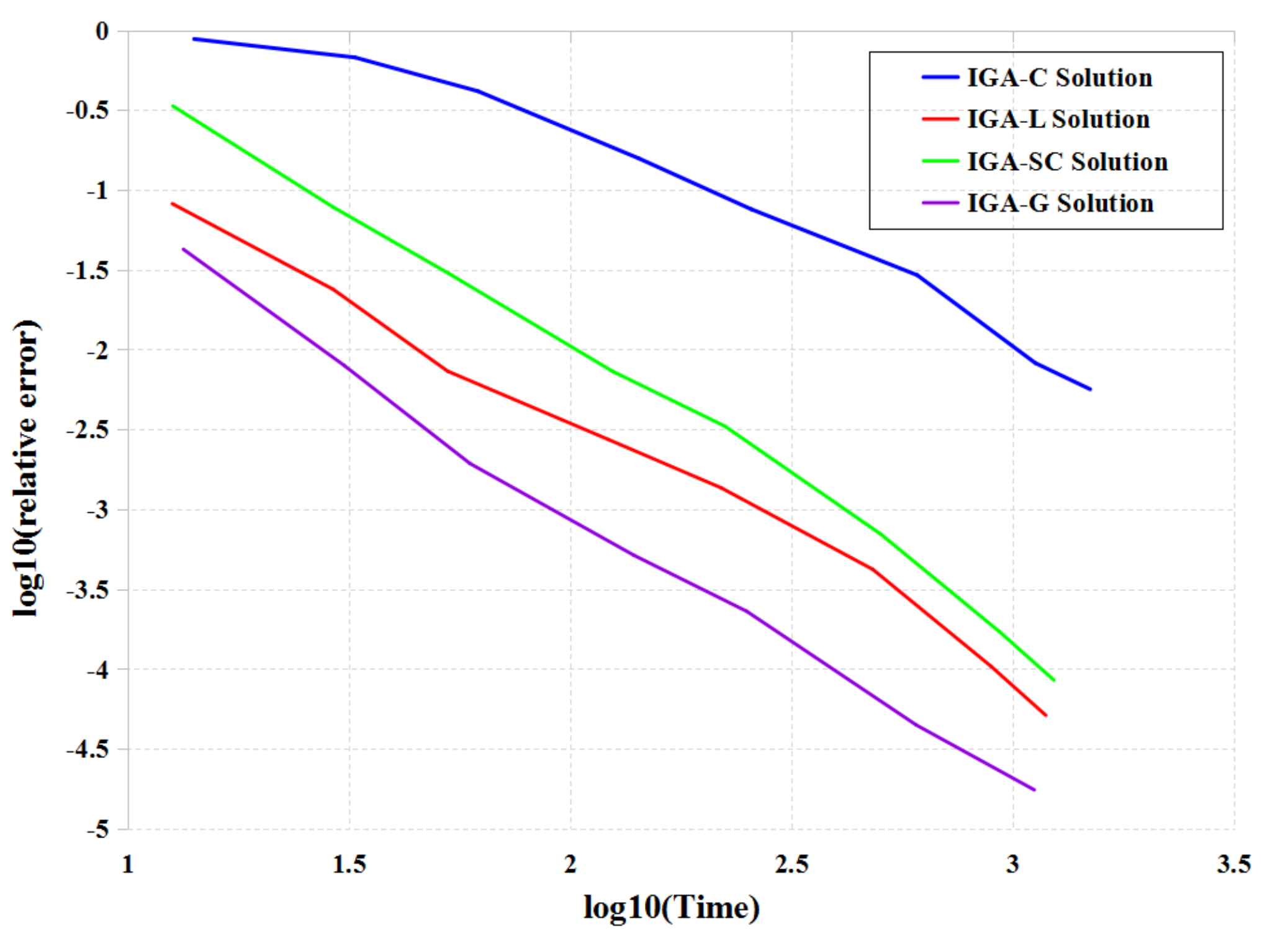}}
  \caption{Diagrams of $\log_{10}(h)$ v.s. $\log_{10}$(relative error) for
        IGA-L (a), IGA-G (b), IGA-C (c), and IGA-SC (d) methods, respectively.
  And, diagram of $\log_{10}$(time) v.s. $\log_{10}$(relative error) for the
    four methods (e). }
  \label{fig:frame_corner_diagram}
\end{figure}

 Finally, in Fig.~\ref{fig:frame_corner_diagram},
    we present diagrams of $\log_{10}(h)$ v.s. $\log_{10}$(relative error), and diagram of $\log_{10}$(time) v.s. $\log_{10}$(relative error)
    for IGA-L, IGA-G, IGA-C, and IGA-SC methods.
 From diagrams of $\log_{10}(h)$ v.s. $\log_{10}$(relative error)
    (Figs.\ref{subfig:fc_diagram_iga_l}-\ref{subfig:fc_diagram_iga_sc}),
    we can see that the convergence rates of IGA-L, IGA-G, and IGA-SC are the same, i.e., $O(h^2)$ for degree $p=3$, $O(h^6)$ for degrees $p=4$ and $5$,
    $O(h^8)$ for degree $p=6$ and $7$.
 Lastly, from diagram of $\log_{10}$(time) v.s. $\log_{10}$(relative error),
    IGA-L method is better than IGA-SC and IGA-C methods.

\section{Conclusion}
\label{sec:conclusion}

 IGA approximates the solution of a boundary value
    problem (or initial value problem) by a NURBS function.
 In this paper, we developed the IGA-L method to determine the unknown coefficients of the approximate
    NURBS function
    by fitting the sampling values in the
    least-squares sense.
 We proved the consistency and convergence properties of IGA-L.
 Moreover, the many numerical examples presented in this paper show that
    a small computational increase in IGA-L leads to
    large improvements in the accuracy,
    and furthermore that IGA-L is more flexible and more stable than IGA-C.

\section*{Acknowledgement} This work is supported by the National Natural Science Foundation of China
    (Grant Nos. 61379072, 61202201, 60970150).
    Dr. Qianqian Hu is also supported by the Open Project
    Program (No. A1305) of the State Key Lab of CAD\&CG, Zhejiang
    University.


\section*{Appendix}

 In Appendix A1-A5,
    we list the control points, knot vector, and weights of the
    NURBS representation of the physical domains in the numerical examples.
 In Appendix A6, the proof to the formula~\pref{eq:two_dim_error} in the two
    dimensional case in Theorem~\ref{thm:error_formula} is presented.


\subsection*{A1: NURBS representation of the physical domain in
Examples I and IV}

 The physical domain in Examples I and IV is represented by a cubic
    B-spline curve.
 Its control points are listed in the following
    Table~\ref{tbl:one-dim-ctrlpnt}.


\begin{table}[!htb]
\caption{Control points of the cubic B-spline curve in Examples I
    and IV} \label{tbl:one-dim-ctrlpnt} \centering
\begin{tabular}{c  c  c  c}
 \hline
$B_1$ & $B_2$ & $B_3$ & $B_4$ \\
 \hline
 0  & $\frac{1}{3}$ & $\frac{2}{3}$ & 1 \\
 \hline
\end{tabular}
\end{table}

The knot vector is
\begin{equation*}
    0\ 0\ 0\ 0\ 1\ 1\ 1\ 1.
\end{equation*}


\subsection*{A2: NURBS representation of the physical domain in
Example II}

The physical domain in Example II is represented by a cubic NURBS
    patch.
Its control points are listed in Table~\ref{tbl:two-dim-ctrlpnt},
    and Table~\ref{tbl:two-dim-weight-ctrlpnt} presents its weights.
The knot vectors along $u-$ and $v-$direction are, respectively,
\begin{equation*}
\begin{split}
0\ 0\ 0\ 0\ 1\ 1\ 1\ 1, \\
0\ 0\ 0\ 0\ 1\ 1\ 1\ 1.
\end{split}
\end{equation*}


\begin{table}[!htb]
\caption{Control points of the quarter of annulus}
\label{tbl:two-dim-ctrlpnt} \centering
\begin{tabular}{c  c  c  c  c}
 \hline
$i$ & $B_{i,1}$ & $B_{i,2}$ & $B_{i,3}$ & $B_{i,4}$ \\
 \hline
1  & (1,0) & (2,0)  & (3,0) & (4,0)  \\

2  & (1,2-$\sqrt{2}$) & (2, 4-2$\sqrt{2}$) & (3,6-3$\sqrt{2}$) & (4,8-4$\sqrt{2}$)\\

3  & (2-$\sqrt{2}$,1) & (4-2$\sqrt{2}$,2)  & (6-3$\sqrt{2}$,3) & (8-4$\sqrt{2}$, 4)\\

4  & (0,1) & (0,2) & (0,3) & (0,4)  \\
 \hline
\end{tabular}
\end{table}



\begin{table}[!htb]
\caption{Weights for the quarter of annulus}
\label{tbl:two-dim-weight-ctrlpnt} \centering
\begin{tabular}{c  c  c  c  c}
 \hline
i & $\omega_{i,1}$ & $\omega_{i,2}$ & $\omega_{i,3}$ & $\omega_{i,4}$ \\
 \hline
1  & 1 & 1  & 1 & 1  \\
2  & $\frac{1+\sqrt{2}}{3}$ & $\frac{1+\sqrt{2}}{3}$ & $\frac{1+\sqrt{2}}{3}$ & $\frac{1+\sqrt{2}}{3}$ \\
3  & $\frac{1+\sqrt{2}}{3}$ & $\frac{1+\sqrt{2}}{3}$  & $\frac{1+\sqrt{2}}{3}$ & $\frac{1+\sqrt{2}}{3}$ \\
4  & 1 & 1 & 1 & 1  \\
 \hline
\end{tabular}
\end{table}



\subsection*{A3: NURBS representation of the physical domain in
Example III}

 The physical domain in Example III is represented by a cubic
    B-spline solid.
 Its control points are listed in Table~\ref{tbl:three-dim-ctrlpnt}.
 The knot vectors along $u-$, $v-$, and $w-$directions are, respectively,
\begin{equation*}
    \begin{split}
    0\ 0\ 0\ 0\ 1\ 1\ 1\ 1, \\
    0\ 0\ 0\ 0\ 1\ 1\ 1\ 1, \\
    0\ 0\ 0\ 0\ 1\ 1\ 1\ 1. \\
    \end{split}
\end{equation*}


\begin{table}[!htb]

\caption{Control points of the cubic B-spline solid in Example III}
\label{tbl:three-dim-ctrlpnt} \centering
\begin{tabular}{c  c c  c  c  c}
 \hline
$i$ & $j$ & $B_{ij,1}$ & $B_{ij,2}$ & $B_{ij,3}$ & $B_{ij,4}$ \\
 \hline
1  & 1 & (0,0,0) & (0,0,1/3)  & (0,0,2/3) & (0,0,1)  \\

1  & 2 & (0,1/3,0) & (0,1/3,1/3) & (0,1/3,2/3) & (0,1/3,1)\\

1  & 3 & (0,2/3,0) & (0,2/3,1/3) & (0,2/3,2/3) & (0,2/3,1)\\

1  & 4 & (0,1,0) & (0,1,1/3) & (0,1,2/3) & (0,1,1)\\
2  & 1 & (1/3,0,0) & (1/3,0,1/3)  & (1/3,0,2/3) & (1/3,0,1)  \\

2  & 2 & (1/3,1/3,0) & (1/3,1/3,1/3) & (1/3,1/3,2/3) & (1/3,1/3,1)\\

2  & 3 & (1/3,2/3,0) & (1/3,2/3,1/3) & (1/3,2/3,2/3) & (1/3,2/3,1)\\

2  & 4 & (1/3,1,0) & (1/3,1,1/3) & (1/3,1,2/3) & (1/3,1,1)\\
3  & 1 & (2/3,0,0) & (2/3,0,1/3)  & (2/3,0,2/3) & (2/3,0,1)  \\

3  & 2 & (2/3,1/3,0) & (2/3,1/3,1/3) & (2/3,1/3,2/3) & (2/3,1/3,1)\\

3  & 3 & (2/3,2/3,0) & (2/3,2/3,1/3) & (2/3,2/3,2/3) & (2/3,2/3,1)\\

3  & 4 & (2/3,1,0) & (2/3,1,1/3) & (2/3,1,2/3) & (2/3,1,1)\\
4  & 1 & (1,0,0) & (1,0,1/3)  & (1,0,2/3) & (1,0,1)  \\

4  & 2 & (1,1/3,0) & (1,1/3,1/3) & (1,1/3,2/3) & (1,1/3,1)\\

4  & 3 & (1,2/3,0) & (1,2/3,1/3) & (1,2/3,2/3) & (1,2/3,1)\\

4  & 4 & (1,1,0) & (1,1,1/3) & (1,1,2/3) & (1,1,1)\\
 \hline

\end{tabular}

\end{table}


\subsection*{A4: NURBS representation of the physical domain in
Example IV}

 The physical domain in Example IV is represented by a cubic B-spline
    patch.
 Its control points are listed in Table~\ref{tbl:ela-ctrlpnt}.
 The knot vectors along $u-$, and $v-$directions are, respectively,
\begin{equation*}
    \begin{split}
    0\ 0\ 0\ 0\ 1\ 1\ 1\ 1, \\
    0\ 0\ 0\ 0\ 1\ 1\ 1\ 1. \\
    \end{split}
\end{equation*}


\begin{table}[!htb]
\caption{Control points of the cubic B-spline patch in Example IV.}
\label{tbl:ela-ctrlpnt} \centering
\begin{tabular}{c  c  c  c  c}
 \hline
$i$ & $B_{i,1}$ & $B_{i,2}$ & $B_{i,3}$ & $B_{i,4}$ \\
 \hline
1  & (-5.0, -1.0) & (-1.67, -1.0)  & (1.67, -1.0) & (5.0, -1.0)  \\

2  & (-5.0, -0.34) & (-1.67, -0.34) & (1.67, -0.34) & (5.0, -0.34)\\

3  & (-5.0, 0.34) & (-1.67, 0.34)  & (1.67, 0.34) & (5.0, 0.34)\\

4  & (-5.0, 1.0) & (-1.67, 1.0) & (1.67, 1.0) & (5.0, 1.0)  \\
 \hline
\end{tabular}
\end{table}


\subsection*{A5: NURBS representation of the physical domain in
Example VI}

 The physical domain of frame corner in Example VI is represented by
    a bi-quadratic NURBS surface.
 Its control points and weights are listed in Tables~\ref{tbl:fc-ctrlpnt}
    and~\ref{tbl:fc-weight}.
 The knot vectors along $u-$, and $v-$directions are, respectively,
\begin{equation*}
    \begin{split}
    0\ 0\ 0\ 0.5\ 1\ 1\ 1, \\
    0\ 0\ 0\ 1\ 1\ 1. \\
    \end{split}
\end{equation*}

\begin{table}[!htb]
\caption{Control points of the physical domain of frame corner.}
\label{tbl:fc-ctrlpnt} \centering
\begin{tabular}{c  c  c  c  c}
 \hline
$i$ & $B_{i,1}$ & $B_{i,2}$ & $B_{i,3}$ & $B_{i,4}$ \\
 \hline
1  & (-4, 0) & (-4, -4)  & (-4, -4) & (0, -4)  \\

2  & (-2.5, 0) & (-2.5, -1.5) & (-1.5, -2.5) & (0, -2.5)\\

3  & (-1, 0) & (-1, -$\sqrt{3}$/3)  & (-$\sqrt{3}$, -1) & (0, -1)\\
 \hline
\end{tabular}
\end{table}

\begin{table}[!htb]
\caption{Weights of the physical domain of frame corner.}
\label{tbl:fc-weight}
\centering
\begin{tabular}{c  c  c  c  c}
 \hline
$i$ & $B_{i,1}$ & $B_{i,2}$ & $B_{i,3}$ & $B_{i,4}$ \\
 \hline
1  & 1 & 1  & 1 & 1  \\

2  & 1 & 1 & 1 & 1 \\

3  & 1 & $\sqrt{3}/2$  & $\sqrt{3}/2$ & 1\\
 \hline
\end{tabular}
\end{table}


\subsection*{A6: Proof to formula~\pref{eq:two_dim_error} in Theorem~\ref{thm:error_formula}}

 \textbf{Proof:} We only show the theorem in the two-dimensional
    case. The proof for the one- and three-dimensional cases is similar.

 In the two-dimensional case,
    suppose the tensor product NURBS function $T_r(u,v)$ of degree $l_u \times l_v$
    is defined on the knot sequences
    \begin{equation} \label{eq:two_dim_knot_seq}
    \begin{split}
        &\{\underbrace{u_0,u_0,\cdots,u_0}_{l_u+1},u_1,\cdots,u_{n_u-1},\underbrace{u_{n_u},u_{n_u},\cdots,u_{n_u}}_{l_u+1}\},\\
        &\{\underbrace{v_0,v_0,\cdots,v_0}_{l_v+1},v_1,\cdots,v_{n_v-1},\underbrace{v_{n_v},v_{n_v},\cdots,v_{n_v}}_{l_v+1}\}.
    \end{split}
    \end{equation}
 Then the corresponding knot grid is
 \begin{equation}\label{eq:two_dim_knot_grid}
    \mathcal{T}^h = \{[u_i,u_{i+1}] \times [v_j,v_{j+1}],\ i = 0,1,\cdots,n_u-1,\ j =
    0,1,\cdots,n_v-1\}.
 \end{equation}

 Denote
 \begin{equation*}
    R(u,v) = (\mathcal{D}T(u,v) - \mathcal{D}T_r(u,v))^2,
    (u,v) \in [u_0,u_{n_u}] \times [v_0,v_{n_v}].
 \end{equation*}
 Note that $T_r(u,v)$ is generated by
    least-squares fitting the values of $\mathcal{D}T(u,v)$ at the collocation points
    $\bm{\vartheta}_k = (\eta_k, \xi_k)$, i.e., $\mathcal{D}T(\bm{\vartheta}_k), k =
    1,2,\cdots,D$, ($D \geq n_u n_v$).
 And suppose the fitting error is
 \begin{equation} \label{eq:2d_fitting_error}
    e_h = \sum_{k=1}^D R(\bm{\vartheta}_k)
    = \sum_{k=1}^D (\mathcal{D}T(\bm{\vartheta}_k) -
    \mathcal{D}T_r(\bm{\vartheta}_k))^2.
 \end{equation}

 First, based on Lemma~\ref{lem:same_knot_intervals},
    the numerical solution $T_r(u,v)$ and the approximate differential operator $\mathcal{D}T_r(u,v)$ have
    the same knot intervals,
    \begin{equation*}
        [u_i,u_{i+1}] \times [v_j,v_{j+1}],\ i = 0,1,\cdots,n_u-1,\ j =
        0,1,\cdots,n_v-1.
    \end{equation*}
 Now, consider the error between $\mathcal{D}T(u,v)$ and
    $\mathcal{D}T_r(u,v)$ in the $L^2$ norm,
    \begin{equation*}
        \norm{\mathcal{D}T(u,v)-\mathcal{D}T_r(u,v)}_{L^2}^2 =
        \int_{v_0}^{v_{n_v}} \int_{u_0}^{u_{n_u}}
        R(u,v) du dv
        = \sum_{j=0}^{n_v-1} \sum_{i=0}^{n_u-1}
        \int_{v_j}^{v_{j+1}} \int_{u_i}^{u_{i+1}}
        R(u,v) du dv.
    \end{equation*}

 Since each knot interval $[u_i,u_{i+1}] \times [v_j,v_{j+1}]$
    contains at least one collocation point,
    we suppose $\bm{\vartheta}_d = (\eta_d,\xi_d) \in [u_i,u_{i+1}] \times
    [v_j,v_{j+1}]$.
 Using the left and right rectangle integral formula repeatedly,
    we get
    \begin{equation*}
        \begin{split}
        & \int_{v_j}^{v_{j+1}} \int_{u_i}^{u_{i+1}} R(u,v) du dv
        = \int_{v_j}^{v_{j+1}} dv \left(
         \int_{u_i}^{\eta_d} R(u,v) du
         + \int_{\eta_d}^{u_{i+1}} R(u,v) du
         \right) \\
        & = \int_{v_j}^{v_{j+1}} \left((\eta_d-u_i) R(\eta_d,v) +
            (u_{i+1}-\eta_d) R(\eta_d,v) + (\eta_d-u_i)^2
            \frac{R'_u(\mu^{(1)}_{i}{(v)},v)}{2}
            + (u_{i+1}-\eta_d)^2 \frac{R'_u(\mu^{(2)}_{i}(v),v)}{2}\right) dv
        \\
        & = (u_{i+1}-{u_i}) \int_{v_j}^{v_{j+1}}  R(\eta_d,v) dv
            + (\eta_d-u_i)^2 \int_{v_j}^{v_{j+1}}  \frac{R'_u(\mu^{(1)}_{i}(v),v)}{2}
            dv + (u_{i+1}-\eta_d)^2 \int_{v_j}^{v_{j+1}} \frac{R'_u(\mu^{(2)}_{i}(v),v)}{2} dv \\
        & = (u_{i+1}-u_i)\left(
                (\xi_d-v_j) R(\eta_d,\xi_d) + (v_{j+1}-\xi_d) R(\eta_d,\xi_d)
                + (\xi_d - v_j)^2
                \frac{R'_v(\eta_d,\omega^{(1)}_{ij})}{2} + (v_{j+1} - \xi_d)^2 \frac{R'_v(\eta_d,\omega^{(2)}_{ij})}{2}
                \right) \\
        & \qquad + (\eta_d-u_i)^2 \int_{v_j}^{v_{j+1}}  \frac{R'_u(\mu^{(1)}_{i}(v),v)}{2} dv
                + (u_{i+1}-\eta_d)^2 \int_{v_j}^{v_{j+1}}  \frac{R'_u(\mu^{(2)}_{i}(v),v)}{2} dv \\
        & = (u_{i+1}-u_i)(v_{j+1}-v_j) R(\eta_d,\xi_d)
            + (u_{i+1}-u_i) \left(
                (\xi_d - v_j)^2 \frac{R'_v(\eta_d,\omega^{(1)}_{ij})}{2}
                + (v_{j+1} - \xi_d)^2 \frac{R'_v(\eta_d,\omega^{(2)}_{ij})}{2}
                \right)\\
                & \qquad + (\eta_d-u_i)^2 \int_{v_j}^{v_{j+1}}  \frac{R'_u(\mu^{(1)}_{i}(v),v)}{2} dv
                + (u_{i+1}-\eta_d)^2 \int_{v_j}^{v_{j+1}}  \frac{R'_u(\mu^{(2)}_{i}(v),v)}{2} dv,
        \end{split}
    \end{equation*}
    where $\mu^{(1)}_{i}(v), \mu^{(2)}_{i}(v) \in (u_i,u_{i+1})$ and
    $\omega_{ij}^{(1)}, \omega_{ij}^{(2)} \in (v_j,v_{j+1})$.
 By the mean value theorem of integral,
    there exist $(\bar{\mu}^{(1)}_{ij},\bar\omega^{(1)}_{ij}),(\bar{\mu}^{(2)}_{ij},\bar\omega^{(2)}_{ij})\in [u_i,u_{i+1}]\times[v_j,v_{j+1}]$
    such that
 \begin{align*}
    & \int_{v_j}^{v_{j+1}} \frac{R'_u(\mu^{(1)}(v),v)}{2} dv
    = (v_{j+1}-v_j)
        \frac{R'_u(\bar{\mu}^{(1)}_{ij},\bar\omega^{(1)}_{ij})}{2},\ \text{and}\\
    & \int_{v_j}^{v_{j+1}} \frac{R'_u(\mu^{(2)}(v),v)}{2} dv
    = (v_{j+1}-v_j)
        \frac{R'_u(\bar{\mu}^{(2)}_{ij},\bar\omega^{(2)}_{ij})}{2},\\
 \end{align*}
 Therefore,
 \begin{equation*}
 \begin{split}
    \int_{v_j}^{v_{j+1}} \int_{u_i}^{u_{i+1}} R(u,v) du dv & =
    (u_{i+1}-u_i)(v_{j+1}-v_j) R(\eta_d,\xi_d) \\
            & + (u_{i+1}-u_i) \left(
                (\xi_d - v_j)^2 \frac{R'_v(\eta_d,\omega^{(1)}_{ij})}{2}
                + (v_{j+1} - \xi_d)^2 \frac{R'_v(\eta_d,\omega^{(2)}_{ij})}{2}
                \right) \\
            & + (v_{j+1}-v_j) \left(
                (\eta_d-u_i)^2 \frac{R'_u(\bar{\mu}^{(1)}_{ij},\bar\omega^{(1)}_{ij})}{2}
                + (u_{i+1}-\eta_d)^2\frac{R'_u(\bar{\mu}^{(2)}_{ij},\bar\omega^{(2)}_{ij})}{2}
                \right).
 \end{split}
 \end{equation*}

 Moreover, we denote $\Xi = [u_0,u_{n_u}] \times [v_0,v_{n_v}]$. It is easy to show that
 \begin{align*}
    & \min_{(u,v)\in \Xi} \abs{R'_v(u,v)} \leq
     \sum_{j=0}^{n_v-1}\sum_{i=0}^{n_u-1} \frac{(u_{i+1}-u_i)(v_{j+1}-v_j)}{(u_{n_u}-u_0)(v_{n_v}-v_0)}
    \frac{\abs{R'_v(\eta_d,\omega^{(1)}_{ij})}+\abs{R'_v(\eta_d,\omega^{(2)}_{ij})}}{2}
    \leq
    \max_{(u,v) \in \Xi} \abs{R'_v(u,v)}, \\
    & \min_{(u,v)\in \Xi} \abs{R'_u(u,v)} \leq
     \sum_{j=0}^{n_v-1}\sum_{i=0}^{n_u-1} \frac{(u_{i+1}-u_i)(v_{j+1}-v_j)}{(u_{n_u}-u_0)(v_{n_v}-v_0)}
    \frac{\abs{R'_u(\bar{\mu}^{(1)}_{ij},\bar\omega^{(1)}_{ij})}+\abs{R'_u(\bar{\mu}^{(2)}_{ij},\bar\omega^{(2)}_{ij})}}{2}
    \leq
    \max_{(u,v) \in \Xi} \abs{R'_u(u,v)}.
 \end{align*}
 And then, based on the intermediate value theorem,
    there exist $(\eta^{(1)},\xi^{(1)}) \in \Xi$ and $(\eta^{(2)},\xi^{(2)}) \in
    \Xi$,
    such that
    \begin{align*}
    & \abs{R'_v(\eta^{(1)},\xi^{(1)})} = \sum_{j=0}^{n_v-1}\sum_{i=0}^{n_u-1} \frac{(u_{i+1}-u_i)(v_{j+1}-v_j)}{(u_{n_u}-u_0)(v_{n_v}-v_0)}
    \frac{\abs{R'_v(\eta_d,\omega^{(1)}_{ij})}+\abs{R'_v(\eta_d,\omega^{(2)}_{ij})}}{2},
    \\
    & \abs{R'_u(\eta^{(2)},\xi^{(2)})} = \sum_{j=0}^{n_v-1}\sum_{i=0}^{n_u-1} \frac{(u_{i+1}-u_i)(v_{j+1}-v_j)}{(u_{n_u}-u_0)(v_{n_v}-v_0)}
    \frac{\abs{R'_u(\bar{\mu}^{(1)}_{ij},\bar\omega^{(1)}_{ij})}+\abs{R'_u(\bar{\mu}^{(2)}_{ij},\bar\omega^{(2)}_{ij})}}{2}.
    \end{align*}

 As a result,
 \begin{align*}
    &\norm{\mathcal{D}T(u,v)-\mathcal{D}T_r(u,v)}_{L^2}^2 =
        \int_{v_0}^{v_{n_v}} \int_{u_0}^{u_{n_u}}
        R(u,v) du dv
        = \sum_{j=0}^{n_v-1} \sum_{i=0}^{n_u-1}
        \int_{v_j}^{v_{j+1}} \int_{u_i}^{u_{i+1}}
        R(u,v) du dv \\
     &= \sum_{j=0}^{n_v-1} \sum_{i=0}^{n_u-1} (u_{i+1}-u_i)(v_{j+1}-v_j) R(\eta_d,\xi_d) \\
            & \qquad + \sum_{j=0}^{n_v-1} \sum_{i=0}^{n_u-1} (u_{i+1}-u_i) \left(
                (\xi_d - v_j)^2 \frac{R'_v(\eta_d,\omega^{(1)}_{ij})}{2}
                + (v_{j+1} - \xi_d)^2 \frac{R'_v(\eta_d,\omega^{(2)}_{ij})}{2}
                \right) \\
            & \qquad + \sum_{j=0}^{n_v-1} \sum_{i=0}^{n_u-1} (v_{j+1}-v_j) \left(
                (\eta_d-u_i)^2 \frac{R'_u(\bar{\mu}^{(1)}_{ij},\bar\omega^{(1)}_{ij})}{2}
                + (u_{i+1}-\eta_d)^2\frac{R'_u(\bar{\mu}^{(2)}_{ij},\bar\omega^{(2)}_{ij})}{2}
                \right) \\
    & \leq h^2 \sum_{d=1}^D R(\eta_d,\xi_d) + h \sum_{j=0}^{n_v-1}\sum_{i=0}^{n_u-1} (u_{i+1}-u_i)(v_{j+1}-v_j)
                        \frac{\abs{R'_v(\eta_d,\omega^{(1)}_{ij})}+\abs{R'_v(\eta_d,\omega^{(2)}_{ij})}}{2}\\
    & \qquad \qquad \qquad \quad + h \sum_{j=0}^{n_v-1}\sum_{i=0}^{n_u-1} (u_{i+1}-u_i)(v_{j+1}-v_j)
            \frac{\abs{R'_u(\bar{\mu}^{(1)}_{ij},\bar\omega^{(1)}_{ij})}+\abs{R'_u(\bar{\mu}^{(2)}_{ij},\bar\omega^{(2)}_{ij})}}{2}
            \\
    & = h^2 e_h + h(u_{n_u}-u_0)(v_{n_v}-v_0)
        \left(\abs{R'_v(\eta^{(1)},\xi^{(1)})} +
            \abs{R'_u(\eta^{(2)},\xi^{(2)})} \right) \\
    & < h^2 \bar{M} + h(u_{n_u}-u_0)(v_{n_v}-v_0)
        \left(\abs{R'_v(\eta^{(1)},\xi^{(1)})} +
            \abs{R'_u(\eta^{(2)},\xi^{(2)})} \right). \\
 \end{align*}


 On the other hand, since the degrees of $u$ and $v$ in $T_r(u,v)$
    are both larger than the maximum orders of the partial derivatives to $u$ and $v$
    appearing in $\mathcal{D}$ (refer to~\pref{eq:original_pde}),
    respectively,
    similar as the one-dimensional case,
    $R'_u(u,v)$ and $R'_v(u,v)$ are both continuous,
    and then bounded on $\Omega \cup \partial \Omega$, i.e.,
    \begin{equation*}
        \abs{R'_v(\eta^{(1)},\xi^{(1)})} \leq \hat{M},\ \text{and},
        \ \abs{R'_u(\eta^{(2)},\xi^{(2)})} \leq \hat{M},
    \end{equation*}
    where $\hat{M}$ is a positive constant.

 In conclusion,
 \begin{equation*}
    \norm{\mathcal{D}T(u,v)-\mathcal{D}T_r(u,v)}_{L^2}^2 \leq h^2 \bar{M} + 2h(u_{n_u}-u_0)(v_{n_v}-v_0) \hat{M},
 \end{equation*}
    and Theorem~\ref{thm:error_formula} is proved. $\Box$

\bibliographystyle{unsrt}

\bibliography{isogeometric}

\end{document}